\documentclass[11pt]{report}

\usepackage{hyperref}
\usepackage{setspace}
\usepackage{latexsym}
\usepackage{amsmath}
\usepackage{amsfonts}
\usepackage{amssymb}

\newcommand{\cstar}{\bb{C}^\times}
\newcommand{\cf}{\emph{cf.}}
\newcommand{\diag}[1]{\mbox{diag}\left\{#1\right\}}
\newcommand{\eg}{\emph{e.g.}}
\newcommand{\lc}{\emph{loc.\ cit.}}
\newcommand{\eqr}[1]{\mbox{(\ref{eq:#1})}}
\newcommand{\modz}{(mod $\bb{Z}$) }
\newcommand{\ie}{\emph{i.e.}\ }
\newcommand{\qed}{\hfill{$\Box$}\vspace{1cm}}
\newcommand{\spn}[1]{\mbox{span}\left\{#1\right\}}
\newcommand{\Hom}{\mbox{Hom}}
\newcommand{\mult}{\mbox{\emph{Mult}}}
\newcommand{\e}[1]{\textbf{e}(\textstyle#1)}

\newcommand{\mc}[1]{\mathcal{#1}}
\newcommand{\up}{\upsilon}
\newcommand{\eq}[1]{\addtocounter{equation}{#1}(\theequation)\addtocounter{equation}{-#1}}
\newcommand{\G}{\Gamma }

\newcommand{\pf}{\noindent\textbf{Proof.}\hspace{.2cm}}
\newcommand{\abcd}{\begin{pmatrix}a&b\\c&d\end{pmatrix}}
\newcommand{\mat}[4]{\left(\begin{array}{rr}#1&#2\\#3&#4\end{array}\right)}
\newcommand{\cvec}[3]{\begin{pmatrix}#1\\#2\\#3\end{pmatrix}}
\newcommand{\bb}[1]{\mathbb{#1}}
\newcommand{\sln}[2]{SL(#1,\mathbb{#2})}
\newcommand{\gln}[1]{GL(#1,\bb{C})}
\newcommand{\m}{\mc{M}}

\newtheorem{thm}{Theorem}[section]
\newtheorem{rem}[thm]{Remark}
\newtheorem{lem}[thm]{Lemma}
\newtheorem{defn}[thm]{Definition}
\newtheorem{cor}[thm]{Corollary}
\newtheorem{prop}[thm]{Proposition}

\begin{document}
\pagestyle{empty}
\pagenumbering{roman}
\addcontentsline{toc}{chapter}{Abstract}

\begin{center}\textbf{Abstract}\\ \ \\\textbf{CLASSIFICATION OF VECTOR-VALUED}\\
\textbf{MODULAR FORMS OF DIMENSION}\\\textbf{LESS THAN SIX}\\ \ \\by\\ \ \\Christopher Marks
\end{center}\medskip

\noindent The main purpose of this dissertation is to classify spaces of vector-valued modular
forms associated to irreducible, $T$-unitarizable representations of the full modular group, of
dimension less than six. Given such a representation, it is shown that the associated graded
complex linear space of vector-valued modular forms is a free module over the ring
of integral weight modular forms for the full modular group, whose rank is equal to the dimension
of the given representation. An explicit basis is computed for the module structure in each case,
and this basis is used to compute the Hilbert-Poincar\'e series associated to each graded space.
\newpage
\pagestyle{plain}
\tableofcontents

\chapter{Introduction}
\pagenumbering{arabic}

Vector-valued modular forms are a natural generalization, to higher dimension, of the classical theory
of modular forms. From the vector-valued point of view, a modular form (Section~\ref{sec:mf})
is a function $f$ which is holomorphic in the complex upper half-plane $\bb{H}$, which satisfies a
moderate growth condition (Definition~\ref{defn:mg}), and which defines a 1-dimensional invariant
subspace $V=\langle f\rangle$ under a certain right action of the group $\G=SL(2,\bb{Z})$
on the space of functions holomorphic in $\bb{H}$. The vector-valued generalization occurs by again
assuming the invariant subspace $V$ consists of moderate growth functions, but allowing $V$ to have any finite dimension. Given any (finite) spanning set $\beta$ of $V$, one obtains a holomorphic vector-valued function $F$ whose components are the members of $\beta$, and a linear representation $\rho:\G\rightarrow GL(V)$ which determines
how $F$ transforms under the right action of $\G$ on $V$. Any such $F$ is, by definition, a
vector-valued modular form (Section~\ref{sec:vvmf}).\medskip

This point of view is by no means new. Indeed, it was already pointed out in \cite[pp.\ 6-7]{S} that
one may profitably employ the notion of vector-valued modular form in estimating the growth of Fourier
coefficients of cusp forms defined on finite-index subgroups of $\G$. In a similar vein, one sees in
\cite{EZ,Sh1,Sh2} that Jacobi forms, which themselves generalize (among other things) both theta series and the classical $\wp$-function of Weierstrass, can be thought of as vector-valued modular forms of half-integer weight.
And in a completely different setting -- that of rational conformal field theory -- it has been known
for quite some time that each rational Vertex Operator Algebra $V$ (assuming some mild finiteness
conditions) gives rise to a vector-valued modular form whose components are the \emph{graded dimensions}
of the simple $V$-modules associated to $V$ (\cf\ \cite{DLM} and references therein
for an explanation of these terms).\medskip

In spite of the examples and references given above, it seems that no formal study of vector-valued
modular forms had been undertaken until quite recently, when G.\ Mason and M.\ Knopp laid out in
\cite{KM,KM1,KM2} the foundations of the theory expanded upon in this dissertation. These papers
contain, among other things, the following:

\begin{enumerate}
\item A definition of vector-valued modular form equivalent to the one used
in this dissertation.
\item Estimates for the growth of Fourier coefficients of components of vector-valued
modular forms.
\item  A proof that the grading (by weight) of each space of nonzero vector-valued modular forms
is bounded from below.
\item A proof of the existence of vector-valued Poincar\'e series, which implies that
\emph{every} finite-dimensional representation of $\G=SL(2,\bb{Z})$ has associated to it
a nonzero space of vector-valued modular forms.
\end{enumerate}

More recently, Mason has investigated (\cite{M1}) the connection between vector-valued
modular forms and ordinary differential equations in the complex domain. Explicitly, he
uses the well-known modular derivative $D$ to define \emph{Modular Linear Differential
Equations} (MLDEs); these are a certain type of ODE in the punctured unit disk, whose
coefficient functions are integral weight modular forms for $\G$.
Theorem~4.1, \emph{ibid.}, shows that the solution space of any MLDE is invariant under the
slash action of $\G$. Consequently, any spanning set for such a solution space provides a candidate for
a vector-valued modular form. In addition, Mason introduces a modular version of the well-known
Wronskian from ODE theory, and establishes (Theorem~4.3, \lc) a key result, which states
that the modular Wronskian of a vector-valued modular form $F$ is a pure power of Dedekind's eta
function if, and only if, the components of $F$ form a fundamental system of solutions of an MLDE.\medskip

Finally, in \cite{M2} Mason applies the above results to give a complete classification of
spaces of vector-valued modular forms associated to irreducible representations
$$\rho:\G\rightarrow\gln{2}$$
which are \emph{$T$-unitarizable}, \ie such that $\rho\mat{1}{1}{0}{1}$ is similar to a
unitary matrix. In this setting, the following results are established:

\begin{enumerate}
\item The $\bb{C}$-linear space of vector-valued modular forms associated to
$\rho$ is $\bb{Z}$-graded as
$$\mc{H}(\rho)=\bigoplus_{k\geq0}\mc{H}(k_0+2k,\rho),$$
and the minimal weight space $\mc{H}(k_0,\rho)$ contains a vector $F_0$ whose components
form a fundamental system of solutions of an MLDE in weight $k_0$.
\item $\mc{H}(\rho)$ is a free module of rank two over the ring $\m$ of integral weight
modular forms for $\G$, with basis $\{F_0,DF_0\}$. Consequently, $\mc{H}(\rho)$ is a
cyclic module over a certain \emph{skew polynomial ring} $\mc{R}$ (\cf\ Section~\ref{sec:spr}
below), which combines the actions of the modular derivative and $\m$ on $\mc{H}(\rho)$.
\item For each $k\geq0$, $\dim\mc{H}(k_0+2k,\rho)=\left[\frac{k}{3}\right]+1$, so that the
\emph{Hilbert-Poincar\'e series} (\cf\ Section~\ref{sec:hps} below) of $\mc{H}(\rho)$ is
$$\Psi(\rho)=\frac{t^{k_0}(1+t^2)}{(1-t^4)(1-t^6)}.$$
\end{enumerate}

The main purpose of this dissertation is to generalize the preceeding results, to arbitrary real
weight and dimension less than six. It will be seen in Chapter~\ref{ch:lessthan6} that in all
cases considered, the space $\mc{H}(\rho,\up)$ of vector-valued modular forms associated to
a given irreducible, $T$-unitarizable representation $\rho$ and arbitrary multiplier system
$\up$ (\cf\ Section~\ref{sec:ms} below) will again be a free $\m$-module of rank $d$. If $d<4$
then $\mc{H}(\rho,\up)$ is again a cyclic $\mc{R}$-module, but in dimensions 4 and 5
this no longer holds in general, and is instead a very special case. In all cases considered
we are able to exhibit an explicit basis for the $\m$-module structure of $\mc{H}(\rho,\up)$,
and using this we are able to compute the Hilbert-Poincar\'e series of $\mc{H}(\rho,\up)$.

\section{Notation}

We write $\bb{Z}$ for the rational integers, $\bb{R}$ for the real numbers, $\bb{C}$ for the complex
numbers. We denote by $\bb{H}$ the complex upper half-plane $\{x+iy\in\bb{C}\ |\ y>0\}$, and will use
$z$ to denote an arbitrary element of $\bb{H}$. Given $z\in\bb{H}$, we use $q$ to denote
the exponential $q=q(z)=e^{2\pi iz}$, and for $r\in\bb{R}$, we will denote the exponential
$e^{2\pi ir}$ by $\e{r}$. The multiplicative group of nonzero complex numbers will be written
as $\cstar$, and $\bb{S}^1$ will denote the subgroup
$\{x+iy\in\cstar\ |\ x^2+y^2=1\}=\{\e{r}\ |\ r\in\bb{R}\}$, \ie the complex unit circle. We use
the superscript ``t'', on the right, to denote the transpose of a vector, so \eg\
$F=(f_1,f_2)^t$ denotes the column vector
$$F=\begin{pmatrix}f_1\\f_2\end{pmatrix}.$$

\chapter{Background}
\section{The modular group $\G$}

The set
$$\G=\sln{2}{Z}=\left\{\abcd\mid\ a,b,c,d\in\bb{Z},\ ad-bc=1\right\}$$
forms a group under ordinary multiplication of matrices, and will be
referred to as the \textbf{modular group}. We shall make extensive use of the following
three elements of $\G$:
$$R=\mat{-1}{1}{-1}{0},\ S=\mat{0}{-1}{1}{0},\ T=\mat{1}{1}{0}{1},$$
along with the accompanying identities
\begin{equation}\label{eq:modid} R^3=I,\ S^2=-I,\ RS=T.\end{equation}
It is well-known that $\G$ is generated by any two of $R,S,T$, in fact
\begin{equation}\label{eq:gammagen}
\G=\langle\ S,T\ |\ S^4=I,\ S^2=(ST)^3\ \rangle.
\end{equation}
(For two completely different proofs of this fact, \cf\ \cite[Th.\ 1.2.4]{R},
\cite[Ch.\ VII, Th.\ 2]{Se}.)\medskip

From \eq{0} we may easily prove the important
\begin{lem} Let $\G'=\langle\ xyx^{-1}y^{-1}\ |\ x,y\in\G\ \rangle$ denote the commutator
 subgroup of $\G$. Then the quotient $\G/\G'$ is cyclic of order 12, with generator $T\G'$.
\end{lem}

\pf Let $x,y$ denote the cosets $R\G',S\G'$ respectively in $\G/\G'$. Then \eqr{modid}
and \eqr{gammagen} imply that
\begin{eqnarray*}\G/\G'&=&\langle x,y\ |\ x^3=y^4=e,\ xy=yx\ \rangle\\
&=&\{x^ay^b\ |\ 0\leq a\leq2,\ 0\leq b\leq3\}.
\end{eqnarray*}
Thus $\G/\G'$ has order 12, and a direct computation shows that the element $xy=T\G'$
generates $\G/\G'$.\qed

The previous Lemma severely restricts the nature of the \textbf{characters}
of $\G$, \ie the homomorphisms $\G\rightarrow\cstar$, as we see in the important
\begin{cor}\label{cor:chi} The group Hom$(\G,\cstar)$ is cyclic of order 12, with generator
$\chi$ satisfying
$$\chi(T)=\e{\frac{1}{12}},\ \chi(S)=\e{-\frac{1}{4}},\ \chi(R)=\e{\frac{1}{3}}.$$
\end{cor}

\pf Since the commutator quotient $\G/\G'$ is finite, Hom$(\G,\cstar)\cong\G/\G'$, so
we need only exhibit an element of Hom$(\G,\cstar)$ of order 12. It is clear that the
character $\chi$ defined above does the job.\qed

\section{Integral-weight modular forms for $\G$}

In this section, we recall the basic facts we will need regarding integral weight modular forms
for the modular group $\G$. These facts are all standard, and can be found in any textbook
on modular forms, \eg\ \cite{L,R,Se}.\medskip

Recall that $\G$ acts from the left on $\bb{H}$ as \textbf{linear fractional transformations},
where for $\gamma=\abcd$, $z\in\bb{H}$ we have
$$(\gamma,z)\mapsto\gamma z=\frac{az+b}{cz+d}.$$
For each $k\in\bb{Z}$, define the map
$$J_k:\G\times\bb{H}\rightarrow\bb{C},$$
$$J_k\left(\abcd,z\right)=(cz+d)^k.$$
Direct computation shows that $J_k$ is a \textbf{1-cocycle}, \ie for any $\gamma,\sigma\in\G$, $z\in\bb{H}$,
we have
\begin{equation}\label{eq:1cocycle}
J_k(\gamma\sigma,z)=J_k(\gamma,\sigma z)J_k(\sigma,z).
\end{equation}
Let $\mc{H}$ denote the $\bb{C}$-linear space of holomorphic functions $f:\bb{H}\rightarrow\bb{C}$.
Each $k\in\bb{Z}$ defines a \textbf{slash operator} in weight $k$, which is the map
$$|_k:\mc{H}\times\G\rightarrow\mc{H},$$
$$(f,\gamma)\mapsto f|_k\gamma,$$
where for each $z\in\bb{H}$ we set
$$f|_k\gamma(z)=f(\gamma z)J_k(\gamma,z)^{-1}.$$
Because $J_k$ is a 1-cocycle for each $k\in\bb{Z}$, one sees that $\{|_k\}_{k\in\bb{Z}}$ is a
family of right actions of $\G$ on $\mc{H}$. For a given $k\in\bb{Z}$, it is clear that
$|_k\gamma$ defines an invertible linear operator on $\mc{H}$, for each $\gamma\in\G$. Thus $\mc{H}$ is a right
$\G$-module for the $|_k$ action of $\G$ on $\mc{H}$.
\begin{defn}\label{defn:mg}A function $f\in\mc{H}$ is of \textbf{moderate growth} if
there is an integer $N>0$ and a real number $c>0$ such that
$|f(x+iy)|\leq y^N$ whenever $y>c$.\end{defn}\qed

The set of all moderate growth functions on $\bb{H}$ clearly forms a subspace
$\mc{H}_\infty$ of $\mc{H}$.

\begin{defn}\label{defn:mform} Let $k\in\bb{Z}$. A function $f\in\mc{H}_\infty$ is a
\textbf{modular form} of weight $k$ if $f$ is invariant under
the right $\G$ action $|_k$, \ie if for every $\gamma\in\G$,
$$f|_k\gamma=f.$$
We write $\mc{M}_k$ for the $\bb{C}$-linear space of modular forms of weight $k$.
\end{defn}

If $f\in\m_k$, then slashing with the element $T\in\G$ shows that
$$f(z)=f|_kT(z)=f(z+1),$$
\ie $f$ is periodic of period 1. Since $f$ is holomorphic in $\bb{H}$, there is a Fourier development
$$f(z)=\sum_{n\in\bb{Z}}a(n)q^n,$$
valid in the punctured disk $0<|q|<1$, $q=e^{2\pi iz}$. For any fixed
$z_0\in\bb{H}$, the $n^{th}$ Fourier coefficient of $f$ is given by
$$a(n)=\int_{z_0}^{z_0+1}f(z)e^{-2\pi inz}dz$$
and since $f$ is of moderate growth, there is an $N>0$
and a $c>0$ such that the estimate
$$|a(n)|\leq y^Ne^{2\pi yn}$$
obtains whenever $y>c$. Taking the limit as
$y\rightarrow\infty$, we find that $a(n)=0$ for all
$n<0$. Thus
\begin{equation}
f(z)=\sum_{n\geq0}a(n)q^n,
\end{equation}
and $f(z)$ has the finite limit $a(0)$ as $Im(z)\rightarrow\infty$ in
any vertical strip of bounded width in $\bb{H}$. This finiteness property is commonly
referred to via phrases which abuse, to varying degrees, otherwise well-defined mathematical
notions, \eg\ ``$f$ is holomorphic at infinity.'' Similar abuse is understood in referencing the
\emph{order} of a nonzero function $f$ at infinity, which is defined from \eq{0} to be
$$\nu_\infty(f)=\min\{n\ |\ a(n)\neq0\}.$$
If $\nu_\infty(f)>0$ (or if $f\equiv0$ in $\bb{H}$) then $f$ is called a \textbf{cusp form}.
We denote by $\mc{S}_k$ the subspace of $\mc{M}_k$ consisting of weight $k$ cusp forms.

Since a nonzero $f\in\m_k$ is by assumption holomorphic
throughout $\bb{H}$, one may expand $f$ in a convergent power series at any $z\in\bb{H}$ and
obtain the nonnegative integer $\nu_z(f)$, which denotes the order of $f$ at $z$ in the usual sense
of the word. Performing a contour integral around the boundary of the standard
\emph{fundamental domain} for the action of $\G$ on $\bb{H}$, one
obtains (\cf \cite[Th.\ 2.1]{L}) the important \textbf{Valence Formula}, which says
\begin{equation}\label{eq:valform}
\nu_\infty(f)+\frac{1}{3}\nu_\omega(f)+\frac{1}{2}\nu_i(f)+\sum_{z\neq i,\omega}\nu_z(f)=\frac{k}{12},
\end{equation}
where $\omega=\e{\frac{1}{3}}$.

Note that the constant functions $f\equiv c\in\bb{C}$ are modular of weight $0$, so $\dim\mc{M}_0\geq1$.
Furthermore, we conclude from
$$f(z)=f|_kS^2=(-1)^kf(z)$$
that a nonzero modular form is necessarily  of \emph{even} weight. From the Valence Formula
one sees that $M_k=\{0\}$ for all $k<0$, and the space of \emph{all} integral weight modular
forms (for $\G$ and the trivial multiplier system, see below) is therefore
$$\mc{M}=\bigoplus_{k\geq0}\mc{M}_{2k}.$$

Since $J_k(\gamma,z)J_m(\gamma,z)=J_{k+m}(\gamma,z)$ holds for any $k,m\in\bb{Z}$, $\gamma\in\G$,
$z\in\bb{H}$, it is clear from Definition~\ref{defn:mform} that $\mc{M}$ is a graded ring with respect to
the $|_k$ action of $\G$, \ie if $f\in\mc{M}_k$, $g\in\mc{M}_m$ then $fg\in\mc{M}_{k+m}$.

For each even integer $k\geq4$, one defines the normalized \textbf{Eisenstein series} of weight $k$ to be
\begin{equation}\label{eq:eseries}
E_k(q)=1-\frac{2k}{B_k}\sum_{n\geq1}\sigma_{k-1}(n)q^n,
\end{equation}
where $B_k$ are the \emph{Bernoulli numbers}, defined via the identity
$$\frac{t}{e^t-1}=1-\frac{t}{2}+\sum_{k\geq2}\frac{B_k}{k!}t^k,$$
and
$$\sigma_k(n)=\sum_{0<d|n}d^k.$$
In particular, one has
$$E_4=1+240\sum_{k\geq1}\sigma_3(n)q^n,$$
$$E_6=1-504\sum_{k\geq1}\sigma_5(n)q^n.$$
For each even $k\geq4$, $E_k\in\mc{M}_k$, and one finds by the Valence Formula that
\begin{equation}\label{eq:dimintwt}
\dim\mc{M}_k=\left\{\begin{array}{lcr}\left[\frac{k}{12}\right]&\ &k\equiv2\pmod{12}\ \ \\ \\
\left[\frac{k}{12}\right]+1&\ &k\not\equiv2\pmod{12}\ .\end{array}\right.
\end{equation}
Explicitly, one has $\mc{M}_0=\bb{C}$, $\mc{M}_2=\{0\}$, $\mc{M}_k=\bb{C}E_k$ for
$k=4,6,8,10$, and for $k\geq12$ one has the recursive formula
\begin{equation}\label{eq:deltarecursion}
\mc{M}_k=\mc{S}_{k}\oplus\bb{C}E_k=\Delta\mc{M}_{k-12}\oplus\bb{C}E_k,
\end{equation}
where $\Delta$ denotes the weight 12 cusp form
\begin{eqnarray}\label{eq:Delta}
\Delta(q)&=&\frac{E_4^3(q)-E_6^2(q)}{1728}\nonumber\\
\nonumber\\
 &=&q\prod_{n=1}^\infty(1-q^n)^{24}\\
\nonumber\\
 &=&q+\sum_{n\geq2}\tau(n)q^n,\nonumber
\end{eqnarray}
known as the \textbf{discriminant} function. In many ways, this function embodies the
classical arithmetic theory of modular forms, as it was Ramanujan's initial enquiries into the
nature of $\Delta$ -- resulting in his conjectures (proven later by Mordell, Deligne resp.)
that $i)$ $\tau$ is a multiplicative function and $ii)$ the growth condition
$\tau(n)=O(n^{\frac{11}{2}+\epsilon})$ holds for any $\epsilon>0$
 -- which motivated the subsequent period of intense development of the subject.\medskip

One also finds from the Valence formula that $z=\omega$ (resp.\ $z=i$) is the only zero of
$E_4$ (resp.\ $E_6$) in $\bb{H}$. Using this fact and the modularity of the functions, one
proves easily that $E_4$ and $E_6$ are algebraically independent over $\bb{C}$. From this
and \eqr{deltarecursion}, \eqr{Delta}, one sees that $\mc{M}=\bb{C}[E_4,E_6]$, \ie $\m$
is isomorphic (as graded $\bb{C}$-algebra) to the polynomial algebra $\bb{C}[X,Y]$ in two
variables, via the identifications $X\leftrightarrow E_4$, $Y\leftrightarrow E_6$.

\section{Automorphy factors and multiplier systems for $\G$}\label{sec:ms}

In this section, we discuss automorphy factors, which generalize the 1-cocycle $J$ to arbitrary
real weight, and allow one to define modular forms of arbitrary real weight. All that follows can
be found in \cite[Chs.\ 3, 6]{R}, however it should be noted that, unlike Rankin, we do \emph{not}
require our automorphy factors to be defined on $PSL(2,\bb{Z})=\G/\{\pm I\}$;
thus we will have twelve multiplier systems for a given real weight, instead of the six described
in \lc\medskip

Let $\log z=\ln|z|+i\arg z$, $-\pi<\arg z\leq\pi$ denote the principal branch of the complex
logarithm function, \ie the choice of $\log$ which restricts to a holomorphic function in the
domain $\{z=re^{i\theta}\ |\ r>0,\ -\pi<\theta<\pi\}$. Define, as usual, $a^b=e^{b\log a}$ for
$a\in\cstar$, $b\in\bb{C}$. Having determined a fixed branch of $\log$, we are free to define
as before the map
$$J_k:\G\times\bb{H}\rightarrow\bb{C},$$
$$J_k\left(\abcd,z\right)=J_k\left(\abcd,z\right)=(cz+d)^k,$$
where $k$ is now allowed to be any real number. We have seen in the last section how
$J_k$ is a 1-cocycle when $k\in\bb{Z}$. On the other hand, if $k$ is not an integer it is
easy to find $\gamma,\sigma$ and $z$ such that \eqr{1cocycle} no longer holds. In order to
recover the cocycle structure in these cases, one is led to the following
\begin{defn}\label{defn:af}
Let $k\in\bb{R}$.  An \textbf{automorphy factor} of weight $k$ is a function
$$\nu:\G\times\bb{H}\rightarrow\bb{C}$$
satisfying the following conditions:\medskip

\begin{tabbing}\hspace{1cm}\=$i)$ For each $\gamma\in\G$, $\nu(\gamma,z)$ is holomorphic as a
function of $z$.\\ \\
\>$ii)$  $|\nu(\gamma,z)|=|J(\gamma,z)|^k$ for all
$(\gamma,z)\in\G\times\bb{H}$.\\ \\
\>$iii)$ $\nu$ is a $1-$cocycle.\end{tabbing}
\end{defn}

Note that $J_k(\gamma,z)$ is nonzero for any $\gamma\in\G$, $z\in\bb{H}$, so $ii)$ and
the maximum modulus principle imply that every automorphy factor of weight $k$ has the form
\begin{equation}\label{eq:afms}
\nu(\gamma,z)=\upsilon(\gamma)J_k(\gamma,z),
\end{equation}
for some function $\up:\G\rightarrow\bb{S}^1$. Any such $\up$ will be called a
\textbf{multiplier system} (for $\G$) of weight $k$. Using property $iii)$, we find that
\begin{equation}
\up(\gamma\alpha)=\up(\gamma)\up(\alpha)\sigma_k(\gamma,\alpha),
\end{equation}
where
$$\sigma_k(\gamma,\alpha)=\frac{J_k(\gamma,\alpha
z)J_k(\alpha,z)}{J_k(\gamma\alpha,z)}.$$
Note that $\sigma_k$ is independent of $\up$, so \eq{0} implies
that the ratio of two multiplier systems is a homomorphism
$\G\rightarrow\cstar$. This implies, via Corollary~\ref{cor:chi}, that \emph{if} there exists one
multiplier system $\up$ of weight $k$, then there are in fact 12, namely $\up\chi^n$, where
$0\leq n\leq11$ and $\chi$ is as in Corollary~\ref{cor:chi}.
The following shows that for each $k\in\bb{R}$, such an $\up$ does exist:
\begin{prop}\label{prop:upk}
Let $k\in\bb{R}$. There exists a multiplier system $\up_k$ for
$\G$ of weight $k$, satisfying
\begin{equation}\label{eq:upk}
\up_k(T)=\e{\frac{k}{12}},\ \up_k(S)=\e{-\frac{k}{4}},\ \up_k(ST)=\e{-\frac{k}{6}}.
\end{equation}
\end{prop}

\pf Since the discriminant function
\begin{equation}
\Delta(q)=q+\sum_{n=2}^\infty\tau(n)q^n
\end{equation}
is holomorphic in $|q|<1$ and vanishes only at
$q=0$, the function
$$D:\bb{H}\rightarrow\bb{C},$$
$$D(z)=\frac{\Delta(q)}{q}=1+\tau(2)q+\cdots$$
is holomorphic and nonzero in $\bb{H}$. By Cauchy's Theorem we may choose a well-defined branch
of $\log(D(q))$ and, after setting $\kappa=\frac{k}{12}$, obtain the function
$$D^\kappa(z)=e^{\kappa\log(D(q))},$$
which is holomorphic in $\bb{H}$. Multiplying by $q^\kappa$ yields a function
$$\phi_\kappa(z)=q^\kappa D^\kappa(q)$$
which remains holomorphic in $\bb{H}$. Our proposed automorphy factor is then
$$\nu:\G\times\bb{H}\rightarrow\bb{C}\,,$$
$$\nu(\gamma,z)=\frac{\phi_\kappa(\gamma z)}{\phi_\kappa(z)}.$$
It is clear that statement $i)$ of Definition~\ref{defn:af} holds. As for statement $ii)$, we have
$$|\nu(\gamma,z)|=\left|\frac{\phi_\kappa(\gamma z)}{\phi_\kappa(z)}\right|=
\left|\frac{q(\gamma z)^\kappa\left(\frac{\Delta(\gamma z)}{q(\gamma z)}\right)^\kappa}
{q(z)^\kappa\left(\frac{\Delta(z)}{q(z)}\right)^\kappa}\right|.$$
Using the fact that $|(ab)^r|=|a^rb^r|=|a|^r|b|^r$ for any $a,b\in\cstar$, $r\in\bb{R}$, and recalling
that $\Delta$ is modular of weight 12, we obtain
$$|\nu(\gamma,z)|=\left|\frac{\Delta(\gamma z)}{\Delta(z)}\right|^\kappa=|J(\gamma,z)|^{12\kappa},$$
so statement $ii)$ is satisfied by $\nu$, for the weight $12\kappa=k$.  To show $iii)$, we
compute:
$$\nu(\gamma\sigma,z)\phi_\kappa(z)=\phi_\kappa(\gamma\sigma
z)=\nu(\gamma,\sigma z)\phi_\kappa(\sigma z)=\nu(\gamma,\sigma
z)\nu(\sigma,z)\phi_\kappa(z).$$
Since $\phi_\kappa(z)$ is non-zero in $\bb{H}$, we find that $iii)$ holds, and $\nu$ is
an automorphy factor of weight $k$. By $\eqr{afms}$, we have
$$\nu(\gamma,z)=\up_k(\gamma)J_k(\gamma,z),$$
for some $\up_k$ of weight $k$. Note that
$$\nu(T,z)=\frac{\phi_k(Tz)}{\phi_k(z)}=\frac{\e{\kappa}\phi_k(z)}{\phi_k(z)}=\e{\kappa},$$
so $\up_k(T)=\e{\kappa}=\e{\frac{k}{12}}$.  Also, we may use the fact that $Si=i$ to deduce
$1=\nu(S,i)=\up_k(S)(i)^k$, \ie $\up_k(S)=(i)^{-k}=\e{-\frac{k}{4}}$. Similarly, the fact
that $ST$ fixes $\rho=\e{\frac{2}{3}}$ yields $\up_k(ST)=(\omega+1)^{-k}=\e{-\frac{k}{3}}$.\qed

For each $k\in\bb{R}$, we write $\mult(k)$ for the set of multiplier systems of weight $k$. By
Corollary~\ref{cor:chi}, we know
$$\mult(k)=\{\up_k\chi^n\ |\ 0\leq n\leq11\},$$
and that Hom$(\G,\cstar)$ acts transitively by (say) right multiplication on $\mult(k)$. Furthermore,
we have
\begin{lem}\label{lem:multprod}
For every $k,m\in\bb{R}$, $\up_k\up_m=\up_{k+m}\in\mult(k+m)$.
\end{lem}

\pf Let $\nu_k$, $\nu_m$ denote the automorphy factors associated with $\up_k$, $\up_m$, and
set $\nu(\gamma,z)=\nu_k(\gamma,z)\nu_m(\gamma,z)$ for each $\gamma\in\G$, $z\in\bb{H}$. It
is clear that $\nu$ is an automorphy factor in weight $k+m$, with associated multiplier system
$\up_k\up_m$. On the other hand, it is equally clear from \eqr{upk} that the identity
$$\frac{\up_k\up_m}{\up_{k+m}}(\gamma)=\frac{\up_k(\gamma)\up_m(\gamma)}{\up_{k+m}(\gamma)}=1$$
holds for $\gamma=S,T$. Since the ratio of two multiplier systems is a character of
$\G=\langle S,T\rangle$, we find that the above ratio is identically 1 on $\G$,
and the Lemma is proved.\qed

\begin{cor} If $\up\in\mult(k)$, $\mu\in\mult(m)$, then $\up\mu\in\mult(k+m)$.\end{cor}

\pf Write $\up=\up_k\chi^a$, $\mu=\up_m\chi^b$ for some integers $a,b$, and use the
previous Lemma.\qed

\begin{cor}\label{cor:equivmult} For any real numbers $k,m$ we have
$$\mult(k)=\mult(m)\ \leftrightarrow\ k\equiv m\pmod{\bb{Z}}.$$
\end{cor}

\pf Letting $\up=\up_k$, $\mu=\chi^n$ in the previous Corollary shows that if
$k\equiv m\pmod{\bb{Z}}$ then $\mult(k)=\mult(m)$. Conversely, if
$\mult(k)=\mult(m)$, then $\up_k=\up_m\chi^n$ for some integer $n$, so by Lemma~\ref{lem:multprod}
we have $\up_{k-m}=\chi^n$. Therefore
$$\up_{k-m}(T)=\e{\frac{k-m}{12}}=\e{\frac{n}{12}}=\chi^n(T),$$
so $k-m\equiv n\pmod{12\bb{Z}}$. In particular, $k-m\in\bb{Z}$.\qed

We will henceforth write \textbf{1} for the trivial multiplier system, which
satisfies $\textbf{1}(\gamma)=1$ for each $\gamma\in\G$. Clearly
$\textbf{1}=\up_0\in\mult(0)$, so for each $k\in\bb{Z}$ we have
$$\mult(k)=\Hom(\G,\cstar)=\langle\chi\rangle,$$
\ie the integral weight multiplier systems are exactly the characters of $\G$.\medskip

Following \cite{R} (up to a factor of 12), we define the \textbf{cusp parameter} of $\up$
to be the unique $m\in\bb{R}$, $0\leq m<12$, such that $\up(T)=\e{\frac{m}{12}}$.

\section{Modular forms of arbitrary real weight}\label{sec:mf}

Using the results of the last Section, we may extend the \emph{slash} action of $\G$ to arbitrary real weight,
and this allows us to define and classify modular forms for arbitrary real weight, as we now discuss.\medskip

Let $\nu$ be an automorphy factor of weight $k\in\bb{R}$, with associated multiplier system $\up$.
We define the map
$$|_k^\up:\mc{H}\times\G\rightarrow\mc{H},$$
$$(f,\gamma)\mapsto f|_k^\up\gamma,$$
where
$$f|_k^\up\gamma(z):=\nu(\gamma,z)^{-1}f(\gamma z).$$
Explicitly, if $\gamma=\abcd\in\G$, then for each $z\in\bb{H}$ we have
\begin{eqnarray*}
f|_k^\up\gamma(z)&=&\up(\gamma)^{-1}J_k(\gamma,z)^{-1}f(\gamma z)\\ \\
&=&\up(\gamma)^{-1}(cz+d)^{-k}f\left(\frac{az+b}{cz+d}\right).
\end{eqnarray*}
By property $iii)$ of Definition~\ref{defn:af}, $|_k^\up$ defines a right action of $\G$ on $\mc{H}$,
and as before a modular form is defined to be a moderate growth invariant of this action:

\begin{defn}\label{defn:realmod} Let $k\in\bb{R}$, and $\up\in\mult(k)$. A \textbf{modular form} of weight $k$
(for the multiplier system $\up$) is a function $f\in\mc{H}_\infty$ which satisfies
$$f|_k^\up\gamma=f$$
for every $\gamma\in\G$.
\end{defn}

We denote by $\mc{H}(k,\up)$ the space of weight $k$ modular forms for $\up$, but
when $k\in\bb{Z}$ and $\up=\textbf{1}$, we will continue to write $\mc{M}_k$ instead
of $\mc{H}(k,\textbf{1})$. Similarly, if $\up\in\mult(k)$, then we write
$$\mc{H}(\up)=\bigoplus_{n\in\bb{Z}}\mc{H}(k+n,\up)$$
to denote the space of modular forms for $\up$, but continue to write $\mc{M}$ in place of
$\mc{H}(\textbf{1})$.

As in the integral weight case, we may slash an $f\in\mc{H}(k,\up)$ with $T\in\G$ and obtain
$$f(z)=f|_k^\up T(z)=\up(T)^{-1}f(z+1)=\e{-\frac{m}{12}}f(z+1),$$
where $m$ denotes the cusp parameter of $\up$. This shows that the function
$\frac{f}{q^\frac{m}{12}}$ is periodic of period 1, and as before the moderate growth
condition implies a left-finite Fourier development
\begin{equation}\label{eq:qexpn}
f(z)=q^\lambda\sum_{n\geq0}a(n)q^n,
\end{equation}
for some $\lambda\geq0$ which satisfies $\lambda\equiv\frac{m}{12}\pmod{\bb{Z}}$. If $\lambda>0$, we
again call $f$ a \textbf{cusp form}, and write $f\in\mc{S}(k,\up)$.\medskip

Note that the proof of Proposition~\ref{prop:upk} shows that the function
$\phi_\frac{k}{12}\in\mc{H}\left(k,\up_k\right)$. In particular, for $k=\frac{1}{2}$ we have
$$\phi_\frac{1}{24}(z)=q^\frac{1}{24}\left(\frac{\Delta(q)}{q}\right)^\frac{1}{24},$$
a function which plays a decisive role in all that follows. Using \eqr{Delta}, we define
\textbf{Dedekind's eta function} $\eta(z)$ to be the scalar multiple of $\phi_\frac{1}{24}$
whose $q$-expansion begins $q^\frac{1}{24}+\cdots$, \ie whose product expansion is
\begin{equation}\label{eq:etaprod}
\eta(z)=q^\frac{1}{24}\prod_{n\geq1}(1-q^n).
\end{equation}
We have that $\eta\in\mc{H}\big(\frac{1}{2},\up_\frac{1}{2}\big)$ and, more generally, for each
$k\geq0$, $\eta^{2k}\in\mc{H}(k,\up_k)$.  By the Weierstra\ss\ theory of infinite products (\cf, \eg,
\cite[Ch.\ 5, Sec.\ 2.2]{A}) there is an $\alpha(k)\in\bb{S}^1$ such that
\begin{equation}\label{eq:eta2k}
\eta^{2k}(z)=\alpha(k)q^\frac{k}{12}\prod_{n\geq1}(1-q^n)^{2k}=\alpha(k)q^{\frac{k}{12}}+\cdots.
\end{equation}
Note that $\alpha(k)=1$ if and only if $k\in\bb{Z}$.
\begin{prop}
Suppose $\up\in\mult(k)$ has cusp parameter $m$. Then
$$\mc{H}(k,\up)=\mc{M}_{k-m}\eta^{2m}.$$
\end{prop}
\pf  Writing $\up=\up_k\chi^N$ for some $N\in\bb{Z}$, we find that the cusp parameter satisfies
the relation $k+N=12x+m$ for a unique integer $x$, so in particular $k-m$ is an integer.
By Lemma~\ref{lem:multprod}, we obtain
$$\up=\up_{k+N}=\up_m\chi^{12x}=\up_m,$$
since $\chi^{12}=\textbf{1}$. Thus for any $g\in\mc{M}_{k-m}$, we have
\begin{eqnarray*}\hspace{2.7cm}
(\eta^{2m}g)|_k^\up\gamma(z)&=&\eta^{2m}(\gamma z)g(\gamma z)\up(\gamma)^{-1}J_k(\gamma,z)^{-1}\\ \\
&=&\eta^{2m}(\gamma z)\up_m(\gamma)^{-1}J_m(\gamma,z)^{-1}g(\gamma z)J_{k-m}(\gamma,z)^{-1}\\ \\
&=&\big(\eta^{2m}|_m^\up\gamma(z)\big)\big(g|_{k-m}\gamma(z)\big)\\ \\
&=&\eta^{2m}(z)g(z),
\end{eqnarray*}
so $\mc{M}_{k-m}\eta^{2m}\subseteq\mc{H}(k,\up)$.
Conversely, for each $f\in\mc{H}(k,\up)$ we have
\begin{eqnarray*}\frac{f}{\eta^{2m}}|_{k-m}\gamma(z)&=
&\frac{f(\gamma z)J_k(\gamma,z)^{-1}}{\eta^{2m}(\gamma z)J_m(\gamma,z)^{-1}}\\ \\
&=&\frac{f(\gamma z)J_k(\gamma,z)^{-1}\up(\gamma)^{-1}}{\eta^{2m}(\gamma z)J_m(\gamma,z)^{-1}
\up_m(\gamma)^{-1}}\\ \\
&=&\frac{f|_k^\up\gamma(z)}{\eta^{2m}|_m^{\up_m}\gamma(z)}\\ \\
&=&\frac{f(z)}{\eta^{2m}(z)}.
\end{eqnarray*}
Using \eqr{qexpn}, \eqr{eta2k}, we find that the quotient has a $q$-expansion of the form
$$\frac{f(q)}{\eta^{2m}(q)}=\frac{q^\lambda\sum_{n\geq0}a(n)q^n}{\alpha(m)q^\frac{m}{12}+\cdots}=
q^{\lambda-\frac{m}{12}}\sum_{n\geq0}b(n)q^n$$ with $\lambda-\frac{m}{12}$ some integer.
Since $f$ is holomorphic, we have $\lambda\geq0$, and by definition
$0\leq m<12$. Therefore $\lambda-\frac{m}{12}>-1$ and is integral, \ie $\lambda-\frac{m}{12}\geq0$
and $\frac{f}{\eta^{2m}}$ is holomorphic at infinity.\qed

The previous Proposition characterizes $\mc{H}(\up)$ for a given multiplier system $\up$, and we
record this as
\begin{cor}\label{cor:modmult}
Let $\up$ be a multiplier system for $\G$, with cusp parameter $m$. Then
\begin{eqnarray*}
\mc{H}(\up)&=&\bigoplus_{k\geq0}\mc{H}(m+2k,\up)\\ \\
&=&\bigoplus_{k\geq0}\mc{M}_{2k}\eta^{2m}\\ \\
&=&\mc{M}\eta^{2m}.
\end{eqnarray*}\qed
\end{cor}

\section{Vector-valued modular forms}\label{sec:vvmf}

Suppose that $\rho\in\Hom(\G,\cstar)$ is a character of $\G$, $\up\in\mult(k)$, and $f\in\mc{H}(k,\rho\,\up)$
is a modular form of weight $k$ for the multiplier system $\rho\,\up$. Then by definition, $f$ is of
moderate growth, and for each $\gamma\in\G$, we have
$$f|_k^{\rho\,\up}\gamma=f.$$
Put another way, there is a 1-dimensional subspace $V=\langle f\rangle$ of $\mc{H}_\infty$,
which is a right $\G$-module under the $|_k^\up$ action of $\G$ on $\mc{H}$; using $\{f\}$
as a basis for $V$ produces the 1-dimensional representation $\rho:\G\rightarrow\gln{1}$, arising from
the $|_k^\up$ action of $\G$ on $V$, so that
$$f|_k^\up=\rho(\gamma)f$$
for each $\gamma\in\G$. To generalize the concept of modular form, one merely assumes that
the dimension of $V$ (and $\rho$) is allowed to be greater than 1:

\begin{defn}\label{defn:vvmf}Let $\rho:\G\rightarrow\gln{d}$ be a representation of $\G$ of
dimension $d$, $k$ a real number, and $\up\in\mult(k)$. A tuple of functions
$$F=\cvec{f_1}{\vdots}{f_d}$$
is a \textbf{vector-valued modular form} of weight $k$ (for $\rho$ and $\up$), if the following
conditions are satisfied:\medskip

\begin{tabbing}\hspace{3.7cm}\=\ i) $f_j\in\mc{H}_\infty$ for all $1\leq j\leq d$.\\ \\
\>ii) For each $\gamma\in\G$, $F|_k^\up\gamma=\rho(\gamma)F$,\\ \\
\end{tabbing}
where the $|_k^\up$ operation is defined componentwise on $F$. The $\bb{C}$-linear space of
vector-valued modular forms of weight $k$ for $\rho$ and $\up$ will be
denoted $\mc{H}(k,\rho,\up)$, and in light of Corollary~\ref{cor:equivmult}, we define
$$\mc{H}(\rho,\up)=\bigoplus_{n\in\bb{Z}}\mc{H}(k+n,\rho,\up).$$
\end{defn}

By definition, the components of a vector-valued modular form span a finite-dimensional
subspace of $\mc{H}_\infty$ which is a right $\G$-module for the $|_k^\up$ action of
$\G$ on $\mc{H}$. Conversely, we have

\begin{prop} Let $k\in\bb{R}$, $\up\in\mult(k)$, and suppose that $V$ is a finite-dimensional
subspace of $\mc{H}_\infty$ which is a right $\G$-module for the $|_k^\up$ action of $\G$ on $\mc{H}$.
If $\{f_1,\cdots,f_n\}$ is a spanning set for $V$, then there is a representation
$\rho:\G\rightarrow\gln{n}$ such that $F=(f_1,\cdots,f_n)^t$ is a vector-valued modular form
of weight $k$ for $\rho$ and $\up$.
\end{prop}

\pf Set $d:=\dim V$ and assume for now that $d=n$, so that the $f_i$ are linearly independent.
Then for each $i$, and for each $\gamma\in\G$, there are unique scalars $\rho_{i,j}(\gamma)$
such that
\begin{equation}f_i|_k^\up\gamma=\sum_{j=1}^n\rho_{i,j}(\gamma)f_j,\end{equation}
and $\rho(\gamma):=\big(\rho_{i,j}(\gamma)\big)$ belongs to $\gln{n}$. Let $F$ be as in the
statement of the Proposition. Then condition $i)$ of Definition~\ref{defn:vvmf} holds since $V\leq\mc{H}_\infty$,
and for each $\gamma\in\G$ we have
$$F|_k^\up\gamma=\rho(\gamma)F,$$
so condition $ii)$ holds as well. Thus we merely need to check that $\rho:\G\rightarrow\gln{n}$
is a homomorphism. This follows directly from the fact that $|_k^\up$ is a group action, \ie if
$\gamma_1,\gamma_2\in\G$, then we have
\begin{eqnarray*}0&=&F|_k^\up\gamma_1\gamma_2-F|_k^\up\gamma_1|_k^\up\gamma_2\\ \\
&=&\rho(\gamma_1\gamma_2)F-\rho(\gamma_1)\rho(\gamma_2)F\\ \\
&=&\big(\rho(\gamma_1\gamma_2)-\rho(\gamma_1)\rho(\gamma_2)\big)F.
\end{eqnarray*}
Since the components of $F$ are linearly independent, we conclude that
$$\rho(\gamma_1\gamma_2)=\rho(\gamma_1)\rho(\gamma_2),$$
as required.\smallskip

Now we may tackle the general case, wherein $\dim V=d=n-r$ for some $r\geq0$. Since the $f_i$ span
$V$, we may relabel if necessary, and assume that $\{f_1\cdots,f_d\}$ form a basis of $V$. Set
$G=(f_1,\cdots,f_d)^t$. By the above argument, there is a representation $\sigma:\G\rightarrow\gln{d}$
such that
$$G|_k^\up\gamma=\sigma(\gamma)G$$
for each $\gamma\in\G$. Choose any $R\in\gln{n}$ such that
$$RF=\left(\begin{array}{c}G\\0\end{array}\right),$$
and define
$$\rho:\G\rightarrow\gln{n},$$
$$\rho(\gamma)=R^{-1}\mat{\sigma(\gamma)}{0}{0}{I_r}R,$$
$$\ $$
where $I_r$ denotes the identity matrix of $\gln{r}$. It is clear that $\rho$ is a representation
of $\G$, and for each $\gamma\in\G$ we have

\begin{eqnarray*}F|_k^\up\gamma&=&\left(R^{-1}\begin{pmatrix}G\\0\end{pmatrix}\right)|_k^\up\gamma\\ \\
&=&R^{-1}\left(\begin{pmatrix}G\\0\end{pmatrix}|_k^\up\gamma\right)\\ \\
&=&R^{-1}\begin{pmatrix}\sigma(\gamma)G\\0\end{pmatrix}\\ \\
&=&R^{-1}\begin{pmatrix}\sigma(\gamma)&0\\0&I_r\end{pmatrix}\begin{pmatrix}G\\0\end{pmatrix}\\ \\
&=&R^{-1}\begin{pmatrix}\sigma(\gamma)&0\\0&I_r\end{pmatrix}RF\\ \\
&=&\rho(\gamma)F.
\end{eqnarray*}\qed

\begin{lem}\label{lem:2zindec}Let $k\in\bb{R}$, $\up\in\mult(k)$, and assume $\rho$ is an
indecomposable representation such that $\mc{H}(k,\rho,\up)\neq\{0\}$.  Then $\rho(S^2)=\pm I$, and
$$\mc{H}(\rho,\up)=\bigoplus_{n\in\bb{Z}}\mc{H}(k+2n,\rho,\up).$$
\end{lem}

\pf The identity $S^2=-I$ shows that $S^2$ is in the center of $\G$, so $\rho(S^2)$
is in the center of $\rho(\G)$, and has eigenvalues in the set $\{1,-1\}$. But each eigenspace
of an element in the center of $\rho(\G)$ is preserved by every element of $\rho(\G)$, so if
$\rho(S^2)$ has two nontrivial eigenspaces, then $\rho$ can be decomposed into a direct sum of
sub-representations corresponding to these eigenspaces, \ie $\rho$ is \emph{not} indecomposable.
Therefore $\rho(S^2)$ must have only one eigenspace, meaning $\rho(S^2)=\epsilon I$, with $\epsilon=1$
or $-1$.

Let $n\in\bb{Z}$ and assume there is a nonzero $F$ in $\mc{H}(k+n,\rho,\up)$. Then
\begin{equation}
\epsilon F(z)=\rho(S^2)F(z)=F|_{k+n}^\up
S^2(z)=\up(S^2)^{-1}(-1)^{k+n}F(z),
\end{equation}
from which we deduce
$$\epsilon=\up(S^2)^{-1}(-1)^{k+n}.$$
But the $n=0$ case of \eq{0} holds by assumption, so we also have
$$\epsilon=\up(S^2)^{-1}(-1)^k.$$
Thus $\mc{H}(k+n,\rho,\up)=\{0\}$ whenever $n$ is odd.\qed

It is clear that each space $\mc{H}(\rho,\up)=\bigoplus_{n\in\bb{Z}}\mc{H}(m+n,\rho,\up)$
of vector-valued modular forms is a graded $\m$-module via componentwise multiplication, \ie if
$g\in\m_k$, $F=(f_1,\cdots,f_d)^t\in\mc{H}(m+n,\rho,\up)$, then
$$gF=\cvec{gf_1}{\vdots}{gf_d}\in\mc{H}(m+n+k,\rho,\up).$$
\begin{rem}
One of the main purposes of this dissertation is to flesh out the $\m$-module structure for
spaces $\mc{H}(\rho,\up)$ for irreducible, $T$-unitarizable $\rho$ of dimension less than six. It
transpires that in all cases considered (\cf\ Ch.~\ref{ch:lessthan6} below), $\mc{H}(\rho,\up)$
is a \emph{free} $\m$-module of rank equal to the dimension of $\rho$. Recently, it has been shown
(\cite[Th.\ 1]{MM}) in the integral weight, trivial
multiplier system setting that $\mc{H}(\rho,\textbf{1})$ is free of rank $d$ over $\m$ for
\emph{any} $T$-unitarizable representation $\rho:\G\rightarrow\gln{d}$ which satisfies
$\rho(S^2)=\pm I$, but the nature of the basis for an arbitrary $\rho$ remains quite mysterious.\qed
\end{rem}

A key point to bear in mind regarding vector-valued modular forms is that the main objects of
study are the $|_k^\up$-invariant subspaces $V\leq\mc{H}_\infty$, and \emph{not} the actual
column vectors $F$. Explicitly, suppose $V=\spn{f_1,\cdots,f_d}$ and $$F=\cvec{f_1}{\vdots}{f_d}\in\mc{H}(k,\rho,\up).$$
Any $U\in\gln{d}$ will transform $F$ to a vector $UF$ whose components again span $V$, and
if we define the representation
$$\rho_U:\G\rightarrow\gln{d},$$
$$\rho_U(\gamma)=U\rho(\gamma)U^{-1},$$
then we have
\begin{eqnarray*}\rho_U(\gamma)\big(UF\big)&=&U\rho(\gamma)F\\
&=&U\big(F|_k^\up\gamma\big)\\
&=&\big(UF\big)|_k^\up\gamma.\end{eqnarray*}
Thus
$$F\in\mc{H}(k,\rho,\up)\ \Longleftrightarrow\ UF\in\mc{H}(k,\rho_U,\up),$$
and in this way a single $|_k^\up$-invariant subspace $V$ gives rise to a whole family of
graded spaces of vector-valued modular forms.\smallskip

As usual, we say that $\rho$ is \textbf{equivalent} to $\rho'$ if
$\rho'=\rho_U$ for some invertible matrix $U$. From the above discussion, we find that
\begin{lem}\label{lem:vvequiv}If $\rho,\rho':\G\rightarrow\gln{d}$ are equivalent representations
of $\G$, then for each multiplier system $\up$, $\mc{H}(\rho,\up)$ and $\mc{H}(\rho',\up)$
are isomorphic as graded $\bb{C}$-linear spaces, and as graded $\m$-modules.\qed
\end{lem}

Although Definition~\ref{defn:vvmf} places no restriction on the nature of the representation $\rho$,
the ``purest'' generalization of the classical theory of modular forms would certainly include
the condition that the component functions of any form have classical $q$-expansions. (note, however (\cf\ Sec.~\ref{sec:ode} below), 
that a more general class of \emph{logarithmic} vector-valued modular forms are well worth studying, but shall not be treated in this
dissertation.) It turns out that the \textbf{$T$-unitarizable} representations, \ie those $\rho$ for which $\rho(T)$ is conjugate to 
a unitary matrix, are the ones which ensure that components of vector-valued modular forms have familiar $q$-expansions:

\begin{lem}\label{lem:tun}Suppose that $F=(f_1,\cdots,f_d)^t\in\mc{H}(k,\rho,\up)$ has linearly independent
components, and let $m$ denote the cusp parameter of $\up$. Then the following are equivalent:
\begin{tabbing}\hspace{3cm}\=\ i) \=For each \=$1\leq j\leq d$, there is a convergent $q$-expansion\\ \\
\>\>\>$f_j(q)=q^{\lambda_j}\sum_{n\geq0}a_j(n)q^n$,\ \ $\lambda_j\geq0$.\\ \\
\>ii) $\rho(T)=\diag{\e{r_1},\cdots,\e{r_d}}$, for some real numbers $r_j$.
\end{tabbing}
\end{lem}

\pf By definition, if $F\in\mc{H}(k,\rho,\up)$ then
\begin{equation}\rho(T)F(z)=F|_k^\up T(z)=\up(T)^{-1}F(z+1).\end{equation}
Suppose $i)$ holds. Then \eq{0} says
$$\rho(T)F(z)=\e{-\frac{m}{12}}\diag{\e{\lambda_1},\cdots,\e{\lambda_j}}F(z),$$
and since the components of $F$ are linearly independent, we conclude that $ii)$ holds, for
the real numbers $r_j=\lambda_j-\frac{m}{12}$.\smallskip

Conversely, suppose that $ii)$ holds. Then \eq{0} says that
$$f_j(z+1)=\e{\frac{m}{12}}\e{r_j}f_j(z)$$
for each $j$. Let $\lambda_j$ be the unique real number satisfying $0\leq\lambda_j<1$,
$\lambda_j\equiv r_j+\frac{m}{12}$ (mod $\bb{Z}$). Then for each $j$, the function
$\frac{f_j(z)}{q^{\lambda_j}}$ has period 1 and is holomorphic in $\bb{H}$.
This implies a Fourier development
$$f_j(z)=q^{\lambda_j}\sum_{n\in\bb{Z}}a_j(n)q^n,$$
valid in the punctured disk $0<|q|<1$. For any fixed
$z_0\in\bb{H}$, the $n^{th}$ Fourier coefficient of $f_j$ is given by
$$a_j(n)=\int_{z_0}^{z_0+1}f_j(z)e^{-2\pi i(n+\lambda_j)z}dz$$
and since $f_j$ is of moderate growth, there is an $N>0$
and a $c>0$ such that the estimate
$$|a_j(n)|\leq y^Ne^{2\pi y(n+\lambda_j)}$$
obtains whenever $y>c$. Taking the limit as
$y\rightarrow\infty$, we find that $a_j(n)=0$ for all $n<0$, so $i)$ holds.
\qed

\begin{rem}\label{rem:vvequiv} If $\rho:\G\rightarrow\gln{d}$ is $T$-unitarizable but
$\rho(T)$ is not diagonal, then the components of an $F\in\mc{H}(\rho,\up)$ will
be $\bb{C}$-linear combinations of $q$-expansions of the form given by the Lemma.
This can be inconvenient for computational purposes. However, we may still assume without
loss of generality that the conditions of the Lemma hold, for the following reason:\medskip

Since $\rho$ is $T$-unitarizable, there is a $U\in\gln{d}$ such that $U\rho(T)U^{-1}$
has the form given by $ii)$ above. If $\up$ is any multiplier system for $\G$, then
Lemma~\ref{lem:vvequiv} implies that the graded spaces $\mc{H}(\rho,\up)$ and $\mc{H}(\rho_U,\up)$,
are isomorphic, and the components of any $F\in\mc{H}(\rho_U,\up)$ have $q$-expansions given by
$i)$ above.\medskip

Thus when $\rho$ is $T$-unitarizable, we will frequently say ``assume without loss that $\rho(T)$ is
diagonal'', etc.\qed
\end{rem}

\begin{lem}\label{lem:12r}Suppose $\rho:\G\rightarrow\gln{d}$ is a $T$-unitarizable
representation of $\G$ such that $\rho(T)$ has the eigenvalues $\e{r_j}$ for some real numbers $r_j$,
and set $r=r_1+\cdots+r_d$. Then $12r$ is an integer, and $\det\rho=\chi^{12r}$.\end{lem}

\pf It is clear that $\det\rho$ is a character of $\G$, so by Corollary~\ref{cor:chi} we
know $\det\rho=\chi^n$ for some integer $n$. On the other hand, we have
$$\e{r}=\det\rho(T)=\chi^n(T)=\e{\frac{n}{12}},$$
so $12r\equiv n\pmod{12\bb{Z}}$. Since $\chi^{12}=\textbf{1}$, we have
$$\det\rho=\chi^n=\chi^{12r}.$$
\qed

\begin{cor}\label{cor:detrho} Assume that $F$ and $\rho$ satisfy statements $i)$ and $ii)$ of
Lemma~\ref{lem:tun}, and set $\lambda=\lambda_1+\cdots+\lambda_d$. Then $\det\rho=\chi^{12\lambda-md}$.
\end{cor}

\pf As in the proof of the Lemma, we have for each $j$ that $r_j\equiv\lambda_j-\frac{m}{12}\pmod{\bb{Z}}$.
Setting $r=r_1+\cdots+r_d$, we have $12\lambda-md\equiv12r$ (mod $12\bb{Z}$), so by Lemma~\ref{lem:12r}
we have
$$\det\rho=\chi^{12r}=\chi^{12\lambda-md}.$$
\qed

\begin{lem}\label{lem:detmat}
Let $d\geq1$, and for $1\leq j\leq d$ suppose $k_j\in\bb{R}$, $\up_j\in\mult(k_j)$, and $F_j=(f_{1,j},\cdots,f_{d,j})^t\in\mc{H}(k_j,\rho,\up_j)$. Then
$$\det(f_{i,j})\in\mc{H}(k,\det\rho,\up),$$
where $\up=\up_1\cdots\up_d$ and $k=k_1+\cdots+k_d$.
\end{lem}

\pf We have
\begin{eqnarray*}\det(f_{ij})|_k^\up\gamma&=&\det(F_1|_{k_1}^{\up_1}\gamma,\cdots,F_d|_{k_d}^{\up_d}\gamma)\\ \\
&=&\det(\rho(\gamma)F_1,\cdots,\rho(\gamma)F_d)\\ \\
&=&\det\rho(\gamma)\det(F_1,\cdots,F_d)\\ \\
&=&\det\rho(\gamma)\det(f_{i,j}),
\end{eqnarray*}
for any $\gamma\in\G$.  It is clear that $\det(f_{i,j})$ is holomorphic in $\bb{H}$
and of moderate growth, since the $f_{i,j}$ are.\qed

\section{Hilbert-Poincar\'e series}\label{sec:hps}

Suppose that $R=\bigoplus_{j\geq0}R_j$ is a $\bb{Z}$-graded commutative ring which is a
finitely generated $\bb{C}$-algebra, with $R_0=\bb{C}$. Let the weights of the generators be $k_1,\cdots,k_d$.
Suppose also that $\lambda\in\bb{R}$ and $M=\bigoplus_{k\geq0}M_{\lambda+k}$ is a $\bb{Z}$-graded,
finitely generated $R$-module, which satisfies $\dim M_{\lambda+k}<\infty$ for all $k\geq0$. The
\textbf{Hilbert-Poincar\'e series} (\cf\ \cite[Ch.\ 2]{Be}) of $M$ is the formal power series
$$\Psi(M)(t)=\sum_{k\geq0}t^{\lambda+k}\dim M_{\lambda+k}.$$
By Theorem 2.1.1, \lc, there is a polynomial $f(t)\in\bb{Z}[t]$ such that
$$\Psi(M)(t)=\frac{t^\lambda f(t)}{\prod_{j=1}^d(1-t^{k_j})},$$
where $k_1,\cdots,k_d$ denote the weights of the generators of $M$ as $R$-module.
Now specialize the above setting to the case $R=\m$, the ring of integral weight modular forms
for $\G$, and $M=\mc{H}(\rho,\up)=\bigoplus_{k\geq0}\mc{H}(\lambda+2k,\rho,\up)$ a finitely
generated $\m$-module of vector-valued modular forms for some representation $\rho$ and
multiplier system $\up$. Bearing in mind that $\m=\bb{C}[E_4,E_6]$,
we find that the Hilbert-Poincar\'e series for $\mc{H}(\rho,\up)$ has the form
$$\Psi(\rho,\up)(t)=\frac{t^\lambda f(t)}{(1-t^4)(1-t^6)}$$
for some polynomial $f(t)\in\bb{Z}[t]$. If we assume further that $\mc{H}(\rho,\up)$ is actually
free of rank $d$ over $\m$, say with basis $\{F_1,\cdots,F_d\}$ such that
$F_j\in\mc{H}(\lambda+2k_j,\rho,\up)$, then in fact
\begin{equation}\label{eq:hp}
\Psi(\rho,\up)(t)=\frac{t^\lambda\sum_{j=1}^dt^{2k_j}}{(1-t^4)(1-t^6)}.
\end{equation}
In Chapter~4, we shall compute the Hilbert-Poincar\'e series for many different $\mc{H}(\rho,\up)$,
all of which will be free of finite rank as $\m$-module; thus we will refer back to
\eq{0} quite frequently.

\chapter{Differential Equations and Vector-Valued Modular Forms}

In this chapter, we recall briefly the theory of ordinary differential equations at regular singular
points, due to Fuchs and Frobenius, and explain the connection with vector-valued modular
forms. The discussion in Section~\ref{sec:ode} follows that of \cite[Ch.\ 5,9]{H}.  Most of
the remainder of this Chapter is a (trivial) generalization of the results of \cite{M1} to arbitrary
real weight.

\section{Ordinary differential equations at regular singular points}\label{sec:ode}

Let $\Omega$ be an open, simply connected subset of $\bb{C}$ which contains 0, and consider
an $n^{th}$ order linear ODE of the form
\begin{equation}\label{eq:ode}
q^n\left(\frac{d}{dq}\right)^nf(q)+q^{n-1}A_{n-1}(q)\left(\frac{d}{dq}\right)^{n-1}f(q)+\cdots+qA_1(q)
\frac{d}{dq}f(q)+A_0(q)f(q)=0,
\end{equation}
where the coefficient functions $A_{n-k}(q)$ are holomorphic in $\Omega$, say
$$A_{n-k}(q)=\sum_{j\geq n_k}a_k(j)q^j$$
for some nonnegative integers $n_k$. If $n_k\geq k$ for every $k$, then the solutions of
\eqr{ode} form an $n$-dimensional vector space of functions which are holomorphic throughout $\Omega$.
On the other hand, if $n_k<k$ for at least one $k$, then $q=0$ is said to be a \textbf{regular
singular point} of \eqr{ode}. In this case the nature of the set of solutions is a bit more complicated
than when \eqr{ode} is free of singularities. Given $q\in\Omega-\{0\}$ and any open disk
$U_q\subseteq\Omega-\{0\}$ centered at $q$, the solution space of \eqr{ode}, viewed as a differential equation
defined on $U_q$, is free of singularities. Thus there exists an $n$-dimensional linear space $V_q$
of functions which are holomorphic in $U_q$ and are the set of solutions of \eqr{ode} in $U_q$.
If one takes an open cover of $\Omega-\{0\}$ consisting of such disks $U_q$, then the corresponding
solution spaces $V_q$, considered over all $q\in\Omega-\{0\}$, give a locally defined system
of solutions in $\Omega-\{0\}$. Note that nothing is assured regarding the nature of a given solution
at $q=0$; it could have a removable singularity, a pole, or an essential singularity at $q=0$, depending
on the nature of the coefficient functions $A_k(q)$.\medskip

In any event, to find solutions to \eq{0} one assumes that $f(q)=\sum_{j\geq0}b(j)q^{r+j}$
for some complex number $r$, and then attempts to determine the $b(j)$ recursively. Performing the
operations indicated in \eqr{ode} and gathering like powers of $q$ yields the equations

\begin{tabbing}\hspace{1cm}\=$b(0)P(r)=0$,\\ \\
\>$b(1)P(r+1)$\=$+b(0)P_1(r)=0$,\\ \\
\>\>$\vdots$\\ \\
\>$b(k)P(r+k)+b(k-1)P_1(r+(k-1))+\cdots+b(1)P_{k-1}(r+1)+b(0)P_k(r)=0$,\\ \\
\>\>\vdots\end{tabbing}

where
\begin{equation}
P(r)=r(r-1)\cdots(r-(n-1))+a_{n-1}(0)r(r-1)\cdots(r-(n-2))+\cdots+a_1(0)r+a_0(0)
\label{eq:ind}
\end{equation}
is the so-called \textbf{indicial polynomial} for \eqr{ode}, and for $j\in\bb{N}$ we define
$$P_j(r)=a_{n-1}(j)r(r-1)\cdots(r-(n-2))+a_{n-2}(j)r(r-1)\cdots(r-(n-3))+\cdots+a_1(j)r+a_0(j).$$
Taking as an initial condition that $b(0)\neq0$, one finds from the first equation that the
exponent $r$ must be an \textbf{indicial root}, \ie a solution of $P(r)=0$. Assuming this is the case,
it is clear that the other $b(j)$ will be determined recursively from this initial data,
so long as $P(r+k)\neq0$ for $k\in\bb{N}$, \ie so long as no other indicial root has the form
$r+k$, $k\in\bb{N}$. If this condition is satisfied, then each $b(j)$ is determined uniquely by
the above formulae, and some standard arguments from analysis show that the formal series $f(q)=b(0)q^r+b(1)q^{r+1}+\cdots$ does in fact converge in some punctured disk $0<|q|<R_r$ of
positive radius, and is therefore a solution to \eq{0}. If \emph{none} of the indicial roots
differ by an integer, one may repeat this same process with each of the $n$ (distinct) indicial
roots, and produce in this way a \textbf{fundamental system} of solutions of \eqr{ode},
\ie a basis $\{f_1,\cdots,f_n\}$ of the solution space of \eqr{ode} such that
$$f_j(q)=q^{r_j}+\sum_{k\geq1}a_j(k)q^{r_j+k}$$
and $r_1,\cdots,r_n$ are the indicial roots of \eqr{ode}. This fundamental system is holomorphic
in the punctured disk $0<|q|<\min\{R_{r_1},\cdots,R_{r_n}\}$, and starting from any $q_0$ in this
disk, one may obtain the local solution spaces $V_q$ as above by analytic continuation of the
fundamental system throughout $\Omega-\{0\}$.\medskip

Note that in the unlucky event that $P(r)=P(r+k)=0$ for some $k\in\bb{N}$, this recursive process will almost
surely be unsuccessful, unless what we like to call a ``Frobenius miracle'' occurs and it happens that
$$b(k-1)P_1(r+(k-1))+\cdots+b(1)P_{k-1}(r+1)+b(0)P_k(r)=0.$$
In this case, one may choose any desired value for $b(k)$, and proceed to the next step in the recursion,
unscathed as it were. So long as there is no integer $k'>k$ such that $P(r+k')=0$,
one emerges from the process with a solution $f(q)=q^r\sum_{j\geq0}b(j)q^{r+j}$, as before.\medskip

If two indicial roots are in fact congruent (mod ${\bb{Z}}$), it is almost certain there will be no Frobenius miracle (else the name would be unjustified!), and in this case the attempt at recursion terminates abruptly. One may still determine explicitly a fundamental set of solutions of \eqr{ode}, but now the general solution will involve terms of the form $q^r\log^k qf(q)$, where $f$ is holomorphic in a neighborhood of $q=0$. In what follows, we shall always
assume that the indicial roots of our differential equations are incongruent (mod $\bb{Z}$), and as such we
provide no further details regarding this part of the theory. However, it should be pointed out that
this case, as it pertains to conformal field theory, is certainly of interest (\cf\ \cite{Mil}), and leads
to a theory of so-called \emph{logarithmic} vector-valued modular forms, wherein the accompanying
representation $\rho$ no longer satisfies the condition that $\rho(T)$ is semisimple.

\section{The modular derivative}\label{sec:dk}

In this section, we prove some basic facts about the modular derivative. Most of these results are stated,
with sketchy proofs, in \cite[Sec.\ 10.5]{L}, for the trivial multiplier system only.\bigskip

Recall the Dedekind eta function $\eta\in\mc{H}(\frac{1}{2},\up_{\frac{1}{2}})$. It follows directly
from Proposition~\ref{prop:upk} and Definition~\ref{defn:realmod} that

\begin{lem}\label{lem:etatrans}
\begin{tabbing}\hspace{3cm}\=$\eta(Tz)=\e{\frac{1}{24}}\eta(z)$,\\ \\
\>$\eta(Sz)=\sqrt{\frac{z}{i}}\eta(z)$.
\end{tabbing}\qed
\end{lem}

Define $P(z)=\frac{i}{\pi}\frac{d}{dz}\log(\eta(z))=\frac{i}{\pi}\frac{\eta'(z)}{\eta(z)}$. Then
$P$ is holomorphic in $\bb{H}$ and, using the product expansion \eqr{etaprod} for $\eta$,
we find that the $q$-expansion of $P$ is of the form
\begin{equation}\label{eq:pexpn}
P(q)=-\frac{1}{12}+\cdots,
\end{equation}
so $P$ is holomorphic and nonzero at infinity. In fact $P=-\frac{1}{12}E_2$, where $E_2$ denotes
the Eisenstein series of weight 2, \ie \eqr{eseries} with $k=2$; as we shall not use this
fact, we omit the proof, but \cf\ \cite[Th.\ 5.2]{L} for (very few) details.

\begin{lem}\label{lem:pfunc} For each $\gamma=\abcd\in\G$, we have
\begin{equation}\label{eq:pfunc}
P|_2\gamma(z)=P(z)-\frac{c}{2\pi i}(cz+d)^{-1}.
\end{equation}
\end{lem}

\pf We will establish \eqr{pfunc} for $\gamma=T,S$, and obtain the general result by induction on
the length of a word $\gamma\in\G=\langle T,S\rangle$ in $T$ and $S$. If $\gamma=T$, we find from Lemma~
\ref{lem:etatrans} that
$$P(Tz)=\frac{i}{\pi}\frac{\eta'(Tz)}{\eta(Tz)}=\frac{i}{\pi}\frac{\e{\frac{1}{24}}\eta'(z)}
{\e{\frac{1}{24}}\eta(z)}=P(z),$$
so $P$ is $T$-invariant. Since $T=\abcd$ with $c=0$, $d=1$, we have \eqr{pfunc}.

Similarly, using Lemma~\ref{lem:etatrans} and the chain rule, we find that
$$\eta'(Sz)=z^2\left(\sqrt{\frac{z}{i}}\eta'(z)+\frac{\eta(z)}{2\sqrt{iz}}\right),$$
thus we have
\begin{eqnarray*}P(Sz)&=&\frac{i}{\pi}\frac{\eta'(Sz)}{\eta(Sz)}\\ \\
&=&\frac{i}{\pi}\frac{z^2\left(\sqrt{\frac{z}{i}}\eta'(z)+\frac{\eta(z)}{2\sqrt{iz}}\right)}
{\sqrt{\frac{z}{i}}\eta(z)}\\ \\
&=&\frac{iz^2}{\pi}\frac{\eta'(z)}{\eta(z)}+\frac{iz^2}{\pi}\frac{1}{\sqrt{\frac{z}{i}}2\sqrt{iz}}\\ \\
&=&z^2P(z)+\frac{iz}{2\pi}.
\end{eqnarray*}
Therefore

\begin{eqnarray*}P|_2S(z)&=&P(Sz)z^{-2}\\ \\
&=&P(z)+\frac{i}{2\pi}z^{-1},
\end{eqnarray*}
so we have \eqr{pfunc} for $\gamma=S$.

Assume inductively that \eqr{pfunc} holds for any $\gamma\in\G$ which, as a word in the
generators $S$ and $T$, has length less than $n$. Suppose $\sigma=\gamma T$ for
some $\gamma=\abcd$ of length $n-1$. Then we have
\begin{eqnarray*}
P|_2\sigma T(z)&=&(P|_2\gamma)|_2T(z)\\ \\
&=&(P|_2\gamma)(Tz)=P(Tz)-\frac{c}{2\pi i}(cTz+d)^{-1}\\ \\
&=&P(z)\-\frac{c}{2\pi i}\big(cz+(c+d)\big)^{-1},
\end{eqnarray*}
where we use the $T$-invariance of $P$ proven above.  Noting that
$$\sigma=\abcd\mat{1}{1}{0}{1}=\mat{a}{a+b}{c}{c+d},$$
we have \eqr{pfunc}. If $\sigma=\gamma S$, then we have
\begin{eqnarray*}
P|_2\sigma(z)&=&(P|_2\gamma)|_2S(z)\\ \\
&=&(P|_2\gamma)(Sz)z^{-2}\\ \\
&=&\left[P(Sz)-\frac{c}{2\pi i}\big(c(Sz)+d\big)^{-1}\right]z^{-2}\\ \\
&=&\left[P(z)z^2-\frac{z}{2\pi i}-\frac{c}{2\pi i}\left(-\frac{c}{z}+d\right)^{-1}\right]z^{-2}\\ \\
&=&P(z)-\frac{1}{2\pi i}z^{-1}-\frac{c}{2\pi i}z^{-1}(dz-c)^{-1}\\ \\
&=&P(z)-\frac{1}{2\pi i}\left[\frac{1}{z}+\frac{c}{z(dz-c)}\right]\\ \\
&=&P(z)-\frac{d}{2\pi i}(dz-c)^{-1}.
\end{eqnarray*}
Since
$$\sigma=\gamma S=\abcd\mat{0}{-1}{1}{0}=\mat{b}{-a}{d}{-c},$$

\noindent ones sees that \eqr{pfunc} obtains in this case.  Note that an analogous proof as above shows
that \eqr{pfunc} holds for $\gamma=T^{-1}$, and from this one easily verifies, as above, that
\eqr{pfunc} holds for $\sigma=\gamma T^{-1}$.  Since $S^{-1}=-S$, it is clear that \eqr{pfunc}
holds for any $\sigma=\gamma S^{-1}$ as well. Since $S$ and $T$ generate $\G$, we are done.\qed

\begin{defn} Let $k\in\bb{R}$. The \textbf{modular derivative} in weight $k$ is the linear operator $$D_k:\mc{H}\rightarrow\mc{H},$$
\begin{eqnarray*}
D_kf(z)&:=&\frac{1}{2\pi i}\frac{d}{dz}f(z)+kP(z)f(z)\\ \\
&=&q\frac{d}{dq}f(q)+kP(q)f(q).
\end{eqnarray*}
\end{defn}\medskip

\begin{lem} Let $k\in\bb{R}$, $\up\in\mult(k)$. Then for each $f\in\mc{H}(k,\up)$,
$D_kf\in\mc{H}(k+2,\up)$.\end{lem}

\pf For $\gamma=\abcd\in\G$, $z\in\bb{H}$, we have
\begin{equation}
(D_kf)|_{k+2}^\up\gamma(z)=\up(\gamma)^{-1}(cz+d)^{-(k+2)}\left(\frac{1}{2\pi i}f'(\gamma z)+
kP(\gamma z)f(\gamma z)\right).\end{equation}

The chain rule says that
\begin{eqnarray*}
(f(\gamma z))'&=&f'(\gamma z)\gamma'(z)\\ \\
&=&f'(\gamma z)(cz+d)^{-2},
\end{eqnarray*}

\noindent but $f\in\mc{H}(k,\up)$, so we also have

\begin{eqnarray*}
(f(\gamma z))'&=&(\up(\gamma)(cz+d)^kf(z))'\\ \\
&=&\up(\gamma)(cz+d)^kf'(z)+kc\up(\gamma)(cz+d)^{k-1}f(z).
\end{eqnarray*}
Therefore

$$f'(\gamma z)=\up(\gamma)(cz+d)^{k+2}f'(z)+kc\up(\gamma)(cz+d)^{k+1}f(z).$$\smallskip

\noindent Substituting this into \eq{0} and using \eqr{pfunc} (and again using the fact
that $f|_k^\up\gamma=f$) yields
\begin{eqnarray*}(D_kf)|_{k+2}^\up\gamma(z)&=&\up(\gamma)^{-1}(cz+d)^{-(k+2)}
\Big[\frac{1}{2\pi i}\Big(\up(\gamma)(cz+d)^{k+2}f'(z)+kc\up(\gamma)(cz+d)^{k+1}f(z)\Big)\\ \\
&\ &+\ k\Big(P(z)(cz+d)^2-\frac{c}{2\pi i}(cz+d)\Big)\up(\gamma)(cz+d)^kf(z)\Big]\\ \\
&=&\frac{1}{2\pi i}\Big(f'(z)+kc(cz+d)^{-1}f(z)\Big)+\Big(kP(z)-\frac{kc}{2\pi i}(cz+d)^{-1}\Big)f(z)\\ \\
&=&\frac{1}{2\pi i}f'(z)+kP(z)f(z)\\ \\
&=&(D_kf)(z).
\end{eqnarray*}\qed

The previous Lemma shows that for each multiplier system $\up$, there is a weight two graded
operator $D$ on the space $\mc{H}(\up)$ of modular forms for $\up$, which acts as $D_k$ on
$\mc{H}(k,\up)$. It is clear that $D$ is a derivation on $\mc{H}$, in the sense that for any $f,g\in\mc{H}$,
$k,m\in\bb{R}$,
$$D_{k+m}(fg)=(D_mf)g+f(D_kg).$$
Furthermore, for $\up\in\mult(k)$, D is a graded derivation with respect to the $\m$-module structure
of $\mc{H}(\up)$, \ie if $f\in\m_m$, $g\in\mc{H}(k,\up)$, then
$$D(fg)=D_{k+m}(fg)=(D_mf)g+f(D_kg)\in\mc{H}(k+m+2,\up).$$
Most importantly, $D$ is compatible with the \emph{slash} action of $\G$ on $\mc{H}$:

\begin{lem}\label{lem:dkwt2} Let $\up\in\mult(k)$, $f\in\mc{H}$. Then for each $\gamma\in\G$,
\begin{equation}
D_k(f|_k^\up\gamma)=(D_kf)|_{k+2}^\up\gamma.
\end{equation}
\end{lem}

\pf Suppose $\gamma=\abcd$. By definition, we have

\begin{tabbing}\hspace{2cm}$D_k(f|_k^\up\gamma)(z)\ $\=$=\frac{1}{2\pi i}(f|_k^\up\gamma)'(z)+kP(z)f|_k^\up\gamma(z)$
\\ \\
\>$=\frac{1}{2\pi i}\Big(\up(\gamma)^{-1}(cz+d)^{-k}f(\gamma z)\Big)'+kP(z)f|_k^\up\gamma(z)$\\ \\
\>$=\frac{\up(\gamma)^{-1}}{2\pi i}\Big((cz+d)^{-(k+2)}f'(\gamma z)-kc(cz+d)^{-(k+1)}f(\gamma z)\Big)
+kP(z)f|_k^\up\gamma(z)$.
\end{tabbing}

On the other hand, using Lemma~\ref{lem:pfunc} we obtain

\begin{tabbing}\hspace{1cm}$(D_kf)|_{k+2}^\up\gamma(z)\ $\=$=\Big(\frac{1}{2\pi i}f'+kPf\Big)|_{k+2}^\up\gamma(z)$
\\ \\
\>$=\frac{1}{2\pi i}f'|_{k+2}^\up\gamma(z)+k(Pf)|_{k+2}^\up\gamma(z)$\\ \\
\>$=\frac{\up(\gamma)^{-1}}{2\pi i}(cz+d)^{-(k+2)}f'(\gamma z)+kP|_2\gamma(z)f|_k^\up\gamma(z)$\\ \\
\>$=\frac{\up(\gamma)^{-1}}{2\pi i}(cz+d)^{-(k+2)}f'(\gamma z)+k\Big(P(z)-\frac{c}{2\pi i}(cz+d)^{-1}\Big)
f|_k^\up\gamma(z)$\\ \\
\>$=\frac{\up(\gamma)^{-1}}{2\pi i}\Big((cz+d)^{-(k+2)}f'(\gamma z)-kc(cz+d)^{-(k+1)}f(\gamma z)\Big)
+kP(z)f|_k^\up\gamma(z)$,
\end{tabbing}
so the two sides of \eq{0} are the same.
\qed

For each $k\in\bb{R}$, $n\geq1$, we define the composition
$$D_k^n=D_{k+2(n-1)}\cdots D_k.$$
A simple inductive argument generalizes the previous Lemma to
\begin{cor}\label{cor:dkn} Let $\up\in\mult(k)$, $n\geq1$, $f\in\mc{H}$. Then for each $\gamma\in\G$,
$$(D_k^nf)|_{k+2n}^\up\gamma=D_k^n(f|_k^\up\gamma).$$
\end{cor}\qed

Using Lemma~\ref{lem:dkwt2} we may extend $D$, componentwise, to a weight two graded operator
on spaces of vector-valued modular forms:
\begin{cor}\label{cor:dkwt2vv}
If $F=(f_1,\cdots,f_n)^t\in\mc{H}(k,\rho,\up)$, then $D_kF\in\mc{H}(k+2,\rho,\up)$.
\end{cor}

\pf By Lemma~\ref{lem:dkwt2}, we have

\begin{tabbing}\hspace{2cm}$(D_kF)|_{k+2}^\up\gamma\ $\=$=\cvec{D_kf_1}{\vdots}{D_kf_n}|_{k+2}^\up\gamma
=\cvec{(D_kf_1)|_{k+2}^\up\gamma}{\vdots}{(D_kf_n)|_{k+2}^\up\gamma}
=\cvec{D_k(f_1|_k^\up\gamma)}{\vdots}{D_k(f_n|_k^\up\gamma)}$\\ \\ \\
\>$=D_k\cvec{f_1|_k^\up\gamma}{\vdots}{f_n|_k^\up\gamma}=D_k\big(F|_k^\up\gamma\big)
=D_k\big(\rho(\gamma)F\big)$\\ \\ \\
\>$=\rho(\gamma)\big(D_kF\big).$
\end{tabbing}
\qed

As above in the scalar case, it is clear that if $f\in\m_k$, $F\in\mc{H}(m,\rho,\up)$, then
$$D_{k+m}\big(fF\big)=\big(D_kf\big)F+f\big(D_mF\big)\in\mc{H}(k+m+2,\rho,\up),$$
so the extended version of $D$ continues to act as a graded derivation, which is
compatible with the $\m$-module structure of $\mc{H}(\rho,\up)$.\medskip

Finally, we note that the kernel of $D_k:\mc{H}\rightarrow\mc{H}$ is generated by $\eta^{2k}$:
\begin{lem}\label{lem:dker}Suppose $k\in\bb{R}$ and $D_kf=0$ for some $f\in\mc{H}$. Then
$f=c\eta^{2k}$ for some $c\in\bb{C}$.\end{lem}

\pf Using the fact that $P=\frac{i}{\pi}\frac{\eta'}{\eta}$, we find that
\begin{eqnarray*}D_k\eta^{2k}&=&\frac{1}{2\pi i}\frac{d}{dz}\eta^{2k}+kP\eta^{2k}\\ \\
&=&\frac{2k}{2\pi i}\eta^{2k-1}\eta'+k\frac{i}{\pi}\frac{\eta'}{\eta}\eta^{2k}\\ \\
&=&\eta'\eta^{2k-1}\left(\frac{k}{\pi i}+\frac{ki}{\pi}\right)\\ \\
&=&0.
\end{eqnarray*}
By the Fuchsian theory of linear ODEs, the $1^{st}$ order equation $D_kf=0$ has
a 1-dimensional solution space, so $\ker D_k=\langle\eta^{2k}\rangle$.\qed

\begin{cor}\label{cor:dkinj}Suppose $\rho:\G\rightarrow\gln{d}$ is irreducible, and $d\geq2$. Then for
any multiplier system $\up$, $D:\mc{H}(\rho,\up)\rightarrow\mc{H}(\rho,\up)$ is injective.\end{cor}

\pf This follows directly from Lemma~\ref{lem:dker} and the fact that any $F\in\mc{H}(\rho,\up)$
must have linearly independent components, since $\rho$ is irreducible.\qed

Finally, we formally record a fact which is obvious, but extremely useful:
\begin{lem}\label{lem:dfind} Suppose $F=(f_1,\cdots,f_n)^t\in\mc{H}(k,\rho,\up)$ has
linearly independent components. Then the set $\{F,D_kF,\cdots,D_k^{n-1}F\}$ is independent
over $\m$, so that
$$\bigoplus_{j=0}^{n-1}D_k^jF$$
forms a rank $n$ free $\m$-submodule of $\mc{H}(\rho,\up)$.
\end{lem}

\pf Suppose there is a relation
\begin{equation}
M_{n-1}D_k^{n-1}F+M_{n-2}D_k^{n-2}F+\cdots+M_1D_kF+M_0F=0
\end{equation}
with $M_j\in\m$ for each $j$. Rewriting \eq{0} in terms of the ordinary derivative $\frac{d}{dq}$
(\cf\ \eqr{dkrewrite} below) yields a differential equation $L[f]=0$ of order at most $n-1$, for
which each of the $n$ linearly independent components of $F$ is a solution. This is impossible
unless $L$ is identically 0, and one sees easily that this forces $M_j=0$ for each $j$.\qed

\section{The skew polynomial ring $\mc{R}$}\label{sec:spr}

We may combine the actions of the operator $D$ and the ring $\m$ into a single \textbf{skew polynomial
ring} $\mc{R}$, as follows. Let $d$ be a formal variable, and set
$$\mc{R}=\{M_nd^n+M_{n-1}d^{n-1}+\cdots+M_1d+M_0\ |\ n\geq0,\ M_k\in\m\}.$$
We turn $\mc{R}$ into a noncommutative ring by defining addition as for the ordinary polynomial ring $\m[d]$,
and defining multiplication via the relation
$$dM:=Md+DM,$$
for any $M\in\m$. Given a space $\mc{H}(\rho,\up)$ of vector-valued modular forms, $\mc{R}$ acts
on an $F\in\mc{H}(k,\rho,\up)$ via the relations
$$d\cdot F:=D_kF\in\mc{H}(k+2,\rho,\up),$$
$$M_n\cdot F:=M_nF\in\mc{H}(k+n,\rho,\up),$$
for any $M_n\in\m_n$. It is clear that these definitions turn each space $\mc{H}(\rho,\up)$
into a graded $\mc{R}$-module, and determining the structure of a given space $\mc{H}(\rho,\up)$,
as $\mc{R}$-module, is an important step in the classification of spaces of vector-valued modular forms.

\section{Modular linear differential equations}

\begin{defn} Let $n$ be a nonnegative integer, $k\in\bb{R}$. An $n^{th}$ order
\textbf{modular linear differential equation (MLDE)} in weight $k$ is a differential equation
in the disk $|q|<1$, which has the form
\begin{equation}\label{eq:mlde}
L[f]=D_k^nf+M_4D_k^{n-2}f+\cdots+M_{2(n-1)}D_kf+M_{2n}f=0,
\end{equation}
with $M_j\in\mc{M}_j$ for each $j$.\qed
\end{defn}

\begin{thm}\label{thm:mldevvmf} An MLDE \eq{0} in weight $k$ is either free of singularities
in $|q|<1$, or has $q=0$ as regular singular point. For each $\up\in\mult(k)$,
the solution space of \eq{0} is invariant under the $|_k^\up$ action of $\G$ on $\mc{H}$.
\end{thm}

\pf Using the fact that $P$ (\cf\ \eqr{pexpn}) and each of the $M_j$ are holomorphic in
$|q|<1$, one finds after an elementary induction argument that \eqr{mlde} may be rewritten in the form
\begin{equation}\label{eq:rs}
q^n\frac{d^nf}{dq^n}+q^{n-1}g_{n-1}(q)\frac{d^{n-1}f}{dq^{n-1}}+\cdots+qg_1(q)\frac{df}{dq}+g_0(q)f=0,
\end{equation}
for some functions $g_j(q)$ which are holomorphic in $|q|<1$. This shows that, at worst, $q=0$ is
a regular singular point of \eq{-1}. (Note that $q=0$ is \emph{not} a regular singular point of \eq{0}
if, and only if, $g_{n-j}(q)$ has a zero of order at least $j$ at $q=0$, for each
$1\leq j\leq n$. In particular, one finds by induction that
$$g_{n-1}(q)=(n(k+n-1))P(q)+\frac{n(n-1)}{2},$$
so by \eqr{pexpn} $g_{n-1}(0)=0$ if and only if $k=5(n-1)$ (see also Lemma~\ref{lem:kdet} below).
So, \eg, if $k$ is not an integer multiple of 5 then certainly \eq{0} is regular singular at $q=0$.
It seems extremely likely that $q=0$ is \emph{always} a regular singular point of \eq{0}, but we
have no proof of this.)\medskip

To establish the modular invariance property, suppose that $f$ is a solution of
\eqr{mlde}, and let $\gamma\in\G$. Using Corollary~\ref{cor:dkn} and the modularity of the $M_j$, we have
\begin{tabbing}\hspace{2cm}$L\big[f|_k^\up\gamma\big]$\=$=L\big[f|_k^\up\gamma\big]-0$\\ \\
\>$=L\big[f|_k^\up\gamma\big]-L[f]|_{k+2n}^\up\gamma$\\ \\
\>$=\sum_{j=2}^nM_{2j}D_k^{n-j}\big(f|_k^\up\gamma\big)-\sum_{j=2}^n\big(M_{2j}D_k^{n-j}f\big)|_{k+2n}^\up\gamma$\\ \\
\>$=\sum_{j=2}^nM_{2j}D_k^{n-j}\big(f|_k^\up\gamma\big)-\sum_{j=2}^n\big(M_{2j}|_{2j}\gamma\big)
\big(D_k^{n-j}f\big)|_{k+2(n-j)}^\up\gamma$\\ \\
\>$=\sum_{j=2}^nM_{2j}D_k^{n-j}\big(f|_k^\up\gamma\big)-\sum_{j=2}^nM_{2j}D_k^{n-j}\big(f|_k^\up\gamma\big)$\\ \\
\>$=0$.
\end{tabbing}
Therefore $f|_k^\up\gamma$ is again a solution of \eqr{mlde}.\qed

\begin{cor}\label{cor:solnvec} Suppose the MLDE \eqr{mlde} has nonnegative indicial roots $r_1,\cdots,r_n$,
such that $r_i\not\equiv r_j$ \modz for $i\neq j$. Then \eqr{mlde} has a fundamental system of solutions
$$f_j(q)=q^{r_j}+\sum_{k\geq1}a_j(k)q^{r_j+k},$$
$1\leq j\leq n$, and for each $\up\in\mult(k)$, there is a representation $\rho:\G\rightarrow\gln{n}$,
arising from the $|_k^\up$ action of $\G$ on the solution space of \eqr{mlde},
such that $F=(f_1,\cdots,f_n)^t$ is a vector-valued modular form for $\rho$ and $\up$.
\end{cor}

\pf As discussed previously, the theory of ODEs at regular singular points implies that \eqr{mlde}
has a fundamental system of solutions of the form
\begin{equation}
f_j(q)=q^{r_j}+\sum_{k\geq1}a_j(k)q^{r_j+k}.
\end{equation}
The $f_j$ define holomorphic (since $r_j\geq0$ for each $j$) functions on $\bb{H}$,
the universal cover of $0<|q|<1$, and are clearly of moderate growth. By Theorem~\ref{thm:mldevvmf},
for each $\up\in\mult(k)$ there is a representation $\rho$ as described in the statement of the Corollary,
and $F\in\mc{H}(k,\rho,\up)$.\qed

As in \eqr{rs}, one sees that for each $n\geq1$, there are
functions $f_{n,n-j}(q)$, holomorphic in $|q|<1$, such that
\begin{equation}\label{eq:dkrewrite}
D_k^n=q^n\frac{d^n}{dq^n}+q^{n-1}f_{n,n-1}(q)\frac{d^{n-1}}{dq^{n-1}}+\cdots+
qf_{n,1}(q)\frac{d}{dq}+f_{n,0}(q).
\end{equation}
Regarding \eq{0}, we will need the following
\begin{lem}\label{lem:kdet}
Let $n\geq1$ and write $D_k^n$ as in \eq{0}.
Then
$$f_{n,n-1}(0)=\frac{n(5(n-1)-k)}{12}.$$
\end{lem}

\pf If $n=1$, then \eq{0} reads
$$D_k=q\frac{d}{dq}+kP(q)=q\frac{d}{dq}+f_{1,0}(q),$$
so $f_{1,0}(0)=kE_2(0)=-\frac{k}{12}$, and we have the base case. Assume $n\geq2$ and write
\begin{eqnarray*}D_k^n&=&D_{k+2(n-1)}D_k^{n-1}\\ \\
&=&\scriptstyle{\left(q\frac{d}{dq}+(k+2(n-1))P\right)
\left(q^{n-1}\frac{d^{n-1}}{dq^{n-1}}+g_{n-1,n-2}(q)q^{n-2}\frac{d^{n-2}}{dq^{n-2}}+\cdots+g_{n-1,0}(q)\right)}
\end{eqnarray*}
for some functions $g_{n-1,j}$ which are holomorphic in $|q|<1$.
Comparing this to \eq{0}, we obtain
$$f_{n,n-1}(q)=(n-1)+(k+2(n-1))P(q)+g_{n-1,n-2}(q),$$
so that
$$f_{n,n-1}(0)=(n-1)-\frac{k+2(n-1)}{12}+g_{n-1,n-2}(0).$$
Assuming inductively that
$$g_{n-1,n-2}(0)=\frac{(n-1)(5(n-2)-k)}{12}$$ finishes the proof.\qed

\begin{defn} An \textbf{Eisenstein operator} of order $n\geq2$ is a  differential operator of the form
\begin{equation}
L=D_k^n+\alpha_4E_4D_k^{n-2}+\cdots+\alpha_{2(n-1)}D_k+\alpha_{2n}E_{2n},
\end{equation}
for some $\alpha_j\in\bb{C}$.
\end{defn}

\begin{lem}\label{lem:eisop} Let $L$ be an Eisenstein operator \eq{0} and
\begin{equation}\label{eq:eismlde}
L[f]=0
\end{equation}
the corresponding MLDE. Then the $\alpha_j$ and the weight $k$ are
uniquely determined by the indicial roots of \eq{0}.\end{lem}

\pf If we rewrite \eqr{eismlde} and obtain \eqr{rs}, then in the notation of
\eqr{dkrewrite} we have
\begin{equation}
g_{n-j}(q)=\left\{\begin{array}{ll}f_{n,n-1}(q)&\ j=1\\ \\
f_{n,n-2}(q)+\alpha_4E_4(q)&\ j=2\\ \\
\alpha_{2j}E_{2j}(q)+f_{n,n-j}(q)+\sum_{i=2}^{j-1}\alpha_{2i}E_{2i}(q)f_{n-i,n-j}(q)
&\ 3\leq j\leq n\ .\end{array}\right.
\end{equation}
Let $r_1,\cdots,r_n$ denote the indicial roots of
\eqr{eismlde}. The corresponding indicial equation factors as
\begin{equation}
(r-r_1)(r-r_2)\cdots(r-r_n)=\sum_{j=0}^n(-1)^jS_jr^{n-j},
\end{equation}
where $S_j$ denotes the $j$th symmetric polynomial in
$r_1,\cdots,r_n$. On the other hand, if for each $i\in\{1,2,\cdots,n\}$ we
define integers $a_{i,j}$ such that
$$r(r-1)\cdots(r-(i-1))=\sum_{j=0}^ia_{i,j}r^j,$$
then we may write the indicial equation of \eqr{eismlde} as
\begin{equation}
r^n+A_{n-1}r^{n-1}+\cdots+A_1r+A_0=0,
\end{equation}
 where
\begin{equation}
A_{n-j}=g_{n-j}(0)+\sum_{i=0}^{j-1}a_{n-i,n-j}g_{n-i}(0)
\end{equation}
for $j=1,2,\cdots,n-1$. Equating coefficients in \eq{-2} and
\eq{-1}, we obtain the identity
\begin{equation}
(-1)^jS_j=A_{n-j},
\end{equation}
valid for $j=1,\cdots,n$. Taking $j=1$ in \eq{0}, we obtain
$$-(r_1+\cdots+r_n)=g_{n-1}(0)+a_{n,n-1}.$$
Using \eq{-4} and Lemma~\ref{lem:kdet}, we find that the weight $k$ of
\eq{-5} is determined uniquely by the indicial roots of
\eq{-5}.

If $j=2$, \eq{0} says
$$S_2=g_{n-2}(0)+a_{n,n-2}+a_{n-1,n-2}g_{n-2}(0),$$
so by \eq{-4} we have $\alpha_4$ as a function of $k$ and the
indicial roots. The $j=1$ version of \eq{0} then implies that $\alpha_4$ is also
determined uniquely by the $r_j$.

For arbitrary $j\geq3$, \eq{-4}, \eq{-1} and \eq{0} show
that $\alpha_{2j}$ is a function of the $r_j$ and
$k,\alpha_4,\cdots,\alpha_{2(j-1)}$. If we assume inductively that
$k$ and $\alpha_{2i}$, $2\leq i\leq j-1$ are determined uniquely by
the indicial roots of \eqr{eismlde}, then $\alpha_{2j}$ is as well.\qed

\begin{cor}\label{cor:umlde} Let $n\leq5$. For each set $\{\lambda_1,\cdots,\lambda_n\}$
of complex numbers, there is a unique $n^{th}$ order MLDE with indicial roots
$\lambda_1,\cdots,\lambda_n$.\end{cor}

\pf This follows directly from Lemma~\ref{lem:eisop} and the
classical fact that $\m_{2j}$ is spanned by $E_{2j}$ for $j=2,3,4,5$, so that \emph{every}
MLDE of order less than 6 is of the form $L[f]=0$, with $L$ an Eisenstein operator.\qed

\section{The modular Wronskian}

Suppose $F(q)=(f_1(q),\cdots,f_n(q))^t$ is a vector-valued function which is holomorphic
in some region $\Omega\subseteq\bb{C}$, and define the $n\times n$ matrix
\begin{eqnarray}\label{eq:ordwronsk}M_F(q)&=&\begin{pmatrix}F(q),F'(q),\cdots,F^{(n-1)}(q)\end{pmatrix}\nonumber
\\\nonumber \\
&=&\begin{pmatrix}f_1&\frac{d}{dq}f_1&\cdots&\frac{d}{dq}^{\scriptscriptstyle{(n-1)}}f_1\\
f_2&\frac{d}{dq}f_2&\cdots&\frac{d}{dq}^{\scriptscriptstyle{(n-1)}}f_2\\
\vdots&\vdots&\vdots&\vdots\\
f_n&\frac{d}{dq}f_n&\cdots&\frac{d}{dq}^{\scriptscriptstyle{(n-1)}}f_n\end{pmatrix}.
\end{eqnarray}
Recall from the Fuchsian theory that the \textbf{Wronskian} of $F$ is the holomorphic function
$\mc{W}(F)(q):=\det M_F(q)$. As is well-known, the Wronskian provides a way of determining whether
or not the components of $F$ are linearly independent (as members of the linear space of holomorphic
functions in $\Omega$). Indeed, $\mc{W}(F)$ is identically zero in $\Omega$ if, and only if, the
components of $F$ form a linearly dependent set (note that this result is decidedly false
if one replaces ``holomorphic'' with ``real analytic''.) Furthermore, in the case where the components
of $F$ are solutions of an $n^{th}$ order ODE \eqr{ode} with $0\subset\Omega$ as regular singular point,
there is an explicit formula for the Wronskian of $F$, due to Abel:
\begin{equation}
\mc{W}(F)(q)=\mc{W}(F)(q_0)e^{-\int_{q_0}^q\frac{A_{n-1}(q)}{q}dq}\,,
\end{equation}
where $q_0\in\Omega-\{0\}$ is arbitrary. Since the exponential never vanishes, \eq{0} implies
\begin{prop}\label{prop:wnot0} The components of $F=(f_1,\cdots,f_n)^t$ form a fundamental system of
solutions of \eqr{ode} if, and only if, the Wronskian $\mc{W}(F)$ of $F$ is nonzero throughout
$\Omega-\{0\}$.\qed
\end{prop}

As shown in \cite{M1}, the modular version of the Wronskian is an indespensible tool in the study
of vector-valued modular forms. In this Section we restate and/or reprove some key results
from \lc, for arbitrary real weight.

\begin{defn} Let $F(q)=(f_1(q),\cdots,f_n(q))^t$ be a vector-valued function which is holomorphic
in $0<|q|<1$. The \textbf{modular Wronskian} of $F$ (in weight $k\in\bb{R}$) is the function
$$W(F)(q)=\det\big(F(q),D_kF(q),\cdots,D_k^{n-1}F(q)\big).$$
\end{defn}

It is clear that $W(F)$ is again holomorphic in $0<|q|<1$, and if the $f_j$ are all holomorphic
at $q=0$, then so is $W(F)$. The first thing to notice is that the modular and ordinary
Wronskians of $F$ differ by a power of $q$:
\begin{lem}\label{lem:nonvanish} Suppose $F(q)=(f_1(q),\cdots,f_n(q))^t$ is holomorphic in $0<|q|<1$. Then  $$W(F)=q^\frac{n(n-1)}{2}\mc{W}(F).$$
\end{lem}

\pf Recalling that elementary column operations do not change the determinant of a matrix, it follows
directly from the definition of $D_k^j$ that
\begin{equation}
W(F)(q)=\det\begin{pmatrix}f_1&q\frac{d}{dq}f_1&\cdots&q^{n-1}\frac{d}{dq}^{\scriptscriptstyle{(n-1)}}f_1\\
f_2&q\frac{d}{dq}f_2&\cdots&q^{n-1}\frac{d}{dq}^{\scriptscriptstyle{(n-1)}}f_2\\
\vdots&\vdots&\vdots&\vdots\\
f_n&q\frac{d}{dq}f_n&\cdots&q^{n-1}\frac{d}{dq}^{\scriptscriptstyle{(n-1)}}f_n\end{pmatrix}.
\end{equation}
Expanding down the last column of \eq{0}, we obtain
$$W(F)=\sum_{j=1}^n(-1)^{n-1+j}q^{n-1}\left(\frac{d}{dq}^{\scriptscriptstyle{(n-1)}}f_j\right)W(F,j),$$
where $W(F,j)$ denotes the sub-determinant obtained by omitting the $j^{th}$ row and $(n-1)^{st}$
column of \eq{0}. Inductively, we may assume that $W(F,j)=q^{\frac{(n-2)(n-1)}{2}}\mc{W}(F,j)$,
where $\mc{W}(F,j)$ denotes the Wronskian obtained by omitting the $j^{th}$ row and $(n-1)^{st}$
column from \eqr{ordwronsk}. Thus
\begin{eqnarray*}W(F)&=&\sum_{j=1}^n(-1)^{n+j-1}q^{n-1}\left(\frac{d}{dq}^{\scriptscriptstyle{(n-1)}}f_j
\right)q^{\frac{(n-2)(n-1)}{2}}\mc{W}(F,j)\\ \\
&=&q^{\frac{n(n-1)}{2}}\sum_{j=1}^n(-1)^{n-1+j}\left(\frac{d}{dq}^{\scriptscriptstyle{(n-1)}}f_j
\right)\mc{W}(F,j)\\ \\
&=&q^{\frac{n(n-1)}{2}}\mc{W}(F).\end{eqnarray*}\qed

We also record
\begin{lem}[Mason]\label{lem:modwnot0} Assume that
$$F=\cvec{f_1}{\vdots}{f_n}=\cvec{q^{\lambda_1}+\sum_{j\geq0}a_1(j)q^{\lambda_1+j}}{\vdots}
{q^{\lambda_n}+\sum_{j\geq0}a_n(j)q^{\lambda_n+j}}$$
has linearly independent components, for some pairwise distinct, nonnegative real numbers
$\lambda_1,\cdots,\lambda_n$. Then the modular Wronskian of $F$ has the form
$$W(F)=q^{\lambda}\sum_{k\geq0}a(k)q^k,$$
where $\lambda=\lambda_1+\cdots+\lambda_n$ and $a(0)\neq0$.\end{lem}

\pf See the proof of \cite[Lem.\ 3.6]{M1}.\qed

Combining Lemma~\ref{lem:detmat} and Corollary~\ref{cor:dkwt2vv}, we obtain the useful
\begin{cor}\label{cor:wfmod}
If $F\in\mc{H}(k,\rho,\up)$, then
$$W(F)\in\mc{H}(d(k+d-1),\det\rho\,\up^d),$$
where $\dim\rho=d$.\qed
\end{cor}

The following Theorem is especially important, as it provides, given $\rho$ and $\up$, a lower bound for
the weight of a nonzero form in $\mc{H}(\rho,\up)$:
\begin{thm}\label{thm:wetag}Let $\rho:\G\rightarrow\gln{d}$ be a $T$-unitarizable representation,
$\up\in\mult(k)$, and assume that
$$F=\cvec{q^{\lambda_1}\sum_{n\geq0}a_1(n)q^n}{\vdots}{q^{\lambda_d}\sum_{n\geq0}a_d(n)q^n}\in\mc{H}(k,\rho,\up)$$
has linearly independent components, with $a_j(0)\neq0$ for each $j$. Then
$$W(F)=\eta^{24\lambda}g,$$
where $\lambda=\lambda_1+\cdots+\lambda_d$ and
$g\in\mc{M}_{d(k+d-1)-12\lambda}$ is not a cusp form. In particular, the weight $k$ of $F$ must satisfy
$$k\geq\frac{12\lambda}{d}+1-d.$$
\end{thm}

\pf From Lemma~\ref{lem:modwnot0}, one knows that
$$W(F)=cq^\lambda+\cdots$$
for some $c\in\cstar$. From the product expansion ~\ref{eq:eta2k}, we find that
$\eta^{24\lambda}=\alpha(12\lambda)q^\lambda+\cdots$ is nonzero in $\bb{H}$, and the
quotient
$$g(q):=\frac{W(F)}{\eta^{24\lambda}}=c(0)+\sum_{n\geq1}c(n)q^n$$ is
of moderate growth, with $c(0)=\frac{c}{\alpha(12\lambda)}\neq0$. By Corollary~\ref{cor:wfmod},
we know
$$W(F)\in\mc{H}(d(k+d-1),\det\rho\,\up^d),$$ and we recall that
$\eta^{24\lambda}\in\mc{H}(12\lambda,\up_{12\lambda})$. Therefore
$$\frac{W(F)}{\eta^{24\lambda}}=
\frac{W(F)|_{d(k+d-1)}^{\up^d\det\rho}\gamma}{\eta^{24\lambda}
|_{12\lambda}^{\up_{12\lambda}}\gamma}=
\left.\frac{W(F)}{\eta^{24\lambda}}\right|_{d(k+d-1)-12\lambda}
^{\up^d\,\up_{\scriptscriptstyle{-12\lambda}}\det\rho}\gamma,$$ for
each $\gamma\in\G$. Using Lemma~\ref{lem:multprod} and Corollary~\ref{cor:detrho}, we find that
$$\up^d\,\up_{-12\lambda}\det\rho=\up_{md}\chi^{12\lambda-md}\up_{-12\lambda}=\textbf{1},$$
and
$\frac{W(F)}{\eta^{24\lambda}}\in\mc{M}_{d(k+d-1)-12\lambda}$, as
required.\qed

\begin{thm}\label{thm:wmlde} Under the hypotheses of Theorem~\ref{thm:wetag}, the following are equivalent:
\begin{tabbing}\hspace{2cm}\=i) $g$ is a nonzero constant, \ie $W(F)=c\eta^{24\lambda}$ is a pure power of $\eta$.\\ \\
\>ii) \=The components of $F$ form a fundamental system of solutions\\
\>\>of an MLDE in weight $k$.\\ \\
\>iii) $k=\frac{12\lambda}{d}+1-d$.
\end{tabbing}
\end{thm}

\pf Let $V=\spn{f_1,\cdots,f_d}$ be the subspace of $\mc{H}_\infty$ spanned by the components of $F$.
Note that $\dim V=d$ since the $f_i$ are linearly independent. For any $f\in V$, the determinant
\begin{equation}
\left|\begin{array}{cccc}f&D_kf&\ &D_k^df\\f_1&D_kf_1&\ &D_k^df_1\\\vdots&\vdots&\cdots&\vdots\\
f_d&D_kf_d&\ &D_k^df_d\end{array}\right|
\end{equation}
vanishes identically, since each column consists of linearly dependent vectors. Expanding across
the first row, we obtain the equation
\begin{equation}\sum_{j=0}^d(-1)^jW(F,d-j)D_k^{d-j}f=0,\end{equation}
where for each $0\leq j\leq d$, $W(F,d-j)$ denotes the sub-determinant of \eq{-1} obtained by
omitting the $1^{st}$ row and $(d-j)^{th}$ column. Note that $W(F,d)=W(F)$, the modular
Wronskian of $F$, and since $F\in\mc{H}(k,\rho,\up)$, we
find from Lemma~\ref{lem:detmat} and Corollary~\ref{cor:dkwt2vv} that for each $j$,
$$W(F,d-j)\in\mc{H}(d(k+d-1)+2j,\det\rho\,\up^d).$$
Therefore
$$\frac{W(F,d-j)}{W(F,d)}=\frac{W(F,d-j)}{W(F)}$$
is a meromorphic modular form of weight $2j$ (for the trivial multiplier system \textbf{1}).
Dividing and rewriting \eq{0} yields the ODE
\begin{equation}L[f]=D_k^df+\sum_{j=1}^d(-1)^j\frac{W(F,d-j)}{W(F)}D_k^{d-j}f=0,
\end{equation}
which has $V$ as its solution space.\smallskip

Suppose $i)$ holds. Since $\eta^{24\lambda}$ does not vanish in $\bb{H}$, one sees that
each coefficient function in \eq{0} is holomorphic in $\bb{H}$. By an identical argument
as in the proof of Lemma~\ref{lem:modwnot0}, it is found that for each $j$, $W(F,d-j)$
has the form
$$W(F,d-j)=q^{\lambda}\sum_{n\geq0}b_j(n)q^n,$$
where $b_j(0)$ may, or may not, be zero. By \eqr{eta2k}, we have
$\eta^{24\lambda}=\alpha(12\lambda)q^{\lambda}+\cdots$, so we find that
$$\frac{W(F,d-j)}{W(F)}\in\m_{2j}$$
for each $j$ (note in particular that $W(F,d-1)=0$, since there are no nonzero holomorphic modular
forms of weight 2 for $\G$). Therefore \eq{0} is in fact an MLDE, with solution space $V$, and
$ii)$ holds.\smallskip

On the other hand, if $ii)$ holds, then the components of $F$ are linearly independent throughout
$0<|q|<1$, so by Corollary~\ref{cor:wfmod}, Proposition~\ref{prop:wnot0}, and Lemma~\ref{lem:nonvanish},
$W(F)$ is a modular form which does not vanish in $\bb{H}$. As is well-known, the only such forms are
pure powers of $\eta$, so $i)$ holds.\smallskip

It is clear from Theorem~\ref{thm:wetag} that $i)$ and $iii)$ are equivalent, so we are done.\qed

Using Theorem~\ref{thm:wmlde}, we obtain the following important fact concerning
representations arising from MLDEs:

\begin{cor}\label{cor:indec}
Suppose $\rho:\G\rightarrow\gln{d}$ is $T$-unitarizable, $\up\in\mult(k)$, and $\mc{H}(\rho,\up)$
contains a vector
$$F=\cvec{f_1}{\vdots}{f_d}=\cvec{q^{\lambda_1}+\sum_{n\geq1}a_1(n)q^n}{\vdots}{q^{\lambda_d}+
\sum_{n\geq1}a_d(n)q^n}$$
whose components form a fundamental system of solutions of an MLDE in weight $k$.  Then $\rho$
is indecomposable.
\end{cor}

\pf Suppose $\rho$ decomposes into a direct sum
$\rho=\rho_1\oplus\rho_2$. We may assume, up to equivalence of representation,
that the $|_k^\up$-invariant subspaces corresponding to $\rho_1$ and
$\rho_2$ are spanned by $\{f_1,\cdots,f_{d_1}\}$, $\{f_{\scriptscriptstyle{d_1+1}},\cdots,f_d\}$
respectively, for some $1\leq d_1\leq d$. Set $\dim\rho_2=d_2=d-d_1$, and
$\Lambda_1=\lambda_1+\cdots+\lambda_{d_1}$,
$\Lambda_2=\lambda_{d_1+1}+\cdots+\lambda_d$. By Theorem~\ref{thm:wmlde}, we
have
\begin{equation}
d(k+d-1)=12(\Lambda_1+\Lambda_2).
\end{equation}
On the other hand, if we define $F_1=(f_1,\cdots,f_{d_1})^t$,
$F_2=(f_{d_1+1},\cdots,f_d)^t$, then $F_1\in\mc{H}(k,\rho_1,\up)$,
$F_2\in\mc{H}(k,\rho_2,\up)$, and Theorem~\ref{thm:wetag} yields the
inequalities
$$d_1(k+d_1-1)\geq12\Lambda_1,$$
$$d_2(k+d_2-1)\geq12\Lambda_2.$$
Adding these inequalities and using $\eq{0}$ yields the
inequality
$$2d_1d_2\leq0,$$
so that $d_1=d$, $d_2=0$.\qed

\chapter{Vector-Valued Modular Forms of Dimension Less Than Six}\label{ch:lessthan6}

This chapter brings us to the main results of this dissertation, which are the classification
theorems for spaces of irreducible, $T$-unitarizable vector-valued modular forms of dimension
less than six.

\section{Preliminaries}

\begin{lem}\label{lem:div}
Let $\rho:\G\rightarrow\gln{d}$ be an irreducible representation of dimension
$d\leq4$, and assume the eigenvalues of $\rho(T)$ are $\e{r_j}$, for some
real numbers $r_1,\cdots,r_d$. Then $d\mid12(r_1+\cdots+r_d)$.\end{lem}

\pf Set $r=r_1+\cdots+r_d$.  By Lemma~\ref{lem:12r}, $12r\in\bb{Z}$, and this
establishes the Lemma if $d=1$.\medskip

Suppose $d=2$.  As in the proof of Lemma~\ref{lem:2zindec}, the irreducibility of
$\rho$ implies that $\rho(S^2)=\pm I$, so $\det\rho(S^2)=1$ since $d=2$.
The identities $RS=T$, $R^3=I$ then imply
$$\e{6r}=\det\rho(T^6)=\det\rho((RS)^6)=\det\rho(R^6)\det\rho(S^6)=1,$$
so that $6r$ is an integer, as required.\medskip

Suppose $d=3$.  The identity $R^{3}=I$ implies that $\rho(R)$ is diagonalizable,
with each eigenvalue a cube root of 1. If there is an eigenspace $U$ of $\rho(R)$
satisfying $\dim U>1$, then using $\rho(S^2)=\pm I$ and $\dim\rho=3$, we find that
$U\cap\rho(S)U$ is a nonzero subspace which is invariant under both $\rho(R)$
and $\rho(S)$. The existence of such a $U$ would therefore contradict the
irreducibility of $\rho$, since $R$ and $S$ generate $\G$. This means $\rho(R)$
has the distinct eigenvalues $\omega=\e{\frac{1}{3}},\omega^2,\omega^3$, and
$\det\rho(R)=\omega^6=1$. Since $S^4=I$ and $RS=T$, we have
$$\e{4r}=\det\rho(T^4)=\det\rho((RS)^4)=\det\rho(R)=1,$$
and $4r$ is an integer.\medskip

If $d=4$, we again use $\rho(S^2)=\pm I$, and conclude that the
eigenvalues of $\rho(S)$ are $\pm1$, $\pm i$ respectively. Note that
in either case both eigenvalues occur, since $\rho$ is irreducible.
Similar to the $d=3$ case, if $\rho(S)$ had an eigenspace $U$ of dimension 3,
then the nonzero subspace $U\cap\rho(R)U\cap\rho(R^2)U$ would be invariant under both
$\rho(R)$ and $\rho(S)$, violating the irreducibility of $\rho$. Therefore the eigenvalues
of $\rho(S)$ are either $\{1,1,-1,-1\}$ or $\{i,i,-i,-i\}$, and
either way we have $\det\rho(S)=1$. This implies
$$\e{3r}=\det\rho(T^3)=\det\rho((RS)^3)=\det\rho(S)=1,$$ so $3r$ is
an integer.\qed

We next record the content of the Corollaries following Proposition 2.2 and Theorem
2.9, respectively, in \cite{TW}:
\begin{thm}[Tuba and Wenzl]\label{thm:tw} Let $\rho:\G\rightarrow\gln{d}$ be an irreducible
representation of $\G$, $d<6$. Then the following hold:

\begin{tabbing}\=\ i) The minimal and characteristic polynomials of $\rho(T)$ coincide.\\ \\
\>ii) \=If $d\neq4$, the eigenvalues of $\rho(T)$ determine a unique equivalence class of irreducible\\
\>\>representations of $\G$. If $d=4$, there are at most two equivalence classes of irreducible\\
\>\>representations, with representatives $\rho_0$ and $\rho_1$, such that $\rho_0(T)$ and $\rho_1(T)$ have the\\
\>\>same eigenvalues as $\rho(T)$.
\end{tabbing}\qed
\end{thm}

\begin{cor}\label{cor:tw}If $\rho:\G\rightarrow\gln{d}$ is irreducible and $T$-unitarizable,
$d\leq5$, then $\rho(T)$ has distinct eigenvalues $\e{r_j}$, for some real numbers
$0\leq r_1<\cdots<r_d<1$.\end{cor}

\pf The only point to be made is that $T$-unitarizable implies that $\rho(T)$ can be diagonalized,
so the minimal polynomial of $\rho(T)$ must be separable. Since the minimal and characteristic
polynomials of $\rho(T)$ are the same, $\rho(T)$ must have $d$ distinct eigenvalues.

\qed

We will use Theorem~\ref{thm:tw} in a crucial way in the proofs of the main results below.
Unfortunately, \cite{TW} does not classify \emph{indecomposable} representations
of $\G$, and this presents a slight problem for us, since the method utilized in \cite{M2}
requires one to deduce the irreducibility of a given indecomposable $\rho$ from the
eigenvalues of $\rho(T)$. To accommodate this situation, we make the following

\begin{defn}
An irreducible representation $\rho:\G\rightarrow\gln{d}$ is
\textbf{T-determined} if the following condition holds:
\begin{tabbing}\hspace{2cm}\=If $\rho':\G\rightarrow\gln{d}$ is indecomposable and $\rho'(T)$
has the same\\
\>eigenvalues as $\rho(T)$, then $\rho'$ is irreducible.
\end{tabbing}
\end{defn}

Thus the $T$-determined representations are exactly the ones for which we may
utilize together the results of \cite{TW} and the techniques of \cite{M2}. Of course,
there is a slight loss of generality which accompanies the use of Definition \thedefn,
but qualitatively it is clear that almost all irreducible representations are $T$-determined.
Indeed, from Corollary~\ref{cor:chi} one concludes that, regardless of dimension, a representation
$\rho$ is $T$-determined if no proper sub-product of the eigenvalues of $\rho(T)$ is an $12^{th}$ root
of unity, since the determinant of any sub-representation is again a character of $\G$. See also
\cite[Th.\ 3.1]{M2}, which implies that \emph{every} irreducible $\rho:\G\rightarrow\gln{2}$
for which $\rho(T)$ is semisimple is in fact $T$-determined.

\section{Dimension one}

This is the classical setting, but using the vector-valued language to describe the known results
will serve as a useful warm-up to establishing analogous results in higher dimension.\medskip

Suppose $\rho:\G\rightarrow\gln{1}$ is a representation of $\G$, and let $\up$ be a multiplier
system with cusp parameter $m$. Obviously $\rho$ is irreducible (in fact $\rho\in\Hom(\G,\cstar)$
by definition), so by Corollary~\ref{cor:chi} $\rho=\chi^n$ for some integer $n$. Note that
\begin{eqnarray*}F\in\mc{H}(k,\rho,\up)&\Longleftrightarrow&\forall\gamma\in\G,\ F|_k^\up\gamma=\rho(\gamma)F\\ \\
&\Longleftrightarrow&\forall\gamma\in\G,\ F|_k^{\up\chi^n}\gamma=F\\ \\
&\Longleftrightarrow&F\in\mc{H}(k,\up\chi^n).
\end{eqnarray*}
In other words, a 1-dimensional vector-valued modular form for $\rho$ and $\up$ is exactly a classical
modular form for the multiplier system $\up\chi^n$; note that the definition of $\chi$ also implies that
$\rho$ is $T$-unitarizable, since $\rho(T)=\chi^n(T)=\e{\frac{n}{12}}$. By Corollary~\ref{cor:modmult}
we have
\begin{eqnarray*}
\mc{H}(\rho,\up)&=&\mc{H}(\up\chi^n)\\ \\
&=&\bigoplus_{k\geq0}\mc{H}(k_0+2k,\up\chi^n)\\ \\
&=&\m\eta^{2k_0},
\end{eqnarray*}
where $k_0$ denotes the cusp parameter of the multiplier system $\up\chi^n$. Thus
$\mc{H}(\rho,\up)$ is a free $\m$-module of rank 1, and by Lemma~\ref{lem:dker} the
generator $\eta^{2k_0}$ is a solution to the $1^{st}$ order MLDE $D_{k_0}f=0$. Furthermore,
$\mc{H}(\rho,\up)=\mc{R}\eta^{2k_0}$ is a cyclic $\mc{R}$-module, trivially so since
$D_{k_0}\eta^{2k}=0$, \ie in the dimension 1 setting the $\mc{R}$- and $\mc{M}$-module
structures of $\mc{H}(\rho,\up)$ coincide.\medskip

The remainder of this chapter is largely devoted to determining to what extent the above observations
generalize to higher dimension.

\section{Dimensions two and three}

\begin{prop}\label{prop:minwt23} Let $\rho:\G\rightarrow\gln{d}$ be a $T$-determined, $T$-unitarizable
representation of dimension $d=2$ or $3$, and $\up$ any multiplier system for $\G$. There
is a real number $k_0$, depending on $\rho$ and $\up$, such that
$$\mc{H}(\rho,\up)=\bigoplus_{k\geq0}\mc{H}(k_0+2k,\rho,\up),$$
and $\mc{H}(k_0,\rho,\up)$ contains a nonzero vector $F_0$, whose
components form a fundamental system of solutions of a $d^{th}$ order MLDE in weight $k_0$.
\end{prop}

\pf By Lemma~\ref{lem:vvequiv} (see also Remark~\ref{rem:vvequiv}) together with Corollary~\ref{cor:tw},
we are free to assume that
$$\rho(T)=\diag{\e{r_1},\cdots,\e{r_d}}$$
for some real numbers
\begin{equation}\label{eq:distinct}
0\leq r_1<\cdots<r_d<1.
\end{equation}
Let $m$ denote the cusp parameter of $\up$. A set
$\Lambda=\{\lambda_1,\cdots,\lambda_d\}$ of nonnegative real numbers will be called
\textbf{admissible} if
\begin{equation}\label{eq:adms}
\lambda_j\equiv r_j+\frac{m}{12}\pmod{\bb{Z}}
\end{equation}
for $j=1,\cdots,d$; the admissibility of $\Lambda$ is therefore a
necessary, but not sufficient, condition for a vector
$$\cvec{q^{\lambda_1}\sum_{n\geq0}a_1(n)q^n}{\vdots}{q^{\lambda_d}\sum_{n\geq0}a_d(n)q^n}$$
to be a nonzero form in $\mc{H}(\rho,\up)$. By Corollary~\ref{cor:umlde}, each admissible set
$\Lambda$ yields a unique $d^{th}$ order MLDE
\begin{equation}\label{eq:lambdamlde}
L_\Lambda[f]=0
\end{equation}
whose set of indicial roots is $\Lambda$. Note that \eqr{distinct} implies
$$\lambda_i\equiv\lambda_j\pmod{\bb{Z}}\ \Longleftrightarrow\ i=j.$$
By Corollary~\ref{cor:solnvec}, there is a vector
\begin{equation}\label{eq:Flambda}
F_\Lambda=\cvec{q^{\lambda_1}+\sum_{n\geq1}a_1(n)q^{\lambda_1+n}}{\vdots}{q^{\lambda_d}+
\sum_{n\geq1}a_d(n)q^{\lambda_d+n}}
\end{equation}
whose components are a fundamental system of solutions of \eqr{lambdamlde}, and if we
denote the weight of \eqr{lambdamlde} by $k_\lambda$, then for each $\up'\in\mult(k_\lambda)$,
there is a representation $\rho':\G\rightarrow\gln{d}$, arising from the $|_{k_\lambda}^{\up'}$-invariance
of the solution space of \eqr{lambdamlde}, such that $F_\Lambda\in\mc{H}(k_\lambda,\rho',\up')$.
Setting $\lambda=\lambda_1+\cdots+\lambda_d$ and applying Theorem~\ref{thm:wmlde} tells us that
the weight of \eqr{lambdamlde} is
\begin{equation}
k_\lambda=\frac{12\lambda}{d}+1-d.
\end{equation}
By Corollary~\ref{cor:equivmult}, $\up\in\mult(k_\lambda)$ iff $k_\lambda\equiv m$ (mod
$\bb{Z}$), and by \eq{0} this holds iff $d|(12\lambda-md)$. By Corollary~\ref{cor:detrho},
$12\lambda-md\equiv12r$ (mod $12\bb{Z}$), where $r=r_1+\cdots+r_d$.
Furthermore, Lemma~\ref{lem:div} says that $d|12r$, and since $d=2$ or 3, we
also know that $d|12$.  Therefore $d|(12\lambda-md)$, so that $\up\in\mult(k_\lambda)$. Thus
we may take $\up'=\up$, and we have $F_\Lambda\in\mc{H}(k_\lambda,\rho',\up)$, where $\rho'$
arises from the $|_{k_\lambda}^\up$ action of $\G$ on the solution space of
\eqr{lambdamlde}. Corollary~\ref{cor:indec} implies that $\rho'$ is indecomposable, and clearly
(or see Remark~\ref{rem:monodromy})
$$\rho'(T)=\diag{\e{r_1},\cdots,\e{r_d}}=\rho(T).$$
Since $\rho$ is assumed to be $T$-determined, we have by Theorem~\ref{thm:tw} that $\rho\cong\rho'$.
Again by Lemma~\ref{lem:vvequiv}, we are free to assume that $F_\Lambda\in\mc{H}(k_\lambda,\rho,\up)$.

Since each admissible set consists of nonnegative real numbers satisfying congruence conditions
(mod $\bb{Z}$), there is a unique admissible set $\Lambda_0$ such that each $\lambda_j$ is as
small as possible (note also that these $\lambda_j$ represent the minimal possible leading
exponents of the components of any nonzero vector in $\mc{H}(\rho,\up)$). If the sum of the
indicial roots in $\Lambda_0$ is denoted by $\lambda_0$, we obtain via
Theorems~\ref{thm:wetag},~\ref{thm:wmlde} the minimal possible weight
$$k_0=\frac{12\lambda_0}{d}+1-d$$
of a nonzero vector-valued modular form in $\mc{H}(\rho,\up)$, and we have that
$\mc{H}(k_0,\rho,\up)$ contains a nonzero vector whose components form a fundamental
system of solutions of a $d^{th}$ order MLDE in weight $k_0$.\qed

\begin{lem}\label{lem:23bound}
Under the hypotheses of Proposition~\ref{prop:minwt23}, we have
$$\dim\mc{H}(k_0+2k,\rho,\up)\leq\left[\frac{dk}{6}\right]+1,$$
for each $k\geq0$.
\end{lem}

\pf  We again assume, as we may, that $\rho(T)$ is diagonal. As before, let
$\Lambda_0=\{\lambda_1,\cdots,\lambda_d\}$ denote the set of minimal indicial
roots of the minimal weight MLDE $L_{\Lambda_0}[f]=0$, and
\begin{equation}
k_0=\frac{12\lambda_0}{d}+1-d
\end{equation}
the weight of $L_{\Lambda_0}$, where $\lambda_0=\lambda_1+\cdots+\lambda_d$.
Then an arbitrary vector in $\mc{H}(\rho,\up)$ has the form
\begin{equation}
F(z)=\cvec{q^{\lambda_1+n_1}\sum_{n\geq0}a_1(n)q^n}{\vdots}
{q^{\lambda_d+n_d}\sum_{n\geq0}a_d(n)q^n},
\end{equation}
for some nonnegative integers $n_1,\cdots,n_d$.\medskip

Suppose $\mc{H}(k_0+2k,\rho,\up)$ contains a nonzero $F$ as in
\eq{0}, and set $n=n_1+\cdots+n_d$. By Theorem~\ref{thm:wetag}, the modular
Wronskian of $F$ is
$$W(F)=\eta^{24(\lambda+n)}g,$$ for a unique modular form
$$g\in\mc{M}_{d(k_0+2k+d-1)-12(\lambda+n)}=\mc{M}_{2kd-12n}$$
which is not a cusp form. In particular, $g$ is nonzero, and this
means
$$2kd-12n\geq0,\neq2.$$
Since $d=2$ or $3$, it is clear that $2kd-12n\neq2$,
regardless of $k$ and $n$. Thus the only constraint on $n$ is given by the inequality
\begin{equation}
n\leq\left[\frac{dk}{6}\right].
\end{equation}
Fix $k\geq0$, and set $a(k)=\dim\mc{H}(k_0+2k,\rho,\up)$. For each
integer $j\geq0$, let
$$\phi_j:\mc{H}(k_0+2k,\rho,\up)\rightarrow\bb{C}$$
denote the linear functional which takes $F$ in \eq{-1} to
$\phi_j(F)=a_1(j)$, the $j^{th}$ Fourier coefficient of the first
component of $F$. Then $\dim\ker\phi_j\geq a(k)-1$ for all $j$,
which means
$$\bigcap_{j=0}^{a(k)-2}\ker\phi_j\neq\{0\}.$$
This is equivalent to saying there is a nonzero $F$ as in \eq{-1}
such that $n_1\geq a(k)-1$. Combining this inequality with \eq{0}
finishes the proof.\qed

\begin{thm}\label{thm:2mmod}
Assume $d=2$ in Proposition~\ref{prop:minwt23}. Then the following hold:
\begin{tabbing}\hspace{3cm}\=\ i) \=$\mc{H}(\rho,\up)$ is a free $\mc{M}$-module of rank 2, with basis
$\{F_0, DF_0\}$.\\ \\
\>ii) $\mc{H}(\rho,\up)=\mc{R}F_0$ is cyclic as $\mc{R}$-module.\\ \\
\>iii) For all $k\geq0$, $\dim{H}(k_0+2k,\rho,\up)=\left[\frac{k}{3}\right]+1$.
\end{tabbing}
\end{thm}

\pf By Lemma~\ref{lem:dfind}, the set $\{F_0, DF_0\}$ is independent over $\mc{M}$.
Using Lemma~\ref{lem:23bound}, we immediately obtain the following information:\smallskip

\begin{tabbing}\hspace{3cm}\=$\mc{H}(k_0,\rho,\up)=\langle F_0\rangle$ is 1-dimensional.\\ \\
\>$\mc{H}(k_0+2,\rho,\up)=\langle DF_0\rangle$ is 1-dimensional.\\ \\
\>$\mc{H}(k_0+4,\rho,\up)=\langle E_4F_0\rangle$ is 1-dimensional.\\ \\
\>$\mc{H}(k_0+6,\rho,\up)=\langle E_6F_0, E_4DF_0\rangle$ is 2-dimensional.\\ \\
\>$\mc{H}(k_0+8,\rho,\up)=\langle E_8F_0, E_6DF_0\rangle$ is 2-dimensional.\\ \\
\>$\mc{H}(k_0+10,\rho,\up)=\langle E_{10}F_0, E_8DF_0\rangle$ is 2-dimensional.
\end{tabbing}\smallskip

We now argue as in the proof of Lemma~\ref{lem:23bound}, defining linear functionals
$$\phi_j:\mc{H}(k_0+2k,\rho,\up)\rightarrow\bb{C},\hspace{.6cm} j=1,2$$
where $\phi_j$ takes a vector
$$F=\begin{pmatrix}q^{\lambda_1}\sum_{n\geq0}a_1(n)q^n\\ \\
q^{\lambda_2}\sum_{n\geq0}a_2(n)q^n\end{pmatrix}$$
to the Fourier coefficient $\phi_j(F)=a_j(0)$. If we set
$a(k)=\dim\mc{M}(k_0+2k,\rho,\up)$, then for each $j$ we have
$\dim\ker\phi_j\geq a(k)-1$.  Thus for any $k$ such that $a(k)\geq3$, we have
$$\ker\phi_1\cap\ker\phi_2\neq\{0\}.$$
In other words, $\Delta\mc{H}(k_0+2k-12,\rho,\up)$ has codimension at most
$2$ in $\mc{H}(k_0+2k,\rho,\up)$, for all $k\geq0$. Let
$$U_{2k}=\left\{\begin{array}{llc}\mc{H}(k_0+2k,\rho,\up),& &k<6\\ \\
\langle E_{2k}F_0, E_{2(k-1)}DF_0\rangle,& &k\geq6.\end{array}\right.$$
For $k<6$, we clearly have
\begin{equation}
\mc{H}(k_0+2k,\rho,\up)=U_{2k}\oplus\Delta\mc{H}(k_0+2(k-6),\rho,\up),
\end{equation}
and this forms the base case for the inductive assumption that \eq{0} holds for all $k'$
less than an arbitrary $k\geq6$. Since $\dim U_{2k}=2$ for $k\geq6$, to establish \eq{0} for
our $k$ it is sufficient to show that
\begin{equation}
U_{2k}\cap\Delta\mc{H}(k_0+2(k-6),\rho,\up)=\{0\}.
\end{equation}
Repeated use of \eq{-1} allows us to write
$$\Delta\mc{H}(k_0+2(k-6),\rho)=\Delta U_{2(k-6)}\oplus\cdots\oplus\Delta^nU_{2(k-6n)}$$
for some integer $n\geq1$, so verifying \eq{0} amounts to proving there is no nontrivial relation
\begin{equation}
\alpha_1E_{2k}F_0+\alpha_2E_{2(k-1)}DF_0=\Delta(M_{2(k-6)}F_0+M_{2(k-7)}DF_0),
\end{equation}
for some $\alpha_j\in\bb{C}$ and $M_j\in\mc{M}_j$. But again, $F_0$ and $DF_0$ are independent over
$\mc{M}$, so we must have
\begin{eqnarray*}\alpha_1E_{2k}=\Delta M_{2(k-6)},\\ \\
\alpha_2E_{2(k-1)}=\Delta M_{2(k-7)}.
\end{eqnarray*}
Consideration of $q$-expansions shows that all terms in the above equations must be identically zero,
and this proves that \eq{-2} and \eq{-1} hold for all $k\geq0$. Thus
\begin{equation}
\mc{H}(\rho,\up)=\mc{M}F_0\oplus\mc{M}DF_0
\end{equation}
is a free $\mc{M}$-module of rank 2. Statement ii) of the Theorem follows immediately from \eq{0},
and since $\dim U_{2k}=2$ for $k\geq6$, statement iii) follows inductively from \eq{-3}.\qed

\begin{cor} The Hilbert-Poincar\'e series \eqr{hp} of $\mc{H}(\rho,\up)$ is

\begin{eqnarray*}
\Psi(\rho,\up)(t)&=&\frac{t^{k_0}(1+t^2)}{(1-t^4)(1-t^6)}\\ \\
&=&\frac{t^{k_0}}{(1-t^2)(1-t^6)}.
\end{eqnarray*}\qed
\end{cor}

\begin{thm}\label{thm:3mmod}
Assume $d=3$ in Proposition~\ref{prop:minwt23}. Then the following hold:
\begin{tabbing}\hspace{3cm}\=\ \ i) \=$\mc{H}(\rho,\up)$ is a free $\mc{M}$-module of rank 3, with basis
$\{F_0, DF_0, D^2F_0\}$.\\ \\
\>\ ii) $\mc{H}(\rho,\up)=\mc{R}F_0$ is cyclic as $\mc{R}$-module.\\ \\
\>iii) For all $k\geq0$, $\dim{H}(k_0+2k,\rho,\up)=\left[\frac{k}{2}\right]+1$.
\end{tabbing}
\end{thm}

\pf The proof is essentially identical to that of the previous Theorem. One again finds from
Lemma~\ref{lem:dfind} that the set $\{F_0, DF_0, D^2F_0\}$ is independent over $\m$,
and this fact along with Lemma~\ref{lem:23bound} yields\medskip

\begin{tabbing}\hspace{4cm}\=$\mc{H}(k_0,\rho,\up)=\langle F_0\rangle$ is 1-dimensional.\\ \\
\>$\mc{H}(k_0+2,\rho,\up)=\langle DF_0\rangle$ is 1-dimensional.\\ \\
\>$\mc{H}(k_0+4,\rho,\up)=\langle E_4F_0, D^2F_0\rangle$ is 2-dimensional.\\ \\
\>$\mc{H}(k_0+6,\rho,\up)=\langle E_6F_0, E_4DF_0\rangle$ is 2-dimensional.\\ \\
\>$\mc{H}(k_0+8,\rho,\up)=\langle E_8F_0, E_6DF_0, E_4D^2F_0\rangle$ is 3-dimensional.\\ \\
\>$\mc{H}(k_0+10,\rho,\up)=\langle E_{10}F_0, E_8DF_0, E_6D^2F_0\rangle$ is 3-dimensional.
\end{tabbing}\medskip

We again define linear functionals $\phi_j$, $j=1,2,3$ and conclude that
$$\bigcap_{j=1}^3\ker\phi_j\neq\{0\}$$
whenever $\dim\mc{H}(k_0+2k,\rho,\up)\geq4$. Therefore $\Delta\mc{H}(k_0+2k-12,\rho,\up)$ has
codimension at most 3 in $\mc{H}(k_0+2k,\rho,\up)$, for all $k\geq0$. This time we define
$$U_{2k}=\left\{\begin{array}{llc}\mc{H}(k_0+2k,\rho,\up),& &k<6\\ \\
\langle E_{2k}F_0, E_{2(k-1)}DF_0, E_{2(k-2)}D^2F_0\rangle,& &k\geq6,\end{array}\right.$$
and we again have trivially that
\begin{equation}
\mc{H}(k_0+2k,\rho,\up)=U_{2k}\oplus\Delta\mc{H}(k_0+2(k-6),\rho,\up),
\end{equation}
for $k<6$. This again forms the base case for the inductive assumption that \eq{0} holds for all $k'$
less than an arbitrary $k\geq6$. Noting that $\dim U_{2k}=3$ for $k\geq6$, an identical argument
involving $q$-expansions finishes the induction, proving
\begin{equation}
\mc{H}(\rho,\up)=\mc{M}F_0\oplus\mc{M}DF_0\oplus\mc{M}D^2F_0
\end{equation}
is a free $\mc{M}$-module of rank 3. The other two statements follow immediately, as before.\qed

\begin{cor} The Hilbert-Poincar\'e series \eqr{hp} of $\mc{H}(\rho,\up)$ is

\begin{eqnarray*}
\Psi(\rho,\up)(t)&=&\frac{t^{k_0}(1+t^2+t^4)}{(1-t^4)(1-t^6)}\\ \\
&=&\frac{t^{k_0}}{(1-t^2)(1-t^4)}.
\end{eqnarray*}\qed
\end{cor}

\section{Dimension four}

Suppose $\rho:\G\rightarrow\gln{4}$ is a $T$-determined, $T$-unitarizable representation of $\G$.
By Lemma~\ref{lem:vvequiv} and Corollary~\ref{cor:tw}, we may assume that
$$\rho(T)=\diag{\e{r_1},\cdots,\e{r_4}}$$
for some real numbers
$$0\leq r_1<r_2<r_3<r_4<1.$$
Let $\up$ be a multiplier system with cusp parameter $m$. As before, each set $\Lambda$ of
\emph{admissible} (\cf\ \eqr{adms}) indicial roots, with corresponding sum $\lambda$, yields a unique (by Corollary~\ref{cor:umlde}) MLDE
\begin{equation}L_\Lambda[f]=0\end{equation}
in weight (by Theorem~\ref{thm:wmlde}) $k_\lambda=3\lambda-3$. By Lemma~\ref{lem:div}, we know that
$3r$ is an integer, where $r=r_1+\cdots+r_4$, and we again utilize Corollary~\ref{cor:detrho} to find that $12\lambda-4m\equiv12r\pmod{12\bb{Z}}$. Therefore $3\lambda-m\in\bb{Z}$, so that
$k_\lambda\equiv m\pmod{\bb{Z}}$ and $\up\in\mult(k_\lambda)$.

Again we may form the set $\Lambda_0=\{\lambda_1,\lambda_2,\lambda_3,\lambda_4\}$ of minimal possible admissible indicial roots, the sum of which we denote by $\lambda_0$, and obtain a vector-valued modular form
\begin{equation}\label{eq:f04}
F_0=\cvec{q^{\lambda_1}+\sum_{n\geq1}a_1(n)q^{\lambda_1+n}}{\vdots}
{q^{\lambda_4}+\sum_{n\geq1}a_4(n)q^{\lambda_4+n}}\end{equation}
whose components are a fundamental system of solutions of \eq{-1}. We have $F_0\in\mc{H}(k_0,\rho_0,\up)$, where
\begin{equation}\label{eq:k04}
k_0=3\lambda_0-3
\end{equation}
and $\rho_0$ denotes the (indecomposable by Corollary~\ref{cor:indec}) representation arising from the $|_{k_0}^\up$ action of $\G$ on the solution space of the associated MLDE $L_{\Lambda_0}[f]=0$. It is clear
that $\rho_0(T)=\rho(T)$, so $\rho_0$ is irreducible since $\rho$ is $T$-determined. Note that
$$\rho_0(S^2)F_0=F_0|_{k_0}^\up S^2,$$
so using the identity $S^2=-I$ we obtain
$$\rho_0(-I)=\up(-I)^{-1}(-1)^{-k_0}.$$
On the other hand, we may modify $\Lambda_0$ by fixing an arbitrary $i\in\{1,2,3,4\}$ and forming
a new admissible set $\Lambda_1=\{\lambda_j+\delta_{i,j}\}_{j=1}^4$. To this set of indicial roots
corresponds a second MLDE
\begin{equation}L_{\Lambda_1}[f]=0\end{equation}
in weight $3(\lambda_0+1)-3=k_0+3$. As above, $\up\in\mult(k_0+3)$, so we again may use Corollary~\ref{cor:solnvec}
to obtain a vector
\begin{equation}\label{eq:f1}
F_1=\cvec{q^{\lambda_1+\delta_{i,1}}\sum_{n\geq0}b_1(n)q^n}{\vdots}
{q^{\lambda_4+\delta_{i,4}}\sum_{n\geq0}b_4(n)q^n}\in\mc{H}(k_0+3,\rho_1,\up),
\end{equation}
where $\rho_1$ denotes the representation arising from the $|_{k_0+3}^\up$ action of $\G$ on the solution
space of \eq{0}, and $b_j(0)=1$ for each $j$. Again we have $\rho_1(T)=\rho(T)$, so $\rho_1$ is irreducible,
since $\rho$ is $T$-determined. In this case we observe that
$$\rho_1(S^2)F_1=F_1|_{k_0+3}^\up S^2,$$
so
$$\rho_1(-I)=\up(-I)^{-1}(-1)^{-k_0-3}=-\rho_0(-I).$$
From the proof of Lemma~\ref{lem:div}, we see that $\rho_0(S)$ and $\rho_1(S)$ have different
eigenvalues, either $\{1,1,-1,-1\}$, $\{i,i,-i,-i\}$ respectively, or vice versa. Therefore
$\rho_0$ and $\rho_1$ are inequivalent, and represent the two equivalence classes predicted by
Theorem~\ref{thm:tw} in the dimension four setting. (Alternatively, we know from
Lemma~\ref{lem:2zindec} that
$$\mc{H}(\rho_0,\up)=\bigoplus_{k\in\bb{Z}}\mc{H}(k_0+2k,\rho_0,\up),$$
whereas
$$\mc{H}(\rho_1,\up)=\bigoplus_{k\in\bb{Z}}\mc{H}(k_0+3+2k,\rho_1,\up).$$
Clearly $k_0+3\not\equiv k_0\pmod{2}$, so by Lemma~\ref{lem:vvequiv} we know
that $\rho_0\not\cong\rho_1$.)

From Theorem~\ref{thm:tw}, we conclude that $\rho\cong\rho_j$ for a unique $j\in\{0,1\}$.
Thus to determine the structure of $\mc{H}(\rho,\up)$, it is sufficient (by Lemma~\ref{lem:vvequiv})
to do so for $\mc{H}(\rho_j,\up)$, $j=0,1$. We turn now to this task:

\subsection{Structure of $\mc{H}(\rho_0,\up)$}

This case is completely analogous to that of dimensions two and three, \ie
\begin{thm}\label{thm:40mmod}Let $F_0$, $k_0$ be as in \eqr{f04}, \eqr{k04} respectively. Then
\begin{tabbing}\hspace{1cm}\=\ \ i) \=$\mc{H}(\rho_0,\up)=\bigoplus_{k\geq0}\mc{H}(k_0+2k,\rho_0,\up)$, and
$\mc{H}(k_0,\rho_0,\up)$ contains the vector\\
\>\>$F_0$, whose components form a fundamental system of solutions of an\\
\>\>MLDE in weight $k_0$.\\ \\
\>\ ii) $\mc{H}(\rho_0,\up)$ is a free $\mc{M}$-module of rank 4, with basis $\{F_0,DF_0,D^2F_0,D^3F_0\}$.\\ \\
\>iii) $\mc{H}(\rho_0,\up)=\mc{R}F_0$ is cyclic as $\mc{R}$-module.\\ \\
\>\ iv) For each $k\geq0$, $\dim\mc{H}(k_0+2k,\rho_0,\up)=\left[\frac{2k}{3}\right]+1$.\end{tabbing}\qed
\end{thm}

\pf Suppose that $0\neq F\in\mc{H}(k_0+2k,\rho_0,\up)$ for some $k\in\bb{Z}$. Then the
components of $F$ have leading exponents $\lambda_j+n_j$ for some nonnegative integers $n_1,\cdots,n_4$,
and setting $n=\sum n_j$, we have by Theorem~\ref{thm:wetag} that $W(F)=\eta^{24(\lambda_0+n)}g$,
where the non-cusp form $g$ has weight
\begin{eqnarray*}
wt(g)&=&4(k_0+2k+3)-12(\lambda_0+n)\\ \\
&=&4(3\lambda_0-3+2k+3)-12(\lambda_0+n)\\ \\
&=&8k-12n.
\end{eqnarray*}
Because $g$ is nonzero and $n$ is nonnegative, we find that $k\geq0$ and statement $i)$ holds.
Note also that $4|wt(g)$ regardless of $k$ and $n$, so $wt(g)\neq2$. We therefore obtain the inequality
$$n\leq\left[\frac{2k}{3}\right],$$
and by the same argument as in the proof of Lemma~\ref{lem:23bound} we find that
\begin{equation}
\dim\mc{H}(k_0+2k,\rho_0,\up)\leq\left[\frac{2k}{3}\right]+1,
\end{equation}
for each $k\geq0$.

The remainder of the proof is identical to that of Theorems~\ref{thm:2mmod},~\ref{thm:3mmod}, \ie
after using Lemma~\ref{lem:dfind} to conclude that $\{F_0, DF_0, D^2F_0, D^3F_0\}$ is an independent
set of vectors over $\m$, one applies the bound \eq{0} and produces the base-case data\medskip

\begin{tabbing}\hspace{2cm}\=$\mc{H}(k_0,\rho_0,\up)=\langle F_0\rangle$ is 1-dimensional.\\ \\
\>$\mc{H}(k_0+2,\rho_0,\up)=\langle DF_0\rangle$ is 1-dimensional.\\ \\
\>$\mc{H}(k_0+4,\rho_0,\up)=\langle E_4F_0, D^2F_0\rangle$ is 2-dimensional.\\ \\
\>$\mc{H}(k_0+6,\rho_0,\up)=\langle E_6F_0, E_4DF_0, D^3F_0\rangle$ is 3-dimensional.\\ \\
\>$\mc{H}(k_0+8,\rho_0,\up)=\langle E_8F_0, E_6DF_0, E_4D^2F_0\rangle$ is 3-dimensional.\\ \\
\>$\mc{H}(k_0+10,\rho_0,\up)=\langle E_{10}F_0, E_8DF_0, E_6D^2F_0, E_4D^3F_0\rangle$ is 4-dimensional,
\end{tabbing}\medskip

and statements $ii)-iv)$ follow exactly as in the other Theorems.\qed

\begin{cor} The Hilbert-Poncar\'e series \eqr{hp} of $\mc{H}(\rho_0,\up)$ is
$$\Psi(\rho,\up)(t)=\frac{t^{k_0}(1+t^2+t^4+t^6)}{(1-t^4)(1-t^6)}.$$
\qed
\end{cor}

\subsection{Structure of $\mc{H}(\rho_1,\up)$}

Recall the representation $\rho_1$ from \eqr{f1}. We shall find that $\mc{H}(\rho_1,\up)$ is
again free (of rank 4) as $\mc{M}$-module, \emph{but is no longer cyclic as $\mc{R}$-module.}
Thus $d=4$ is the lowest dimension for which a $T$-unitarizable,
$T$-determined representation $\rho:\G\rightarrow\gln{d}$ is not necessarily cyclic as $\mc{R}$-module
(\cf\ Remark~\ref{rem:cyclic} below). We begin with
\begin{lem}\label{lem:4bound}Let $k_1=k_0+1$, with $k_0$ as in \eqr{k04}. Then

\begin{tabbing}\hspace{2cm}\=\ i) $\mc{H}(\rho_1,\up)=\bigoplus_{k\geq0}\mc{H}(k_1+2k,\rho_1,\up)$.\\ \\
\>ii) For all $k\geq0$, $\dim\mc{H}(k_1+2k,\rho_1,\up)\leq\left[\frac{2k+1}{3}\right]+1$.
\end{tabbing}

\end{lem}

\pf Let $k$ be an integer such that $\mc{H}(k_1+2k,\rho_1,\up)$ contains a nonzero vector $F$. Since
$\{\lambda_1,\lambda_2,\lambda_3,\lambda_4\}$ represent the minimal nonnegative leading exponents
which can occur for such a vector, there are nonnegative integers $p_j$ such that
$$F=\cvec{q^{\lambda_1+p_1}\sum_{n\geq0}c_1(n)q^n}{\vdots}{q^{\lambda_4+p_4}\sum_{n\geq0}c_4(n)q^n}.$$
Set $p=\sum p_j$, and recall the notation $\lambda_0=\sum\lambda_j$. By Theorem~\ref{thm:wetag},
the modular Wronskian of $F$ has the form $\eta^{24(\lambda_0+p)}g$, where $g$ is a modular form of weight $$4(k_1+2k+3)-12(\lambda+p)$$
which is not a cusp form. Substituting $k_1=k_0+1=3\lambda_0-2$, we obtain
$$wt(g)=8k-12p+4\geq0,\neq2.$$
It is clear that no integers $k,p$ will make $wt(g)$ equal to 2, since $4|wt(g)$. Thus for a given
$k$, the only constraint on $p$ is that given by the inequality above, namely
\begin{equation}\label{eq:pbound}
p\leq\left[\frac{2k+1}{3}\right].
\end{equation}
Once again we argue as in the proof of Lemma~\ref{lem:23bound}, via linear functionals
and Fourier coefficients, and this finishes the proof.\qed

\begin{lem}\label{lem:k1base}There is a nonzero vector $G\in\mc{H}(k_1,\rho_1,\up)$ such that

\begin{tabbing}\hspace{3cm}\=\ \ i) $\mc{H}(k_1,\rho_1,\up)=\langle G\rangle$ is 1-dimensional,\\ \\
\>\ ii) $\mc{H}(k_1+2,\rho_1,\up)=\langle F_1,DG\rangle$ is 2-dimensional,\\ \\
\>iii) $\mc{H}(k_1+4,\rho_1,\up)=\langle E_4G,D^2G\rangle$ is 2-dimensional,\end{tabbing}
where $F_1\in\mc{H}(k_1+2,\rho_1,\up)$ is defined in \eqr{f1}.
\end{lem}\smallskip

\pf The subspace
$$V=\langle E_{10}F_1,E_8DF_1,E_6D^2F_1,E_4D^3F_1\rangle\leq\mc{H}(k_1+12,\rho_1,\up)$$
is 4-dimensional by Lemma~\ref{lem:dfind}, and from \eqr{f1} it is clear that each nonzero
vector in $V$ has the form
\begin{equation}
F=\begin{pmatrix}q^{\lambda_1+n_1}\sum_{n\geq0}a_1(n)q^n\\ \ \\q^{\lambda_2+n_2}\sum_{n\geq0}a_2(n)q^n\\ \ \\
q^{\lambda_3+n_3}\sum_{n\geq0}a_3(n)q^n\\ \ \\q^{\lambda_4+n_4}\sum_{n\geq0}a_4(n)q^n\end{pmatrix}
\end{equation}
with $n_j\geq\delta_{i,j}$. For each $j\neq i$, let $\phi_j:V\rightarrow\bb{C}$ denote the linear functional
which takes an $F$ as in \eq{0} to $\phi_j(F)=a_j(0)$. Then $\ker\phi_j$ has codimension
at most 1 in $V$, so $\cap_{j\neq i}\ker\phi_j\neq\{0\}$. In other words,
there is a nonzero $F\in\mc{H}(k_1+12,\rho_1,\up)$ of the form \eq{0} with $n_j\geq1$, $a_j(0)\neq0$
for $j=1,2,3,4$. In particular, such an $F$ is divisible by $\Delta$, say $F=\Delta G$ for some nonzero
$G$ in $\mc{H}(k_1,\rho_1,\up)$. Using the upper bound from Lemma~\ref{lem:4bound} in the $k=0$ case, we
obtain the first statement of the Lemma. Note also that \eqr{pbound} implies that $G$ has the form
\eq{0} with $n_j=0$ for each $j$.

The second statement also follows from Lemma~\ref{lem:4bound}, so long as $DG$ and $F_1$ are
linearly independent. Recalling the q-expansion
$$P(q)=-\frac{1}{12}+\cdots,$$
we compute
$$DG=D_{k_1}G=q\frac{d}{dq}G+k_1PG$$
and find that the $i^{th}$ component of $D_{k_1}G$ has the form
$$g_i(q)=a_i(0)\left(\lambda_i-\frac{k_1}{12}\right)q^{\lambda_i}+\cdots,$$
as opposed to the $i$th component of $F_1$, which is
$$f_i(q)=q^{\lambda_i+1}+\cdots.$$ For $DG$ to be a scalar
multiple of $F_1$, it is therefore required that $\lambda_i=\frac{k_1}{12}$, since $a_i(0)\neq0$ by
assumption. Recalling that the $\lambda_j$ are distinct, there is at most one $j$ for which
$\lambda_j=\frac{k_1}{12}$. Furthermore, the choice of $i\in\{1,2,3,4\}$ was completely arbitrary.
Indeed, all of the calculations above depend only on the sum $\lambda_0$, and although each choice
of $i$ will yield a different MLDE in weight $k_1+2$, it is clear by the definition of $T$-determined
and by Theorem~\ref{thm:tw} that the representations arising from the $|_{k_1+2}^\up$ action
of $\G$ on the respective solution spaces are all equivalent to $\rho_1$. Thus we are free to assume
that $i$ was chosen so that $\lambda_i\neq\frac{k_1}{12}$, and we have statement $ii)$ of the Lemma.

The third statement follows from Lemma~\ref{lem:4bound} and Lemma~\ref{lem:dfind}.\qed

\begin{lem} There is a nontrivial relation
\begin{equation}
DF_1=\alpha_1D^2G+\alpha_2E_4G,\ \ \alpha_1\alpha_2\neq0.
\end{equation}
\end{lem}

\pf Corollary~\ref{cor:dkinj} implies that $DF_1\neq0$, so statement $iii)$ of Lemma~\ref{lem:k1base}
implies the existence of a nontrivial relation \eq{0}. It is clear that $\alpha_1\neq0$, since this
would cause the leading exponents of the $i^{th}$ components of the left and right sides of
\eq{0} to be unequal. Statement $ii)$ of Lemma~\ref{lem:k1base} and the injectivity of $D$
imply $\alpha_2\neq0$.\qed

\begin{lem}\label{lem:41free} The set $\{G,DG,F_1,D^2G\}$ is independent over $\mc{M}$, and generates
a rank four $\mc{M}$-submodule of $\mc{H}(\rho_1,\up)$.
\end{lem}

\pf Suppose there is a (homogeneous) relation
\begin{equation}\label{eq:homeq}
M_kF_1=M_{k-2}D^2G+Q_kDG+M_{k+2}G
\end{equation}
where $M_j,Q_j\in\mc{M}_j$. If $M_k\neq0$, we may divide and obtain
\begin{equation}
F_1=\frac{M_{k-2}}{M_k}D^2G+\frac{Q_k}{M_K}DG+\frac{M_{k+2}}{M_k}G.
\end{equation}
Applying $D$ to each side of \eq{0} and substituting with \eq{-2} shows that $G$ satisfies a
differential equation
$$\scriptstyle{\frac{M_{k-2}}{M_k}D^3G+\left[D\left(\frac{M_{k-2}}{M_k}\right)+\frac{Q_k}{M_k}-\alpha_1\right]D^2G+
\left[D\left(\frac{Q_k}{M_k}\right)+\frac{M_{k+2}}{M_k}\right]DG+
\left[D\left(\frac{M_{k+2}}{M_k}\right)-\alpha_2E_4\right]G=0.}$$
By Lemma~\ref{lem:dfind}, each coefficient function in the above equation is identically zero and, in particular,
$M_{k-2}=0$. Setting the second coefficient equal to zero then yields $Q_k=\alpha_1M_k$. Substituting
into the third coefficient and setting it equal to zero reveals $M_{k+2}=-D(\alpha_1)M_k=0$, and
this in turn gives the data $\alpha_2E_4=0$ from the fourth coefficient. But $\alpha_2\neq0$ by
the previous Lemma, so we have reached an impasse by assuming $M_k\neq0$. On the other hand, if
$M_k=0$ then Lemma~\ref{lem:dfind} again implies that all coefficients in \eqr{homeq} are zero.\qed

Using the result just established along with Lemma~\ref{lem:4bound}, we obtain
\begin{cor}\begin{tabbing}\hspace{4cm}\=$\mc{H}(k_1+6,\rho_1,\up)=\langle E_6G,E_4DG,E_4F_1\rangle$,\\ \\
\>$\mc{H}(k_1+8,\rho_1,\up)=\langle E_8G,E_6DG,E_6F_1,E_4D^2G\rangle$,\\ \\
\>$\mc{H}(k_1+10,\rho_1,\up)=\langle E_{10}G,E_8DG,E_8F_1,E_6D^2G\rangle$.
\end{tabbing}\qed\end{cor}

We are now able to prove
\begin{thm}\label{thm:41mmod}
\begin{tabbing}\hspace{3cm}\=\ i) \=$\mc{H}(\rho_1,\up)$ is free of rank four as $\mc{M}$-module, with\\
\>\>basis $\{G,DG,F_1,D^2G\}$.\\ \\
\>ii) For all $k\geq0$, $\dim\mc{H}(k_1+2k,\rho_1,\up)=\left[\frac{2k+1}{3}\right]+1.$\\ \\
\end{tabbing}
\end{thm}

\pf We continue to use the same reasoning as that of Theorems~\ref{thm:2mmod}, \ref{thm:3mmod},
\ref{thm:40mmod}, (as originally appeared in \cite[Lem.\ 5.4]{M2}), whereby one lifts the
``codimension 1'' argument from the proof of \eqr{deltarecursion} to ``codimension $d$'':\medskip

If $\dim\mc{H}(k_1+2k,\rho_1,\up)\geq5$ for some $k\geq0$, then as with the construction of $G$,
we can intersect kernels of linear functions $\phi_j$, $j=1,2,3,4$ to produce a nonzero $F$
in $\mc{H}(k_1+2k,\rho_1,\up)$ which is divisible by $\Delta$. Therefore $\Delta\mc{H}(k_1+2(k-6),\rho_1,\up)$
has codimension at most 4 in $\mc{H}(k_1+2k,\rho_1,\up)$, for any $k\geq0$. For each $k\geq6$ we define
the 4-dimensional subspace
$$U_{2k}=\langle E_{2k}G,E_{2(k-1)}DG,E_{2(k-1)}DF_1,E_{2(k-2)}D^2G\rangle\leq\mc{H}(k_1+2k,\rho_1,\up),$$
and for $k\leq5$ we set $U_{2k}=\mc{H}(k_1+2k,\rho_1,\up)$. Then
\begin{equation}
U_{2k}\cap\Delta\mc{H}(k_1+2(k-6),\rho_1,\up)=\{0\}
\end{equation}
holds trivially for $k\leq5$, and this is the base case for the inductive assumption
that \eq{0} holds for all $k'$ less than an arbitrary $k\geq6$. Repeated use of this assumption
allows us to write
$$\Delta\mc{H}(k_1+2(k-6),\rho_1,\up)=\Delta U_{2(k-6)}\oplus\cdots\oplus\Delta^nU_{2(k-6n)},$$
so the truth of \eq{0} is equivalent to showing that the coefficients of any relation
$$\scriptstyle{\gamma_1E_{2k}G+\gamma_2E_{2(k-1)}F_1+\gamma_3E_{2(k-1)}DG+\gamma_4E_{2(k-2)}D^2G=\Delta(M_{2(k-6)}G+
M_{2(k-7)}F_1+Q_{2(k-7)}DG+M_{2(k-8)}D^2G)}$$
must be zero, where $M_j,Q_j\in\m_j$. But Lemma~\ref{lem:41free} implies that
\begin{tabbing}\hspace{4cm}\=$\gamma_jE_{2(k-i)}=\Delta M_{2(k-j-6)},\hspace{.5cm}j=1,2,4$\\ \\
\>$\gamma_3E_{2(k-1)}=\Delta Q_{2(k-7)}$\end{tabbing}
and examining $q$-expansions implies immediately that each term is zero in the
above equations. Thus \eq{0} holds for all $k\geq0$, and we have
\begin{equation}
\mc{H}(k_1+2k,\rho_1,\up)=U_{2k}\oplus\Delta\mc{H}(k_1+2(k-6),\rho_1,\up).
\end{equation}
The remaining statements of the Theorem follow immediately from \eq{0}.\qed

\begin{cor}  The Hilbert-Poincar\'e series \eqr{hp} of $\mc{H}(\rho_1,\up)$ is

$$\Psi(\rho_1,\up)(t)=\frac{t^{k_1}(1+2t^2+t^4)}{(1-t^4)(1-t^6)}.$$
\qed
\end{cor}

\begin{rem}\label{rem:cyclic} We have just seen that the representations we have
called $\rho_1$ provide the lowest dimension examples of $T$-unitarizable, $T$-determined
representations of $\G$ whose $\mc{R}$-module of vector-valued modular
forms is not cyclic. This is part of a more general phenomenon, which is addressed in
\cite[Th.\ 3]{MM} (in the trivial multiplier setting): there it is shown that for a
$T$-unitarizable (not necessarily irreducible) space of vector-valued modular forms to be a
cyclic $\mc{R}$-module, it is necessary and sufficient
that the lowest weight space be spanned by a vector whose components form a fundamental system
of solutions of an MLDE.\qed\end{rem}

\section{Dimension five}

Suppose $\rho:\G\rightarrow\gln{5}$ is a $T$-unitarizable, $T$-determined representation, and
fix a multiplier system $\up$ with cusp parameter $m$. Assume without loss that
$\rho(T)=\diag{\e{r_1},\cdots,\e{r_5}}$, where
\begin{equation}\label{eq:evals}
0\leq r_1<\cdots<r_5<1,
\end{equation}
and let $\{\widetilde{\lambda_1},\cdots,
\widetilde{\lambda_5}\}$
denote the minimal admissible set (\cf\ \eqr{adms}) for $\rho$ and $\up$. Then every vector in
$\mc{H}(\rho,\up)$ has the form
\begin{equation}\label{eq:f5}
F=\cvec{q^{\widetilde{\lambda_1}}\sum_{n\geq0}a_1(n)q^n}{\vdots}{q^{\widetilde{\lambda_5}}\sum_{n\geq0}a_5(n)q^n},
\end{equation}
and if $\mc{H}(k,\rho,\up)$ contains a nonzero vector \eq{0}, then by Theorem~\ref{thm:wetag}
we have
$$k\geq\frac{12}{5}\lambda_0-4,$$
where $\lambda_0=\sum\widetilde{\lambda_j}$.

On the other hand, let $\{n_1,\cdots,n_5\}$ be an arbitrary set of nonnegative integers, and set
\begin{equation}\label{eq:lambdaj}
\lambda_j=\widetilde{\lambda_j}+n_j,
\end{equation}
$\Lambda=\{\lambda_1,\cdots,\lambda_5\}$. From \eqr{adms}, \eqr{evals} it is clear that the
$\lambda_j$ are incongruent (mod $\bb{Z}$), and by Corollary~\ref{cor:umlde} there is a unique $5^{th}$ order MLDE
\begin{equation}\label{eq:5mlde}
L_\Lambda[f]=0
\end{equation}
whose set of indicial roots is  $\Lambda$. Let
$$F_\Lambda=\cvec{q^{\lambda_1}+\sum_{n\geq1}b_1(n)q^{\lambda_1+n}}{\vdots}{q^{\lambda_5}+
\sum_{n\geq1}b_5(n)q^{\lambda_5+n}}$$
be a vector whose components form a fundamental system of solutions of \eq{0}, and set $\lambda=\sum
\lambda_j$. By Theorem~\ref{thm:wmlde}, the weight of \eq{0} is
\begin{eqnarray*}
k_\lambda&=&\frac{12}{5}\lambda-4\\ \\
&=&\frac{12}{5}\left(\frac{5m}{12}+r+l+n\right)-4\\ \\
&=&m+\frac{12(r+l+n)}{5}-4,
\end{eqnarray*}
where for each $j$ we have
\begin{eqnarray*}
\lambda_j&=&\widetilde{\lambda_j}+n_j\\
&=&r_j+\frac{m}{12}+l_j
\end{eqnarray*}
for some integers $-1\leq l_j\leq0$, and $l=\sum l_j$, $r=\sum r_j$, $n=\sum n_j$. Note that
\begin{eqnarray}\up\in\mult(k_\lambda)&\Longleftrightarrow&k_\lambda\equiv m\pmod{\bb{Z}}\nonumber\\
\nonumber\\
&\Longleftrightarrow&\frac{12(r+l+n)}{5}\in\bb{Z}.
\end{eqnarray}
Since $12r$ is an integer (Lemma~\ref{lem:12r}) and $(5,12)=1$, there is a unique class of integers
(mod 5) whose members $n$ satisfy \eq{0}, thus there is a unique $N\in\{0,1,2,3,4\}$ such that
\eq{0} is satisfied when $n=N$.  If we now assume in \eqr{lambdaj} that the $n_j$ are chosen
so that $\sum n_j=N$, then the $|_{k_\lambda}^\up$ action of $\G$ on the solution space of
\eqr{5mlde} yields a representation $\rho_\Lambda:\G\rightarrow\gln{5}$ which is indecomposable by
Corollary~\ref{cor:indec}, and satisfies $\rho_\Lambda(T)=\rho(T)$. Since $\rho$ is $T$-determined,
we conclude that $\rho_\Lambda$ is irreducible, so by Theorem~\ref{thm:tw}, $\rho$ is equivalent
to $\rho_\lambda$.

Note that any other choice $n'_j$ of nonnegative integers satisfying $\sum n'_j=N$
yields some other set $\Lambda'=\{\widetilde{\lambda_j}+n'_j\}$ of indicial roots, an MLDE
$L_{\Lambda'}[f]=0$ whose set of indicial roots is $\Lambda'$, and a representation
$\rho_{\Lambda'}:\G\rightarrow\gln{5}$ which, by the arguments just given, is also equivalent to
$\rho$. Thus to classify $\mc{H}(\rho,\up)$, it is sufficient to classify $\mc{H}(\rho_\Lambda,\up)$,
for \emph{any} desired set of admissible indicial roots $\Lambda=\{\widetilde{\lambda_j}+n_j\}$
satisfying $\sum n_j=N$. In the remainder of this Section, we work out the five possible cases,
corresponding to the integers $N=0,1,2,3,4$.

Before proceeding to the various classifications, we prove the following
\begin{lem}\label{lem:5bound} Assume $\rho$, $\up$, $\lambda_0$, and $N$ are as above, and set
$$k_N=\frac{12(\lambda_0+N)}{5}-4.$$
Then the following hold:

\begin{tabbing}\=\ \ i) $\mc{H}(\rho,\up)=\bigoplus_{k\in\bb{Z}}\mc{H}(k_N+2k,\rho,\up).$\\ \\ \\
\>\ ii) $\mc{H}(k_N+2k,\rho,\up)=\{0\}$ for $k<-\frac{6N}{5}$.\\ \\
\>iii) For all $k\geq-\frac{6N}{5}$, $\dim\mc{H}(k_N+2k,\rho,\up)\leq\left\{\begin{array}{lll}\left[\frac{5k}{6}\right]+N&\ &k\equiv 5\pmod{6}\\ \\
\left[\frac{5k}{6}\right]+N+1&\ &k\not\equiv 5\pmod{6}\end{array}\right.$.\\ \\
\end{tabbing}
\end{lem}

\pf Statement $i)$ is clear from the definition of $N$ and the fact that $\mc{H}(\rho,\up)$ is isomorphic
(by Lemma~\ref{lem:vvequiv} and Theorem~\ref{thm:tw}) to $\mc{H}(\rho_{\Lambda},\up)$, for some representation $\rho_\Lambda$ arising from the $|_{k_N}^\up$-invariance of the solution space of an MLDE
$L_\Lambda[f]=0$ in weight $k_N$.

Suppose $k$ is an integer such that $\mc{H}(k_N+2k,\rho,\up)$ contains a nonzero vector $F$. Then there are
nonnegative integers $p_j$ and nonzero complex numbers $a_j(0)$ such that $F$ has the form
\begin{equation}\label{eq:fbound}
F=\cvec{q^{\widetilde{\lambda_1}+p_1}\sum_{n\geq0}a_1(n)q^n}{\vdots}
{q^{\widetilde{\lambda_5}+p_5}\sum_{n\geq0}a_5(n)q^n}.
\end{equation}
Set $p=\sum p_j$. By Theorem~\ref{thm:wetag}, the modular Wronskian
of $F$ has the form $W(F)=\eta^{24(\lambda_0+p)}g$, and the weight of
the non-cusp form $g$ satisfies
\begin{eqnarray*}
wt(g)&=&5(k_N+2k+4)-12(\lambda_0+p)\\ \\
&=&5\left(\frac{12}{5}(\lambda_0+N)-4+2k+4\right)-12(\lambda_0+p)\\ \\
&=&12(N-p)+10k.
\end{eqnarray*}
Since $p\geq0$, setting $wt(g)\geq0$ immediately implies statement $ii)$ of the Lemma.
Furthermore, we clearly have
$$wt(g)\equiv2\pmod{12}\ \Longleftrightarrow\ k\equiv5\pmod{6},$$
and this gives the inequalities
\begin{equation}\label{eq:wbound}
p\leq\left\{\begin{array}{lll}\left[\frac{5k}{6}\right]+N-1&\ &k\equiv 5\pmod{6}\\ \\
\left[\frac{5k}{6}\right]+N&\ &k\not\equiv 5\pmod{6}\end{array}\right..
\end{equation}
Arguing as in Lemma~\ref{lem:23bound} (via linear functionals and Fourier coefficients)
completes the proof of statement $iii)$, and this establishes the Lemma.\qed

\subsection{N=0}

By Lemma~\ref{lem:5bound}, we know that
$$\mc{H}(\rho,\up)=\oplus_{k\geq0}\mc{H}(k_0+2k,\rho,\up),$$
where $k_0=\frac{12}{5}\lambda_0-4$, and we may assume, up to equivalence of representation,
that $\mc{H}(k_0,\rho,\up)$ contains a vector
$$F_0=\cvec{q^{\widetilde{\lambda_1}}+\cdots}{\vdots}{q^{\widetilde{\lambda_5}}+\cdots}$$
whose components are a fundamental system of solutions of a $5^{th}$ order MLDE in weight $k_0$.
As predicted in Remark~\ref{rem:cyclic} above, we will once again find that
$\mc{H}(\rho_0,\up)$ is a cyclic $\mc{R}$-module:
\begin{thm}\label{thm:05} The following hold:
\begin{tabbing}\=\ \ i) $\mc{H}(\rho,\up)$ is a free $\m$-module of rank 5, with basis $\{F_0,DF_0,D^2F_0,D^3F_0,D^4F_0\}$.\\ \\
\>\ ii) For all $k\geq0$, $\dim\mc{H}(k_0+2k,\rho,\up)=\left\{\begin{array}{lll}\left[\frac{5k}{6}\right]&\ &k\equiv 5\pmod{6}\\ \\
\left[\frac{5k}{6}\right]+1&\ &k\not\equiv 5\pmod{6}\end{array}\right.$\\ \\
\>iii) $\mc{H}(\rho,\up)=\mc{R}F_0$ is cyclic as $\mc{R}$-module.\end{tabbing}
\end{thm}

\pf At this point, the method of proof should be all too familiar to the reader. We begin by
noting that the proposed basis of statement $i)$ is independent over $\m$, by Lemma~\ref{lem:dfind}.
We then use Lemma~\ref{lem:5bound} to establish the base data

\begin{tabbing}\hspace{4cm}\=$\mc{H}(k_0,\rho,\up)=\langle F_0\rangle$\\ \\
\>$\mc{H}(k_0+2,\rho,\up)=\langle DF_0\rangle$\\ \\
\>$\mc{H}(k_0+4,\rho,\up)=\langle E_4F_0,D^2F_0\rangle$\\ \\
\>$\mc{H}(k_0+6,\rho,\up)=\langle E_6F_0,E_4DF_0,D^3F_0\rangle$\\ \\
\>$\mc{H}(k_0+8,\rho,\up)=\langle E_8F_0,E_6DF_0,E_4D^2F_0,D^4F_0\rangle$\\ \\
\>$\mc{H}(k_0+10,\rho,\up)=\langle E_{10}F_0,E_8DF_0,E_6D^2F_0,E_4D^3F_0\rangle$.
\end{tabbing}

We define the subspaces $U_{2k}\in\mc{H}(k_0+2k,\rho,\up)$ to be

$$U_{2k}=\left\{\begin{array}{lcl}\mc{H}(k_0+2k,\rho,\up&\ &k<6\\ \\
\langle E_{2(k-j)}D^jF_0\rangle_{j=0}^4&\ &k\geq6\end{array}\right.,$$

and then prove via the standard arguments that
$$\mc{H}(k_0+2k,\rho,\up)=U_{2k}\oplus\Delta\mc{H}(k_0+2(k-6),\rho,\up)$$
for all $k\in\bb{Z}$. This establishes $i)-iii)$ of the Theorem.\qed

\begin{cor} The Hilbert-Poincar\'e series \eqr{hp} for $\mc{H}(\rho,\up)$ is

$$\Psi(\rho,\up)(t)=\frac{t^{k_0}(1+t^2+t^4+t^6+t^8)}{(1-t^4)(1-t^6)}.$$
\end{cor}\qed

\subsection{N=1}

In this case we have by Lemma~\ref{lem:5bound} that
$$\mc{H}(\rho,\up)=\bigoplus_{k\geq0}\mc{H}(k_1+2k,\rho,\up),$$
where $k_1=\frac{12}{5}(\lambda_0+1)-4$. Set $\Lambda=\{\widetilde{\lambda_1}+1,\widetilde{\lambda_2},
\widetilde{\lambda_3},\widetilde{\lambda_4},\widetilde{\lambda_5}\}$. By Corollary~\ref{cor:umlde},
there is a unique $5^{th}$ order MLDE $L_{\Lambda}[f]=0$ whose set of indicial roots is $\Lambda$,
and we may assume (up to equivalence of representation) that $\mc{H}(k_1,\rho,\up)$ contains a
vector-valued modular form
\begin{equation}
F_1=\left(\begin{array}{l}q^{\widetilde{\lambda_1}+1}+\cdots\\
q^{\widetilde{\lambda_2}}+\cdots\\
q^{\widetilde{\lambda_3}}+\cdots\\
q^{\widetilde{\lambda_4}}+\cdots\\
q^{\widetilde{\lambda_5}}+\cdots
\end{array}\right),
\end{equation}
whose components are a fundamental system of solutions of $L_\Lambda[f]=0$.

By the usual argument, we obtain a nonzero vector
\begin{equation}
\widetilde{G}=\beta_1E_4D^4F_1+\beta_2E_6D^3F_1+\beta_3E_8D^2F_1+\beta_4E_{10}DF_1+\beta_5E_{12}F_1
\in\mc{H}(k_1+12,\rho,\up)
\end{equation}
such that $G=\frac{\widetilde{G}_1}{\Delta}$ is a nonzero vector in $\mc{H}(k_1,\rho,\up)$, and we have
\begin{lem}\label{lem:15base}
\begin{tabbing}\hspace{3cm}\=\ i) $\mc{H}(k_1,\rho,\up)=\langle F_1,G\rangle$ is 2-dimensional.\\ \\
\>ii) $\mc{H}(k_1+2,\rho,\up)=\langle DF_1,DG\rangle$ is 2-dimensional.
\end{tabbing}
\end{lem}

\pf If there is a relation $\alpha_1F_1+\alpha_2G=0$, then multiplying by $\Delta$ and
substituting with \eq{0} yields
\begin{eqnarray}
0&=&\scriptstyle{\alpha_2\left(\beta_1E_4D^4F_1+\beta_2E_6D^3F_1+\beta_3E_8D^2F_1+\beta_4E_{10}DF_1+
\beta_5E_{12}F_1\right)+\alpha_1\Delta F_1}\nonumber\\
\nonumber\\
&=&\scriptstyle{\alpha_2\beta_1E_4D^4F_1+\alpha_2\beta_2E_6D^3F_1+\alpha_2\beta_3E_8D^2F_1+
\alpha_2\beta_4E_{10}DF_1+[\alpha_1\Delta+\alpha_2\beta_5E_{12}]F_1.}
\end{eqnarray}
By Lemma~\ref{lem:dfind}, each coefficient function in \eq{0} must be zero. In particular, we have
$\alpha_1\Delta=-\alpha_2\beta_5E_{12}$. Comparing $q$-expansions forces
$\alpha_1=0$, which means $\alpha_2G=0$, \ie $\alpha_2=0$. Using Lemma~\ref{lem:5bound} with $N=1$,
we obtain statement $i)$ of the Lemma. Statement $ii)$ then follows immediately from Lemma~\ref{lem:5bound}
and the injectivity of $D$.\qed

Note that this already shows that $\mc{H}(\rho,\up)$ is \emph{not} a cyclic $\mc{R}$-module, since the
lowest weight space is 2-dimensional.

\begin{lem}\label{lem:15plus4}
$\mc{H}(k_1+4,\rho,\up)=\langle D^2F_1,E_4F_1,E_4G\rangle$ is 3-dimensional.
\end{lem}

\pf Again, we assume a relation $\alpha_1D^2F_1+\alpha_2E_4F_1+\alpha_3E_4G=0$, multiply
by $\Delta$, and substitute with \eq{-1}, obtaining
$$\alpha_3E_4(\beta_1E_4D^4F_1+\cdots+\beta_5E_{12}F_1)+\alpha_1\Delta D^2F_1+\alpha_2\Delta E_4F_1=0.$$
Rewriting the above equation in powers of $D$, we again use Lemma~\ref{lem:dfind} and find that the
resulting coefficient functions must all be zero, and in particular we have
$\alpha_1\Delta=-\alpha_3\beta_3E_4E_8$. Comparing $q$-expansions forces $\alpha_1=0$, so our
relation is reduced to $E_4(\alpha_2F_1+\alpha_3G)=0$. By the previous Lemma, we deduce that
$\alpha_2=\alpha_3=0$.\qed

At this point, we would like to prove that the set $\{F_1,G,DF_1,DG,D^2F_1\}$ is independent over
$\m$, and is in fact a basis for $\mc{H}(\rho,\up)$ as $\m$-module. As we shall see in a moment,
this is in fact the case, but to prove this using our more-or-less \emph{ad hoc} methods seems
to be quite difficult, if not impossible. As such, we are forced to appeal, somewhat reluctantly,
to the quite general result \cite[Th.\ 1]{MM}, which we restate here:

\begin{thm}\label{thm:freemmod} Let $\rho:\G\rightarrow\gln{d}$ be a $T$-unitarizable representation of
arbitrary dimension $d$, $\up$ an arbitrary multiplier system, and assume that $\rho(S^2)=\pm I$.
Then $\mc{H}(\rho,\up)$ is a free $\m$-module of rank $d$.\qed
\end{thm}

It must be noted that Theorem 1 of \cite{MM} only addresses the trivial multiplier system setting.
However, one finds upon inspection of the proof that the result holds for arbitrary real weight and
multiplier system; all that is needed is a $d$-dimensional right $\G$-module, arising from the
slash action of $\G$ on $\mc{H}$. Thus we are free to employ this result here.\medskip

We also note that, as is so common in mathematics, the above Theorem is truly an ``existence only''
result, in that the proof of the Theorem gives absolutely no indication of how one might actually
construct an $\m$-basis for a given space $\mc{H}(\rho,\up)$. From this point of view, the present
dissertation is largely an exercise in constructing explicit bases, for various spaces which are
known by Theorem~\ref{thm:freemmod} to be free over $\m$. On the other hand, the classifications
given in this dissertation (together with \cite{M2}) represent the first-known examples of such bases,
and provided the evidence for the initial conjecture which turned into Theorem~\thethm.\medskip

In any event, we are now able to establish

\begin{thm}\label{thm:51mmod}The following hold:

\begin{tabbing}\hspace{1cm}\=\ i) $\mc{H}(\rho,\up)$ is a\=\ free $\m$-module of rank 5, with basis
$\{F_1,G,DF_1,DG,D^2F_1\}$.\\ \\
\>ii) For each $k\geq0$, $\dim\mc{H}(k_1+2k,\rho,\up)=\left\{\begin{array}{lll}\left[\frac{5k}{6}\right]+1&
\ &k\equiv 5\pmod{6}\\ \\
\left[\frac{5k}{6}\right]+2&\ &k\not\equiv 5\pmod{6}\ .\end{array}\right.$\\ \\
\end{tabbing}
\end{thm}

\pf Note that the condition $\rho(S^2)=\pm I$ is satisfied since $\rho$ is irreducible, as pointed
out in Lemma~\ref{lem:2zindec}. So by Theorem~\ref{thm:freemmod}, $\mc{H}(\rho,\up)=\bigoplus_{k\geq0}\mc{H}(k_1+2k,\rho,\up)$ is free of rank 5 over $\m$. From the
first statement of Lemma~\ref{lem:15base}, it is clear that $\mc{H}(k_1,\rho,\up)$ must contribute
two vectors to any $\m$-basis, so we may as well take $F_1$ and $G$. Since there
are no nonzero modular forms of weight two, the second statement of Lemma~\ref{lem:15base} says that
$\mc{H}(k_1+2,\rho,\up)$ must also contribute two vectors to any $\m$-basis, and again we may
as well take $DF_1$ and $DG$; note that because $\mc{H}(\rho,\up)$ is known to be a free $\m$-module,
the set $\{F_1,G,DF_1,DG\}$ \emph{must} be free over $\m$. Similarly, Lemma~\ref{lem:15plus4} tells us that
the fifth basis vector \emph{must} come from $\mc{H}(k_1+4,\rho,\up)$, since clearly
$\mc{H}(k_1+4,\rho,\up)$ is \emph{not} contained in the $\m$-span of $\{F_1,G,DF_1,DG\}$. Thus we may
take $D^2F_1$ as the fifth basis vector, and statement $i)$ of the Theorem is proved. Using statement
$i)$, one may establish statement $ii)$ by noting the relative weights of the basis vectors, and
using the formulae \eqr{dimintwt} for the dimensions of spaces $\m_k$ of integral weight modular forms
of weight $k$. Alternatively, statement $i)$ implies that
\begin{tabbing}\hspace{3cm}\=$\mc{H}(k_1+6,\rho,\up)=\langle E_6F_1,E_6G,E_4DF_1,E_4DG\rangle$,\\ \\
\>$\mc{H}(k_1+8,\rho,\up)=\langle E_8F_1,E_8G,E_6DF_1,E_6DG,E_4D^2F_1\rangle$,\\ \\
\>$\mc{H}(k_1+10,\rho,\up)=\langle E_{10}F_1,E_{10}G,E_8DF_1,E_8DG,E_6D^2F_1\rangle$,
\end{tabbing}
and defining as usual
$$U_{2k}=\left\{\begin{array}{lcl}\mc{H}(k_1+2k,\rho,\up),&\ &k<6\\ \\
\langle E_{2k}F_1,E_{2k}G,E_{2(k-1)}DF_1,E_{2(k-1)}DG,E_{2(k-2)}D^2F_1\rangle,&\ &k\geq6\ ,
\end{array}\right.$$
one employs the standard ``codimension $d$'' argument and observes that statement $i)$ of the
present Theorem implies immediately that for each $k\in\bb{Z}$,
$$\mc{H}(k_1+2k,\rho,\up)=U_{2k}\oplus\Delta\mc{H}(k_1+2(k-6),\rho,\up).$$
Statement $ii)$ of the Theorem then follows by induction, using the above base data.\qed

\begin{cor} The Hilbert-Poincar\'e series \eqr{hp} for $\mc{H}(\rho,\up)$ is
$$\Psi(\rho,\up)(t)=\frac{t^{k_1}(2+2t^2+t^4)}{(1-t^4)(1-t^6)}.$$
\end{cor}\qed

\subsection{N=2}

In this case, Lemma~\ref{lem:5bound} implies that
$$\mc{H}(\rho,\up)=\bigoplus_{k\geq-2}\mc{H}(k_2+2k,\rho,\up),$$
where $k_2=\frac{12(\lambda_0+2)}{5}-4$. Up to equivalence, we may assume that
$\mc{H}(k_2,\rho,\up)$ contains a nonzero vector
\begin{equation}\label{eq:f2}
F_2=\cvec{q^{\widetilde{\lambda_1}+n_1}+\cdots}{\vdots}{q^{\widetilde{\lambda_5}+n_5}+\cdots}
\end{equation}
whose components form a fundamental set of solutions of an MLDE in weight $k_2$, with
indicial roots $\widetilde{\lambda_j}+n_j$, $\sum n_j=2$. By the standard argument involving
linear functionals and Fourier coefficients, there is a nonzero vector
$$G'\in\langle D^4F_2,E_4D^2F_2,E_6DF_2,E_8F_2\rangle\leq\mc{H}(k_2+8,\rho,\up)$$
such that $G=\frac{G'}{\Delta}$ is again holomorphic. Using Lemma~\ref{lem:5bound}, we obtain\medskip

\begin{lem}\label{lem:25base}
\begin{tabbing}\hspace{3cm}\=\ i) $\mc{H}(k_2-4,\rho,\up)=\langle G\rangle$ is
1-dimensional.\\ \\
\>ii) $\mc{H}(k_2-2,\rho,\up)=\langle DG\rangle$ is 1-dimensional.
\end{tabbing}\qed
\end{lem}

If we write $G$ in the form \eqr{fbound} and examine the modular Wronskian $W(G)$, we find as in
the proof of Lemma~\ref{lem:5bound} that the integer $p=\sum p_j$ satisfies
$$0\leq p\leq-\frac{10}{6}+2,$$
\ie $p=0$ and $G$ has the form
\begin{equation}\label{eq:52g}
G=\cvec{q^{\widetilde{\lambda_1}}\sum_{n\geq0}a_1(n)q^n}{\vdots}{q^{\widetilde{\lambda_5}}\sum_{n\geq0}a_5(n)q^n}
\end{equation}
with $a_j(0)\neq0$ for each $j$. Using this information, we will establish\medskip

\begin{lem}\label{lem:523}\begin{tabbing}\hspace{3cm}\=\ i) $\mc{H}(k_2,\rho,\up)=\langle F_2,D^2G,E_4G\rangle$
is 3-dimensional.\\ \\
\>ii) $\mc{H}(k_2+2,\rho,\up)=\langle D^3G,E_4DG,E_6G\rangle$ is 3-dimensional.
\end{tabbing}
\end{lem}

\pf Statement $ii)$ follows directly from Lemmas~\ref{lem:dfind} and \ref{lem:5bound}. Lemma~\ref{lem:5bound}
also implies statement $i)$, provided there is no relation
\begin{equation} F_2=\alpha_1D^2G+\alpha_2E_4G.\end{equation}
Assuming otherwise, we find that $\alpha_1\neq0$ necessarily. Indeed, by definition $F_2$ is of the form
\eqr{f2} with $n_j\geq1$ for at least one $j$, whereas $E_4G$ has the same form as \eqr{52g} (again with
$a_j(0)\neq0$ for each $j$). We have
\begin{eqnarray}\hspace{2cm}F_2&=&\alpha_1D^2G+\alpha_2E_4G\nonumber\\
\ \nonumber\\
&=&\alpha_1D_{k_2-2}D_{k_2-4}G+\alpha_2E_4G\nonumber\\
\ \nonumber\\
&=&\scriptstyle{\alpha_1\left(q\frac{d}{dq}+(k_2-2)P\right)\left(q\frac{d}{dq}+(k_2-4)P\right)G+\alpha_2E_4G}
\end{eqnarray}
and, using \eqr{52g}, \eqr{pexpn}, we find after an elementary computation that the $j^{th}$ component of \eq{0}
has the form
$$\alpha_1a_j(0)Q(\widetilde{\lambda_j})q^{\widetilde{\lambda_j}}+\cdots,$$
where $Q(z)$ denotes the quadratic polynomial
$$Q(z)=z^2+\left(\frac{3-k_2}{6}\right)z+\left(\frac{(k_2-4)(k_2-2)}{144}+\frac{\alpha_2}{\alpha_1}\right).$$
Because the $\widetilde{\lambda_j}$ are distinct, it is clear that, regardless of $\alpha_1,\alpha_2$,
there exist $j_1,j_2\in\{1,2,3,4,5\}$ such that $\widetilde{\lambda_{j_1}}$ and
$\widetilde{\lambda_{j_2}}$ are \emph{not} roots of $Q(z)$. Furthermore, our initial choice of
indicial roots $\widetilde{\lambda_j}+n_j$ was made arbitrarily, except for the conditions $n_j\geq0$
for each $j$, and $\sum n_j=2$.
Thus we are free to assume that $n_{j_1}=n_{j_2}=1$, $n_j=0$ otherwise, and this implies that the left- and
right-hand sides of \eq{0} must differ in the $j_1$ and $j_2$ components. Therefore no such relation
\eq{-3} can exist, and the proof of statement $i)$ is complete.\qed

Using Lemma \thelem, we will prove
\begin{lem}\label{lem:52free}The set $\{F_2,G,DG,D^2G,D^3G\}$ is independent over $\m$.\end{lem}

\pf It is sufficient to prove there is no nontrivial relation
\begin{equation}
M_kF_2=M_{k-2}D^3G+Q_kD^2G+M_{k+2}DG+M_{k+4}G,
\end{equation}
with $M_j,Q_j\in\m_j$. Assuming otherwise, we
divide by $M_k$ and differentiate, obtaining
$$\scriptstyle{DF_2=\frac{M_{k-2}}{M_k}D^4G+\left[\frac{Q_k}{M_k}+D\left(\frac{M_{k-2}}{M_k}\right)\right]D^3G+
\left[\frac{M_{k+2}}{M_k}+D\left(\frac{Q_k}{M_k}\right)\right]D^2G+
\left[\frac{M_{k+4}}{M_k}+D\left(\frac{M_{k+2}}{M_k}\right)\right]DG+
D\left(\frac{M_{k+4}}{M_k}\right)G.}$$
By statement $ii)$ of the previous Lemma, there is a nontrivial relation
\begin{equation}
DF_2=\alpha_1D^3G+\alpha_2E_4DG+\alpha_3E_6G,
\end{equation}
and substituting this into the above equation yields

\begin{eqnarray}
0&=&\frac{M_{k-2}}{M_k}D^4G+\left[\frac{Q_k}{M_k}+D\left(\frac{M_{k-2}}{M_k}\right)-\alpha_1\right]D^3G\nonumber\\
\nonumber\\
&+&\left[\frac{M_{k+2}}{M_k}+D\left(\frac{Q_k}{M_k}\right)\right]D^2G\nonumber\\
\nonumber\\
&+&\left[\frac{M_{k+4}}{M_k}+D\left(\frac{M_{k+2}}{M_k}\right)-\alpha_2E_4\right]DG\nonumber\\
\nonumber\\
&+&\left[D\left(\frac{M_{k+4}}{M_k}\right)-\alpha_3E_6\right]G.
\end{eqnarray}\medskip

\noindent By Lemma~\ref{lem:dfind}, each of the coefficient functions of \eq{0} are identically zero.
In particular, $M_{k-2}=0$, thus $\frac{Q_k}{M_k}=\alpha_1$ is constant, thus
$M_{k+2}=0$, thus $\frac{M_{k+4}}{M_k}=\alpha_2E_4$, thus $D(\alpha_2E_4)=\alpha_3E_6$. But this last
equality, in light of \eq{-1}, says that
$$DF_2=\alpha_1D^3G+\alpha_2E_4DG+D(\alpha_2E_4)G=D(\alpha_1D^2G+\alpha_2E_4G).$$
Since $D$ is injective on $\mc{H}(\rho,\up)$ (Corollary~\ref{cor:dkinj}), we find that $F_2=\alpha_1D^2G+\alpha_2E_4G$, which
is impossible by the previous Lemma.\qed

Applying the bound from Lemma~\ref{lem:5bound}, we obtain the remainder of the base-case data, namely

\begin{tabbing}\hspace{2cm}\=$\mc{H}(k_2+4,\rho,\up)=\langle E_8G,E_6DG,E_4D^2G,E_4F_2\rangle$
is 4-dimensional.\\ \\
\>$\mc{H}(k_2+6,\rho,\up)=\langle E_{10}G,E_8DG,E_6D^2G,E_4D^3G,E_6F_2\rangle$ is 5-dimensional.
\end{tabbing}

\noindent Now we define as usual the spaces $U_{2k}\leq\mc{H}(k_2+2k,\rho,\up)$, which in this
case take the form
$$U_{2k}=\left\{\begin{array}{lcl}\mc{H}(k_2+2k,\rho,\up)&\ &k<4\\ \\
\langle E_{2(k+2)}G, E_{2(k+1)}DG, E_{2k}D^2G, E_{2(k-1)}D^3G, E_{2k}F_2\rangle&\ &k\geq4\ .
\end{array}\right.$$
Using Lemma~\ref{lem:52free} and identical arguments as those employed previously, we find that
for all $k\in\bb{Z}$,
$$\mc{H}(k_2+2k,\rho,\up)=U_{2k}\oplus\Delta\mc{H}(k_2+2(k-6),\rho,\up),$$
and this establishes
\begin{thm}The following hold:
\begin{tabbing}\hspace{1cm}\=\ i) \=$\mc{H}(\rho,\up)=\bigoplus_{k\geq-2}\mc{H}(k_2+2k,\rho,\up)$
is a free $\m$-module of rank 5,\\
\>\>with basis $\{G,DG,D^2G,D^3G,F_2\}$.\\ \\
\>ii) For each $k\geq-2$, $\dim\mc{H}(k_2+2k,\rho,\up)=\left\{\begin{array}{lll}\left[\frac{5k}{6}\right]+2&
\ &k\equiv 5\pmod{6}\\ \\
\left[\frac{5k}{6}\right]+3&\ &k\not\equiv 5\pmod{6}\ .\end{array}\right.$
\end{tabbing}\qed
\end{thm}

\begin{rem}\label{rem:52cyclic} As in the $N=1$ case, we again observe that $\mc{H}(\rho,\up)$ is \emph{not} cyclic
as $\mc{R}$-module. Note, however, that the two modules fail to be cyclic for different reasons:
in the $N=1$ case, the lowest weight space was 2-dimensional (Lemma~\ref{lem:15base}), which
is an immediate disqualification for cyclicity, whereas the present module \emph{does} have
a lowest weight space that is 1-dimensional (Lemma~\ref{lem:25base}). But unlike the $N=0$ case
(Theorem~\ref{thm:05}), where the module was cyclic, in the present case the lowest
weight vector $G$ does \emph{not} have components which form a fundamental system of solutions of
an MLDE; this is already implicit in the initial discussion which yielded the five cases we are
elucidating, or one may compute directly and find that the Wronskian of $G$ has the form
$c\eta^{24\lambda_0}E_4$, $c\neq0$, and then compare with statement $ii)$ of Theorem~\ref{thm:wmlde}.
By \cite[Th.\ 3]{MM} (at least in the trivial multiplier system case), this means
$\mc{H}(\rho,\up)$ cannot be cyclic as $\mc{R}$-module.\qed
\end{rem}

\begin{cor}The Hilbert-Poincar\'e series \eqr{hp} of $\mc{H}(\rho,\up)$ is
$$\Psi(\rho,\up)(t)=\frac{t^{k_2-4}(1+t^2+2t^4+t^6)}{(1-t^4)(1-t^6)}.$$
\end{cor}\qed

\subsection{N=3}

In this case, we have by Lemma~\ref{lem:5bound} that

$$\mc{H}(\rho,\up)=\bigoplus_{k\geq-3}\mc{H}(k_3+2k,\rho,\up),$$
where $k_3=\frac{12}{5}(\lambda_0+3)-4.$ We shall be more particular in this case about our
choice of indicial roots. First, we single out a $j_1\in\{0,1,\cdots,5\}$ such that
$\widetilde{\lambda_{j_1}}\neq k_3-6$; this is certainly possible, since the $r_1,\cdots,r_5$ are
distinct. Having made this selection, we then choose distinct $j_2,j_3,j_4,j_5\in\{0,1,\cdots,5\}-\{j_1\}$,
and set $n_{j_1}=n_{j_2}=0$, $n_{j_3}=n_{j_4}=n_{j_5}=1$. Using these integers, we construct the
indicial roots $\widetilde{\lambda_{j_1}}+n_{j_1},\cdots,\widetilde{\lambda_{j_5}}+n_{j_5}$, and
obtain via Corollary~\ref{cor:umlde} a unique $5^{th}$ order MLDE in weight $k_3$ with the
above indicial roots. We may assume, up to equivalence of representation, that there is a
vector-valued modular form
\begin{equation}
F_3=\cvec{q^{\widetilde{\lambda_1}+n_1}+\cdots}{\vdots}{q^{\widetilde{\lambda_5}+n_5}+\cdots}
\end{equation}
in $\mc{H}(k_3,\rho,\up)$, whose components form a fundamental set of solutions of the
MLDE just constructed.\medskip

Set $V=\langle D^3F_3,E_4DF_3,E_6F_3\rangle\leq\mc{H}(k_3+6,\rho,\up)$. From \eq{0}, we see that
each nonzero vector in $V$ has the form
\begin{equation}
F=\cvec{q^{\widetilde{\lambda_1}+n_1}\sum_{n\geq0}a_1(n)q^n}{\vdots}
{q^{\widetilde{\lambda_5}+n_5}\sum_{n\geq0}a_5(n)q^n}.
\end{equation}
Define the linear functionals $\phi_1,\phi_2:V\rightarrow\bb{C}$, which take $F$ as in \eq{0}
to $\phi_1(F)=a_{j_1}(0)$, $\phi_2(F)=a_{j_2}(0)$. Then $\ker\phi_1\cap\ker\phi_2\neq\{0\}$,
so there is a nonzero $G_1\in\mc{H}(k_3-6,\rho,\up)$ such that $\Delta G_1\in V$. From the
bound in Lemma~\ref{lem:5bound} for $N=3$, we find that
\begin{lem} $\mc{H}(k_3-6,\rho,\up)=\langle G_1\rangle$ is 1-dimensional.
\end{lem}\qed

Similarly, we define the four-dimensional subspace
$$U=\langle D^4F_3,E_4D^2F_3,E_6DF_3,E_8F_3\rangle\leq\mc{M}(k_3+8,\rho).$$
Again, each $F\in U$ has the form \eq{0}, and we define the linear functionals
$\phi_1,\phi_2,\phi_3:U\rightarrow\bb{C}$, where $\phi_1(F)=a_{j_1}(0)$, $\phi_2(F)=a_{j_1}(1)$,
and $\phi_3(F)=a_{j_2}(0)$. Again we have $\cap_{i=1}^3\ker\phi_i\neq\{0\}$, and we obtain
a nonzero $G_2\in\mc{H}(k_3-4,\rho,\up)$ such that $\Delta G_2\in U$. Note that, by construction,
$G_2$ has the form
\begin{equation}
G_2=\cvec{q^{\widetilde{\lambda_1}+n_1}\sum_{n\geq0}b_1(n)q^n}{\vdots}
{q^{\widetilde{\lambda_5}+n_5}\sum_{n\geq0}b_5(n)q^n},
\end{equation}
with $b_{j_1}(0)=0$. We will need this fact in proving
\begin{lem}\label{lem:53plus2}
\begin{tabbing}\hspace{3cm}\=\ i) $\mc{H}(k_3-4,\rho,\up)=\langle DG_1,G_2\rangle$ is 2-dimensional.\\ \\
\>ii) $\mc{H}(k_3-2,\rho,\up)=\langle D^2G_1,E_4G_1\rangle$ is 2-dimensional.
\end{tabbing}
\end{lem}

\pf Statement $ii)$ follows immediately from Lemmas~\ref{lem:dfind} and \ref{lem:5bound}.
Lemma~\ref{lem:5bound} also implies statement $i)$, provided there is no $\alpha\in\bb{C}$ such
that $DG_1=\alpha G_2$. The proof of this fact follows exactly as in the proof of statement $ii)$ of
Lemma~\ref{lem:k1base}: an examination of its modular Wronskian shows that $G_1$ must
have the form \eq{-1}, with $a_j(0)\neq0$ for each $j$, and a direct computation of $DG_1$,
along with the information $b_{j_1}(0)=0$, $r_{j_1}\neq k_3-6$, completes the proof.\qed

\begin{lem}\label{lem:53free} The set $\{G_1,DG_1,D^2G_1,D^3G_1,G_2\}$ is independent over $\m$.\end{lem}

\pf We need only show there is no relation
$$M_kG_2=M_{k+2}G_1+Q_kDG_1+M_{k-2}D^2G_1+M_{k-4}D^3G_1,$$
with $M_k\neq0$, $M_i,Q_i\in\mc{M}(i,\textbf{1})$ for each $i$. Assuming otherwise, we divide
and differentiate, obtaining
\begin{equation}
\scriptstyle{DG_2=D\left(\frac{M_{k+2}}{M_k}\right)G_1+\left[\frac{M_{k+2}}{M_k}+D\left(\frac{Q_k}{M_k}\right)\right]DG_1+
\left[\frac{Q_{k}}{M_k}+D\left(\frac{M_{k-2}}{M_k}\right)\right]D^2G_1+
\left[\frac{M_{k-2}}{M_k}+D\left(\frac{M_{k-4}}{M_k}\right)\right]D^3G_1+
\frac{M_{k-4}}{M_k}D^4G_1.}
\end{equation}
By Lemma~\ref{lem:53plus2} Corollary~\ref{cor:dkinj}, there is a nontrivial relation
$$DG_2=\alpha_1D^2G_1+\alpha_2E_4G_1,$$
with $\alpha_2\neq0$. Substituting into \eq{0} and using Lemma~\ref{lem:dfind}, we find that
$M_{k-4}=0$, thus $M_{k-2}=0$, thus $\frac{Q_k}{M_k}=\alpha_1$, thus $M_{k+2}=0$, thus $\alpha_2=0$, a contradiction.\qed

Using Lemmas~\ref{lem:5bound}, \ref{lem:53plus2} and \thelem, we find that
\begin{tabbing}\hspace{1cm}\=$\mc{H}(k_3,\rho,\up)=\langle D^3G_1,E_4DG_1,E_6G_1,E_4G_2\rangle$
is 4-dimensional.\\ \\
\>$\mc{H}(k_3+2,\rho,\up)=\langle E_4D^2G_1,E_6DG_1,E_8G_1,E_6G_2\rangle$ is 4-dimensional.\\ \\
\>$\mc{H}(k_2+4,\rho,\up)=\langle E_4D^3G_1,E_6D^2G_1,E_8DG_1,E_{10}G_1,E_8G_2\rangle$ is 5-dimensional.
\end{tabbing}
As usual, we define
$$U_{2k}=\left\{\begin{array}{lcl}\mc{H}(k_3+2k,\rho,\up),&\ &k<3\\ \\
\langle E_{2k}D^3G_1,E_{2(k+1)}D^2G_1,E_{2(k+2)}DG_1,E_{2(k+3)}G_1,E_{2(k+2)}G_2\rangle,&\ &k\geq3\ ,
\end{array}\right.$$
and employ Lemma~\ref{lem:53free} and the standard ``codimension $d$'' argument to find that
$$\mc{H}(k_3+2k,\rho,\up)=U_{2k}\oplus\Delta\mc{H}(k_3+2(k-6),\rho,\up)$$
for all $k\in\bb{Z}$. This establishes
\begin{thm}The following hold:
\begin{tabbing}\hspace{1cm}\=\ i) \=$\mc{H}(\rho,\up)=\bigoplus_{k\geq-3}\mc{M}(k_3+2k,\rho)$
is free of rank 5 as $\m$-module\\
\>\>with basis $\{G_1,DG_1,G_2,D^2G_1,D^3G_1\}$.\\ \\
\>ii) For each $k\geq-3$, $\dim\mc{H}(k_3+2k,\rho,\up)=
\left\{\begin{array}{lll}\left[\frac{5k}{6}\right]+3&\ &k\equiv 5\pmod{6}\\ \\
\left[\frac{5k}{6}\right]+4&\ &k\not\equiv 5\pmod{6}\ .\end{array}\right.$
\end{tabbing}
\end{thm}\qed

\begin{cor} The Hilbert-Poincar\'e series \eqr{hp} of $\mc{H}(\rho,\up)$ is
$$\Psi(\rho,\up)(t)=\frac{t^{k_3-6}(1+2t^2+t^4+t^6)}{(1-t^4)(1-t^6)}.$$
\qed
\end{cor}

\subsection{N=4}

In this case, we have by Lemma~\ref{lem:5bound} that
$$\mc{H}(\rho,\up)=\bigoplus_{k\geq-4}\mc{H}(k_4+2k,\rho,\up),$$
where $k_4=\frac{12}{5}(\lambda_0+4)-4$. Choose an arbitrary $i\in\{0,1,\cdots,5\}$,
define the integers $n_1,\cdots,n_5$ as $n_i=0$, $n_j=1$ otherwise, and set $$\Lambda=\{\widetilde{\lambda_1}+n_1,\cdots,\widetilde{\lambda_5}+n_5\}.$$
By Corollary~\ref{cor:umlde}, there is a unique MLDE $L_\Lambda[f]=0$ in weight $k_4$, whose set of
indicial roots is $\Lambda$. Up to isomorphism, we may assume that $\mc{H}(k_4,\rho,\up)$ contains a nonzero vector
$$F_4=\cvec{q^{\widetilde{\lambda_1}+n_1}+\cdots}{\vdots}{q^{\widetilde{\lambda_5}+n_5}+\cdots}$$
whose components are a fundamental system of solutions of $L_\Lambda[f]=0$.\medskip

Define the 2-dimensional subspace $U=\langle D^2F_4,E_4F_4\rangle\leq\mc{H}(k_4+4,\rho,\up)$. Then
every vector in $U$ has the form
\begin{equation}\label{eq:54f}
F=\cvec{q^{\widetilde{\lambda_1}+1-\delta_{1,i}}\sum_{n\geq0}a_1(n)q^n}{\vdots}
{q^{\widetilde{\lambda_5}+1-\delta_{5,i}\sum_{n\geq0}}a_5(n)q^n},
\end{equation}
and the linear functional $\phi_i:U\rightarrow\bb{C}$, $F\mapsto\phi_i(F):=a_i(0)$ has
a nontrivial kernel. This means there is a nonzero vector $G_1\in\mc{H}(k_4-8,\rho,\up)$ such
that $\Delta G_1\in U$. Assuming $G_1$ has the form \eqr{fbound} and evaluating the inequality
\eqr{wbound} at $k=-4$, $N=4$, we find that
$$0\leq p\leq \left[-\frac{20}{6}\right]+4,$$
so $p=0$ and $G_1$ has the form
\begin{equation}\label{eq:54g1}
G_1=\cvec{q^{\widetilde{\lambda_1}}\sum_{n\geq0}a_1(n)q^n}{\vdots}
{q^{\widetilde{\lambda_5}}\sum_{n\geq0}a_5(n)q^n},
\end{equation}
with $a_j(0)\neq0$ for each $j$. Using Lemma~\ref{lem:5bound}, we obtain
\begin{lem}\label{lem:54base}$\mc{H}(k_4-8,\rho,\up)=\langle G_1\rangle$ is 1-dimensional.
\qed
\end{lem}

We continue on, defining the 3-dimensional subspace $V=\langle D^3F_4,E_4DF_4,E_6F_4\rangle
\leq\mc{H}(k_3+6,\rho,\up)$. Again we have that every nonzero $F\in V$ is as in \eqr{54f}.
From \eqr{54g1} and the definition of $D$, one sees that the $j^{th}$ component of
$DG_1=D_{k_4-8}G_1$ has the form
$$a_j(0)\left(\widetilde{\lambda_j}-\frac{k_4-8}{12}\right)q^{\widetilde{\lambda_j}}+\cdots.$$
Since $\frac{k_4-8}{12}=\frac{\lambda_0+4}{5}-1$ is fixed, there must be an $i_1\in\{0,1,\cdots,5\}$
such that $\widetilde{\lambda_{i_1}}\neq\frac{k_4-8}{12}$. Define the linear functionals $\phi_0,\phi_1:V\rightarrow\bb{C}$, where $F$ as in \eqr{54f} is
sent to $\phi_0(F):=a_i(0)$, $\phi_1(F):=a_{i_1}(1)$ respectively, then it is clear that a
nonzero vector  $G_2\in\mc{H}(k_4-6,\rho,\up)$ can be found such that $\Delta G_2\in V$, and
such that the $i_1$ component of $G_2$ has the form
$$q^{\widetilde{\lambda_{i_1}}+1}\sum_{n\geq0}b_{i_1}(n)q^n.$$
Again assuming that $G_2$ has the form \eqr{fbound} and using \eqr{wbound} with $k=-3$, $N=4$,
we find that
$$1\leq p\leq\left[-\frac{15}{6}\right]+4,$$
so $p=1$ and $G_2$ looks like
\begin{equation}\label{eq:54g2}
G_2=\cvec{q^{\widetilde{\lambda_1}+\delta_{1,i_1}}\sum_{n\geq0}b_1(n)q^n}{\vdots}
{q^{\widetilde{\lambda_5}+\delta_{5,i_1}}\sum_{n\geq0}b_5(n)q^n},
\end{equation}
with $b_j(0)\neq0$ for each $j$.
In particular, one sees that, by construction, $DG_1$ and $G_2$ are linearly independent, since
their respective $i_1$-components have $q$-expansions which differ in leading exponent.
By Lemma~\ref{lem:5bound}, we have
\begin{lem}\label{lem:54plus2}
$\mc{H}(k_4-6,\rho,\up)=\langle DG_1,G_2\rangle$ is 2-dimensional.\qed
\end{lem}

\begin{rem}Lemma~\ref{lem:54plus2} shows that, once again, we are in a case where $\mc{H}(\rho,\up)$ has a lowest
weight space which is 1-dimensional, yet fails nonetheless to be cyclic as $\mc{R}$-module (\cf\ Remark~\ref{rem:52cyclic}).\qed
\end{rem}

If we compute $DG_2=D_{k_4-6}G_2$, then we find, as above, that the $j^{th}$ component has the form
$$b_j(0)\left(\widetilde{\lambda_j}+\delta_{j,i_1}-\frac{k_4-6}{12}\right)
q^{\widetilde{\lambda_j}+\delta_{j,i_1}}+\cdots.$$
Fix an $i_2\in\{1,\cdots,5\}-\{i_1\}$ such that $\widetilde{\lambda_j}\neq\frac{k_4-6}{12}$,
and set
$$W=\langle D^4F_4,E_4D^2F_4,E_6DF_4,E_8F_4\rangle\leq\mc{H}(k_4+8,\rho,\up).$$
This time we define the linear functionals $\phi_1,\phi_2,\phi_3:W\rightarrow\bb{C}$, where $F$ in \eqr{54f}
is mapped to $\phi_1(F):=a_i(0)$, $\phi_2(F):=a_{i_1}(1)$, $\phi_3(F):=a_{i_2}(1)$, respectively.
Again, the intersection $\cap_{j=1}^3\ker\phi_j$ is nontrivial, and this means there is a nonzero
$G_3\in\mc{H}(k_4-4,\rho,\up)$ such that $\Delta G_3\in W$. Once again we assume $G_3$ has the form
\eqr{fbound} and examine \eqr{wbound} with $k=-2$, $N=4$, obtaining the inequality
$$0\leq p\leq\left[-\frac{10}{6}\right]+4,$$
so in fact $p=2$ and $G_3$ has the form
\begin{equation}\label{eq:54g3}
G_3=\cvec{q^{\widetilde{\lambda_1}+p_1}\sum_{n\geq0}c_1(n)q^n}{\vdots}
{q^{\widetilde{\lambda_5}+p_5}\sum_{n\geq0}c_5(n)q^n},
\end{equation}
where $p_{i_1}=p_{i_2}=1$, $p_j=0$ otherwise, and $c_j(0)\neq0$ for each $j$. We can now prove
\begin{lem}\label{lem:54plus4} $\mc{H}(k_4-4,\rho,\up)=\langle E_4G_1,DG_2,G_3\rangle$ is 3-dimensional.
\end{lem}

\pf By Lemma~\ref{lem:5bound}, it suffices to prove there is no nontrivial relation
\begin{equation}
\alpha_1E_4G_1+\alpha_2DG_2+\alpha_3G_3=0.
\end{equation}
Computing using \eqr{54g1}, \eqr{54g2}, \eqr{54g3}, we find that the $j^{th}$ component of \eq{0}
has the form
\begin{equation}
\left(\alpha_1a_j(0)q^{\widetilde{\lambda_j}}+\cdots\right)+
\left(\alpha_2b_j(0)\left(\widetilde{\lambda_j}+\delta_{j,i_1}-\frac{k_4-6}{12}\right)
q^{\widetilde{\lambda_j}+\delta_{j,i_1}}+\cdots\right)+
\left(\alpha_3c_j(0)q^{\widetilde{\lambda_j}+p_j}+\cdots\right).
\end{equation}
Note that when $j=i_1$ in \eq{0}, we find that $\alpha_1=0$, since $a_{i_1}(0)\neq0$ and the second and
third expansions in \eq{0} have no $q^{\widetilde{\lambda_{i_1}}}$ term. Similarly, the $j=i_2$
version of \eq{0} shows that $\alpha_2=0$, since we now know that $\alpha_1=0$, the second expansion
in \eq{0} has leading term $q^{\widetilde{\lambda_{i_2}}}$, and the third has leading term
$q^{\widetilde{\lambda_{i_2}}+1}$. But then $\alpha_3=0$ also, clearly, and we have the result.\qed

We are now forced to proceed as in the $N=1$ case, and invoke Theorem~\ref{thm:freemmod} in order to prove
\begin{thm}The following hold:
\begin{tabbing}\hspace{1cm}\=\ i) \=$\mc{H}(\rho,\up)=\bigoplus_{k\geq-4}\mc{H}(k_4+2k,\rho,\up)$
is free of rank 5 as $\m$-module,\\
\>\>with basis $\{G_1,DG_1,G_2,DG_2,G_3\}$.\\ \\
\>ii) For all $k\geq-4$, $\dim\mc{H}(k_4+2k,\rho,\up)=\left\{\begin{array}{lcl}\left[\frac{5k}{6}\right]+4,&
\ &k\equiv 5\pmod{6}\\ \\
\left[\frac{5k}{6}\right]+5&\ &k\not\equiv 5\pmod{6}\ .\end{array}\right.$\\ \\
\end{tabbing}
\end{thm}

\pf By Theorem~\ref{thm:freemmod}, we know that $\mc{H}(\rho,\up)$ is free of rank 5 over $\m$. It
is clear from Lemmas~\ref{lem:54base} and \ref{lem:54plus2} that we may take $\{G_1,DG_1,G_2\}$ as
part of any basis. But the $\m$-span of $\{G_1,DG_1,G_2\}$ intersects $\mc{H}(k_4-4,\rho,\up)$ only
in $\langle E_4G_1\rangle$, so by Lemma~\ref{lem:54plus4}, there \emph{must} be two more vectors
available in $\mc{H}(k_4-4,\rho,\up)$ to match up with the three we already have. We may as well take
$DG_2$ and $G_3$, and this establishes statement $i)$ of the Theorem.

Using statement $i)$ and Lemma~\ref{lem:5bound}, we immediately obtain the remaining base data
\begin{tabbing}\hspace{3cm}\=$\mc{H}(k_4-2,\rho,\up)=\langle E_6G_1,E_4DG_1,E_4G_2\rangle$,\\ \\
\>$\mc{H}(k_4,\rho,\up)=\langle E_8G_1,E_6DG_1,E_6G_2,E_4DG_2,E_4G_3\rangle$,\\ \\
\>$\mc{H}(k_4+2,\rho,\up)=\langle E_{10}G_1,E_8DG_1,E_8G_2,E_6DG_2,E_6G_3\rangle$,\\ \\
\end{tabbing}
and statement $ii)$ of the Theorem follows in either of the ways described in the proof of the
analogous statement of Theorem~\ref{thm:51mmod}, \ie by using \eqr{dimintwt} or by defining
the spaces
$$U_{2k}=\langle E_{2k}G_1,E_{2(k-1)}DG_1,E_{2(k-1)}G_2,E_{2(k-2)}DG_2,E_{2(k-2)}G_3\rangle$$
and noting that
$$\mc{H}(k_4+2k,\rho,\up)=U_{2k}\oplus\Delta\mc{H}(k_4+2(k-6),\rho,\up).$$
\qed

\begin{cor} The Hilbert-Poincar\'e series \eqr{hp} of $\mc{H}(\rho,\up)$ is
$$\Psi(\rho,\up)(t)=\frac{t^{k_4-8}(1+2t^2+2t^4)}{(1-t^4)(1-t^6)}.$$
\qed
\end{cor}

\bibliographystyle{amsplain}
\addcontentsline{toc}{chapter}{Bibliography}
\bibliography{Marks_dissertation}
\end{document}